\documentstyle[12pt,amssymbols]{article}
\newcommand{\bq}{\begin{equation}} 
\newcommand{\eq}{\end{equation}} 
\newcommand{\bqa}{\begin{eqnarray}} 
\newcommand{\eqa}{\end{eqnarray}} 
\newcommand{\ra}{\rightarrow} 
\def\curlra{\buildrel{\sim}\over\longrightarrow}
\def\qbinom#1#2{{#1}\atopwithdelims[]{#2}}

\def\half{{1 \over 2}}

\def\g{\gamma}

\def\Z{{\Bbb Z}}
\def\C{{\Bbb C}}
\def\N{{\Bbb N}}

\def\ep{\epsilon}
\def\l{\lambda} 
 
\def\ov{\over}

\def\be{\begin{equation}}
\def\ee{\end{equation}}
\def\br{\begin{array}}
\def\er{\end{array}}
\def\bea{\begin{eqnarray}}
\def\eea{\end{eqnarray}}
\def\ba{\begin{equation}\begin{array}}
\def\ea{\end{eqnarray}}
\def\bac{\begin{eqnarray}}
\def\nn{\nonumber}
\def\ra{\rightarrow} 
\def\al{\alpha} 
\def\2pi{1\over 2\pi i} 
\def\q{q-q^{-1}} 

\def\~{\tilde} 
\def\newline{\hfil\break}

\def\ra{\rightarrow} 
\def\va{\varphi} 
 
\def\sq2{\sqrt{2}} 
\def\sqk2{\sqrt{2(k+2}} 
\def\sqk{\sqrt{k}} 
\def\sqs{\sqrt{2\over k}}

\def\be{\begin{equation}} 
\def\ee{\end{equation}} 
\def\br{\begin{array}} 
\def\er{\end{array}} 
\def\bea{\begin{eqnarray}} 
\def\eea{\end{eqnarray}} 
\def\nn{\nonumber}

\newcommand{\uq}{U_q (\widehat{sl(2)})}

\def\vt{\tilde{\Phi}}

\newcommand{\LS}{\widehat{sl(2)}}

\def\ep{\epsilon}
\def\epp{\epsilon^{\prime}}

\begin{document} 
\bibliographystyle{unsrt}
\begin{titlepage}
\rightline{ITP-29} 
\rightline{\today} 
\vbox{\vspace{15mm}} 
\vspace{1.0truecm} 
\begin{center} 
{\LARGE \bf Bosonization of quantum affine groups and its
application to  
the higher spin Heisenberg model  
}\\[8mm] 
{\large A.H. BOUGOURZI}\\ 
[3mm]
{\it Institute for Theoretical Physics\\ 
SUNY at Stony Brook, Stony Brook, NY 11794}\\[15mm] 
\end{center} 
\begin{abstract} 
In this paper, we present a detailed analysis  of the  
diagonalization of the higher spin Heisenberg model 
using its quantum affine symmetry $U_q(\widehat{sl(2)})$.  
In particular, we describe the bosonizations of the latter
algebra, its highest weight representations, vertex operators 
and screening operators.  Finally, we use this bosonization
method to 
compute the vacuum-to-vacuum expectation 
values and the form factors of any  local operator.

\end{abstract} 
\footnotetext[1]{ 
Email: {\tt bougourz@insti.physics.sunysb.edu} 
}

\end{titlepage} 

\section{Introduction}

Ever since the remarkable work of Bethe \cite{Bet31} 
on the Heisenberg model \cite{Hei28},
a tremendous amount of work has been accomplished towards 
its exact solution \cite{Hul38}-\cite{Fleal95}.\footnote{This list is by
no means exhaustive and  many other interesting 
references are available in the literature.} The  reason is
that  this  model exhibits
strong quantum effects due to many-body  interactions that
can be both experimentally measured and exactly calculable. 
In particular, its 
behavior in the antiferromagnetic regime 
is highly nontrivial because each  eigenstate of the
 Heisenberg Hamiltonian is built from 
a superposition of many  spin configurations. 
It is even more interesting in
the thermodynamic limit where its Hilbert space acquires 
a particle-like structure. 
As an example of this quantum behavior, 
the spin of the elementary excitations, called spinons or 
quasi particles, 
is always equal to $s=1/2$ 
regardless of the value of the local spin \cite{JoMc72,FaTa79}. 
Moreover, spinons are always created in pairs by local
operators. While the latter statement can be intuitively  
expected from the former, 
only rigorous  
analysis based on the Bethe ansatz and later confirmed
through the crystal base theory led to the former. 
That spinons are always created  in pairs is due to the
fact that 
the total spin  can only change 
by discrete integer values, which are 
then  shared among 
an even number of spinons.

The Bethe ansatz method has played a major role 
in deriving the bulk of the existing results related to
this model, namely the structure of its Hilbert space, 
its spectrum, S matrix, long distance correlation functions
etc...
To confirm and complete these results another  
method was found in Refs. \cite{JiMi94,Idzal93}. 
Its main new
feature is that it makes use of the quantum group symmetry
of this model from the start. This symmetry 
 appears spontaneously  in the 
thermodynamic limit and in the antiferromagnetic regime. The
latter two cases
are obviously of great interest because the other ones, 
namely ferromagnetic and massless regimes, 
are almost fully solved  
through the Bethe ansatz or conformal field theory approaches.

The above  quantities have been extensively studied
in the framework of the quantum group 
symmetry method for the spin-1/2
case. Its main conceptual  outcome  is 
that it allows for 
the identification of 
the dynamical fields responsible for the
creation of spinons from a vacuum state. 
They have the form of vertex
operators similar to those frequently used in the context
of string theory and conformal field theory. 
In fact, most of the detailed calculation
of correlation functions relies on technics already familiar
in the context of the latter theories.
More precisely, they are analogous to those used 
for the calculation of correlation functions on a torus. 
While, it has
benefited considerably from string theory and conformal field
theory, this method  has also paid back 
in the sense that 
a simpler  expression for the trace of
vertex operators (i.e., correlation functions on a torus) has
been found in the context of the quantum group symmetry 
\cite{Jin92,Idz93}. 
Once the above identification has been made the
whole problem of evaluating  form factors and 
correlation functions reduces to a technical problem
of calculating traces of vertex operators over irreducible
highest weight representations (IHWR's) of the quantum
affine algebra $\uq$. 
Unfortunately, 
these traces
are nontrivial because the latter IHWR's
are infinite-dimensional 
and the vertex operators have
complicated (generalized) commutation relations. 
This 
problem is solved through 
the bosonization method  where the latter 
elements, i.e., $\uq$, its IHWR's, and the vertex
operators,  are all realized in terms 
of harmonic oscillators satisfying simple Heisenberg algebras.

Let us note here that in the spin-1/2 case 
there is an appealing identity 
that
made the use of the bosonization method straightforward.
It is the fact that
each  Fock space of the harmonic oscillators coincides precisely
with one of the IHWR's. Since evaluation of 
traces over Fock spaces
is simple (especially due the trace formula mentioned 
before) the short distance 
form factors and correlation functions have been
easily evaluated in Ref. \cite{JiMi94}.
However, in the higher spin case, i.e., $s\geq 1$ such a 
coincidence does  no longer
hold. While the bosonization technic is still applicable
the Fock spaces are much larger than the IHWR's. 
Since  both of them 
are infinite-dimensional spaces the projection
from the former to the latter is subtle. It requires a 
detailed and relatively complicated 
cohomology analysis of these Fock spaces 
\cite{BeFe90,Kon94}. 
This accounts for the main technical difficulty in the
evaluation of   
the form factors and correlation functions in this case.
Let us note  that 
 the difficulty in the spin-1 case is somewhat intermediary. 
Here, the IHWR's coincide with
the Fock spaces built from a set of harmonic oscillators 
and another set of modes satisfying the
Clifford algebra. This has also allowed for the calculation
of correlation functions \cite{BoWe94}. 

Since all of the spin-1/2 results
obtained through this symmetry method were completed and 
compiled in 
Ref. \cite{JiMi94} we will not touch upon this case here. 
We will
instead focus on the remaining goal of  
evaluating  the correlation functions
and form factors  for $s> 1$. For completeness we will
review and compile most of the developments leading up to this
goal.  
Obviously, we expect the final results to be complicated. 
However, they are not necessarily more complicated 
than their counterparts
in string theory and conformal field theory on genus 1 and higher 
surfaces.
Also, they are still accessible through symbolic or numeric
manipulations. 
Although, a finite chain  model can also be 
solved directly
through these manipulations, it will not  yield the subtle
features that arise just in 
the thermodynamic limit.
Furthermore,
the solution through the  symmetry method allows
for the identification of the contributions to form factors
 from each sector
specified by the number of spinon pairs. As far as we know
this is not yet possible through a direct numerical calculation
which only leads to the total contribution. Also, one can
always derive simple perturbative results from the exact
ones. The perturbation parameter is the anisotropy 
parameter $q$ which coincides with the deformation parameter in
$\uq$.  
Finally, despite
the fact that the bulk of the final results is complicated 
one can still derive  simple selection rules 
and identities among  correlation functions and form
factors.
For example, we will see  that
the fact that spinons arise just in pairs is a consequence
of a selection rule. 
Also from a fundamental
point of view this method allows us to check that 
Smirnov's axioms \cite{Smi92}, which were first intended 
for integrable 
models of quantum field theory, are still satisfied (with
a little modification) by the form factors of the present
lattice model. This fact, might lay the foundation of 
  yet another
method to the evaluation of  form factors using
the modified Smirnov's axioms.
Form factors are also known to satisfy the  q-KZ equation
which plays a crucial role in many fields of physics and 
mathematics. Therefore, the  
expressions that we derive here
 are  explicit solutions to this equation.

This paper is organized as follows. In section one we review 
the quantum affine algebra $\uq$. In section two 
we review its explicit
bosonization, that is, its realization in terms of bosonic
fields. We also introduce its associated screening currents
and their bosonization. In section three we recall the
diagonalization of the higher spin Heisenberg model through 
$\uq$ and its vertex
operators. The main consequence of this diagonalization 
is that 
the correlation functions and form factors are 
expressed as traces of these vertex operators.
For the purpose of evaluating these traces we review 
in section four the bosonizations of
the vertex operators and the $\uq$-IHWM's. In section five
we use this bosonization to evaluate explicitly the 
general N-point correlation functions and N-point form factors
of local operators of this model. This is the main 
new result of this paper. In doing so we derive
certain selection rules and an identity for form factors.
We also make a consistency check using the
known results of the $s=1/2$ case. For completeness,
we compile some other results concerning these quantities.
In particular, we complete the proof that the form
factors satisfy all (modified) Smirnov's axioms. We also review
the fact that they satisfy the q-KZ equation. Finally, 
we present in the appendix some useful relations, 
the operator product expansions and 
the trace formula of the 
operators used in this paper.   
\newpage
\section{The  quantum affine algebra $\uq$ and its
representation theory}

Since the key idea here is to solve the  higher spin 
Heisenberg model
through its quantum affine symmetry $\uq$,
let us briefly review this algebra and its representations. 

\subsection{\bf  $\uq$}

This  algebra  is  generated 
by the  Chevalley elements 
$\{e_i, f_i,  k_i^{\pm 1}, q^{\pm d};\> i=0,1\}$, 
with defining relations  \cite{Dri85,Jim85}
\bac
&&k_ik_j=k_jk_i,\quad k_ik_i^{-1}=k_i^{-1}k_i=1,\nn \\  
&&k_i e_i k^{-1}_i=q^2e_i,\quad k_i e_j k^{-1}_i=
q^{-2}e_j, ~~i\neq j,\nn \\  
&&k_i f_i k^{-1}_i=q^{-2}f_i,\quad k_i f_j k^{-1}_i=
q^2f_j, ~~i\neq j,\nn \\  
&&[e_i, f_j]=\delta_{i,j}{k_i-k_i^{-1}\over q-q^{-1}},\nn \\ 
&&(e_{i})^3e_{j}-[3]
(e_{i})^2e_{j}e_{i}+
[3]e_{i}e_{j}(e_{i})^2-e_{i}(e_{j})^3=0,\nn \\ 
&&(f_{i})^3f_{j}-[3]
(f_{i})^2f_{j}f_{i}+
[3]f_{i}f_{j}(f_{i})^2-f_{i}(f_{j})^3=0,\nn \\ 
&&q^de_iq^{-d}=q^{\delta_{i,0}}e_i,
\quad q^df_iq^{-d}=q^{-\delta_{i,0}}f_i\quad 
q^dk_iq^{-d}=k_i, 
\label{cheval} 
\ea
where $[x]=(q^{x}-q^{-x})/(\q)$ is a $q$-integer and
$q$ is  a deformation parameter.
The main property of this algebra is 
that it is an associative Hopf algebra with 
comultiplication $\Delta$, antipode $a$, and 
counit $\epsilon$ given by 
\bac
\Delta(e_i)&=&e_i\otimes 1+k_i\otimes e_i,\nn\\        
\Delta(f_i)&=&f_i\otimes k_i^{-1} +1\otimes f_i,\nn \\           
\Delta(k_i)&=&k_i\otimes k_i,\qquad \Delta(q^d)=
q^d\otimes q^d,\nn \\ 
a(e_i)&=&-k^{-1}_ie_i,\qquad a(f_i)=-f_ik_i,\nn \\
a(k_i)&=&k_i^{-1},\qquad a(q^d)=q^{-d},\nn \\ 
\epsilon(e_i)&=&\epsilon(f_i)=0,\qquad \epsilon(k_i)=1,\qquad
i=0,1.
\label{coco}
\ea
The
special element $\gamma=k_0k_1$
commutes with the whole   $\uq$ algebra 
and acts as $q^k$ on 
its highest weight
representations. Here, $k$ is referred to as the level of the
representation.

The Chevalley generators are only associated with
the simple roots $\al_0$ and $\al_1$ of $\uq$. 
As in the classical case, it is important to find 
all  $\uq$ generators associated with the 
infinite-dimensional set of roots 
$\{\pm \alpha+n\delta; n\in \Z\}\cup\{n\delta;
n\in \Z\backslash\{0\}\}$, 
with $\al=\al_1$ and $\delta=
\al_0+\al_1$.  The Serre relations
will  become redundant in the commutation relations
of such generators.
These generators, found by Drinfeld, 
\cite{Dri86} are $\{E_{n}^{\pm},
\al_m, K^{\pm 1}, q^{\pm d}, \gamma^{\pm 1/2}; n\in \Z,
m\in\Z^*=\Z\backslash \{0\}\}$
with defining relations
\bac
&&\g^{1/2}\g^{-1/2}=\g^{-1/2}\g^{1/2}=1,\quad
{[\g^{\pm 1/2}, y]}=0,\quad \forall y\in \uq,\label{Eq1}\nn\\
&&KK^{-1}=K^{-1}K=1,\quad K\al_nK^{-1}=  \al_n,\quad 
KE^\pm_n K^{-1} =q^{\pm 2} E^\pm_n,\nn\\
&&q^{d}q^{-d}=q^{-d}q^{d}=1,\quad Kq^{\pm d}K^{-1}=q^{\pm d},
\quad
q^dE^\pm_nq^{-d}=q^nE^\pm_n,\quad 
q^d\al_nq^{-d}=q^n\al_n,\label{Eq7}\nn\\
&&{[\al_n,\al_m]} = {[2n](\g^{n}-\g^{-n}) \ov
n(\q)}\delta_{n+m,0},\label{Eq8}\nn\\
&& {[\al_n,E^{\pm}_m]}=
\pm{\g^{\mp |n|/2}[2n]\over n}
E^\pm_{n+m},\label{Eq9}\nn\\
&&{[E^+_n,E^-_m]} = {\g^{(n-m)/2}\Psi_{n+m}-\g^{(m-n)/2}
\Phi_{n+m}\over q-q^{-1}},\label{Eq10}\nn\\
&&E^\pm_{n+1}E^\pm_m-q^{\pm 2}E^\pm_mE^\pm_{n+1}=
 q^{\pm 2}E^\pm
_nE^\pm_{m+1}-E^\pm_{m+1}E^\pm_n, \label{Eq11}
\ea
where
$\Psi_n$ and $\Phi_n$
are given by
the mode expansions of the currents
$\Psi(z)$ and $\Phi(z)$, 
which are themselves defined by
\bac
\Psi(z)&=&\sum\limits_{n\geq 0}\Psi_nz^{-n}=K
\exp\{(\q)\sum\limits_{n>0}\al_nz^{-n}\},\nn\\
\Phi(z)&=&\sum\limits_{n\leq 0}\Phi_nz^{-n}=K^{-1}
\exp\{-(\q)\sum\limits_{n<0}\al_nz^{-n}\}.
\label{algebra}
\ea
Here $z$ is a formal variable and
\be
K=\Psi_0=\Phi_0^{-1}= q^{\al_0}.
\label{KP}\ee

The algebra isomorphism $\rho$ between   
the Chevalley generators
and the Drinfeld generators is given explicitly by
\bac
&&\rho:\quad k_0\ra \g K^{-1},\quad k_1\ra K,\nn\\
&&\rho:\quad e_{1}\ra E^{+}_0,\quad f_{1}\ra E^{-}_0,\nn \\
&&\rho:\quad e_{0}\ra E^-_{1}K^{-1},\quad  f_{0}\ra KE^+_{-1}.
\label{iso}
\ea
Using this isomorphism
one can re-express the  
comultiplication (\ref{coco}) in terms of the 
Drinfeld generators as
\bac
\!\!\!\!\!\!& &\Delta(E^+_n)=E^+_n\otimes\gamma^{kn}+  
\gamma^{2kn}K\otimes  
E^+_n+ \sum_{i=0}^{n-1}
\gamma^{k(n+3i)/2}\Psi_{n-i}\otimes \gamma^{k(n-i)}      
 E^+_i\: {\rm mod}   \: {N_-}\otimes {N_+^2},\nn \\             
 \!\!\!\!\!\!& &\Delta(E^+_{-m})=
E^+_{-m}\!\otimes\!\gamma^{-km}\!+\!  
K^{-1}\!\otimes\! E^+_{-m}+             
\sum_{i=0}^{m-1}\gamma^{{k(m-i)\over 2}}
\Phi_{-m+i}\otimes \gamma^{k(i-m)}  
 E^+_{-i}           
\: {\rm mod  }\:  
 N_-\otimes N_+^2,\nn \\          
\!\!\!\!\!\!& &\Delta(E^-_{-n})=E^-_{-n}\!\otimes\!
\gamma^{-2kn}K^{-1}\!+\!      
  \gamma^{-kn}\!\otimes \!E^-_{-n}\!+\!\!
\sum_{i=0}^{n-1}\!\gamma^{-k(n-i)} 
E^-_i\!\otimes\!\gamma^{{-k(n+3i)\over 2}}\Phi_{i-n}\: 
{\rm mod  }\:   
N_-^2\!\otimes\! N_+,\nn \\          
\!\!\!\!\!\!& &\Delta(E^-_m)=\gamma^{km}\otimes 
E^-_m+E^-_m\otimes   
K +        
\sum_{i=1}^{m-1}\gamma^{k(m-1)}E^-_m\otimes 
\gamma^{-k(m-i)/2}          
\Psi_{m-i}\:{\rm mod  }\: N_-^2\otimes N_+,\nn \\             
\!\!\!\!\!\!& &\Delta(\al_m)=\al_m\otimes\gamma^{km/2}+            
\gamma^{3km/2}\otimes \al_m\: {\rm mod  }\: 
N_-\otimes N_+,\quad \nn \\             
\!\!\!\!\!\!&
&\Delta(\al_{-m})=\al_{-m}\otimes\gamma^{-3km/2}+
\gamma^{-km/2}     
\otimes  \al_{-m}\: {\rm mod  }\: 
N_-\otimes N_+,\quad\nn\\             
\!\!\!\!\!\!& &\Delta(K^{\pm 1})=K^{\pm 1}\otimes  K^{\pm
1},\quad \nn\\           
\!\!\!\!\!\!& &\Delta(\gamma^{\pm \half})=
\gamma^{\pm\half}\otimes   
\gamma^{\pm \half},\nn \\ 
\!\!\!\!\!\!& &\Delta(q^{\pm d})=q^{\pm d}\otimes q^{\pm d},           
\label{comult}
\ea
where $m>0$, $n\geq 0$, and $N_\pm$ and $N_\pm^2$ are   
left ${\bf Q}(q)[\gamma^\pm,  \Psi_m, \Phi_{-n}; 
\: m, n\in {\bf   
\Z_{\geq 0}}]$ representations   
generated    
 by $\{E^\pm_m; \:m\in {\bf \Z}\}$   
and $\{E^\pm_m E^\pm_n; \:m, n\in {\bf \Z}\}$  respectively   
 \cite{ChPr91,Jimal92}.

Let us now use the formal variable  approach
 to    re-express
$\uq$ as a quantum current algebra with elements
$\{\Psi(z), \Phi(z), E^{\pm}(z), \g^{\pm 1/2}, q^{\pm d}\}$, 
where
\be
E^\pm(z)=\sum_{n\in\Z}E^\pm_nz^{-n},
\ee
and with defining relations
\bea
\g^{1/2}\g^{-1/2}=\g^{-1/2}\g^{1/2}&=&1,
\quad {[\g^{\pm 1/2}, y]}=0,
\quad \forall y\in \uq,\label{eq1}\\
{[\Psi(z),\Psi(w)]}&=&0,\label{eq2}\\
{[\Phi(z),\Phi(w)]}&=&0,\label{eq3}\\
\Psi(z)\Phi(w)&=&g(wz^{-1}\g)g(wz^{-1}\g^{-1})^{-1}
\Phi(w)\Psi(z), \label{eq4}\\
\Psi(z)E^{\ep}(w)&=&
g(wz^{-1}\g^{-\ep /2})^{-\ep}E^\ep(w)\Psi(z),\label{eq5} \\
\Phi(z)E^{\ep}(w)&=&g(zw^{-1}\g^{-\ep/2})^{\ep}
E^\ep(w)\Phi(z),\label{eq6}\\
{[E^\ep(z),E^{-\ep}(w)]}&=&\ep{\delta(zw^{-1}\g^{-\ep})
\Psi(w\g^{\ep/2})- \delta(zw^{-1}\g^\ep)\Phi(z\g^{\ep/2})\over
\q},\label{eq7}\\
(z-w q^{2\ep})E^{\ep}(z)E^{\ep}(w)&=&(z q^{2\ep}-w)
E^{\ep}(w)E^{\ep}(z),\label{eq8}\\
q^{d}E^{\ep}(z)&=&E^\ep(zq^{-1})q^d,\label{eq9}\\
q^{d}\Psi(z)&=&\Psi(zq^{-1})q^d,\label{eq10}\\
q^{d}\Phi(z)&=&\Phi(zq^{-1})q^d.
\label{eq11}
\eea
Here $\ep=\pm 1$ and $g(z)$ is meant to be the following
formal power series      in $z$:
\be
g(z)=\sum_{n\in \N}c_nz^n,
\ee
where the coefficients $c_n, n\in \N$ are determined from 
the Taylor
expansion of the function
\be
f(\xi)={q^2\xi-1\over \xi-q^2}=\sum_{n\in \Z_+}c_n\xi^n
\ee
at $\xi=0$ \cite{FrJi88}.
In the above relations we have
also introduced the 
$\delta$-function $\delta(z)$ which is defined as the
formal Laurent series
\be
\delta(z)=\sum_{n\in \Z}z^n,
\ee
and which plays a key role in the formal calculus approach.
Its main properties are summarized by the
following relations:
\bac
&&\delta(z)=\delta(z^{-1}),\nn \\
&&\delta(z)={1\over z}+{z^{-1}\over 1-z^{-1}},\nn \\
&&G(z,w)\delta(azw^{-1})=G(z,az)\delta(azw^{-1})=
G(a^{-1}w,w)\delta(azw^{-1} ),\quad a\in \C^*,
\label{del}\ea
where $G(z,w)$ is any operator whose
 formal  Laurent  expansion in $z$ and $w$ is given by
\be
G(z,w)=\sum_{n,m\in \Z}G_{n,m}z^nw^m.
\ee
Note that it is crucial that both $\delta(z)$ and
$G(z,w)$ have expansions in integral powers of $z$ and
$w$, otherwise the above properties of the 
$\delta$-function will not  hold.

The three relations  (\ref{eq9}), (\ref{eq10}) and (\ref{eq11})
translate the fact
that $E^{\pm}_n$, $\Psi_n$ and $\Phi_n$ are
homogeneous of the same degree $n$.

Unfortunately, the Hopf algebra structure of the current
algebra (2.9)-(2.19), which would follow from 
(\ref{comult}), is not available yet in a closed form. 
However, Drinfeld has succeeded in deriving
a closed formula for the comultiplication, although it is not
clear yet how it is related to the one defined in (\ref{coco}). 
It is given by \cite{Jin94}:
\bac
\Delta(\Psi(z))&=&\Psi(zq^{b_1})\otimes \Psi(zq^{b_2}),\nn\\
\Delta(\Phi(z))&=&\Phi(zq^{a_1})\otimes \Phi(zq^{a_2}),\nn\\
\Delta(E^{+}(z))&=&E^+(zq^{d_0})\otimes 1+
\Psi(zq^{d_1})\otimes E^+(zq^{d_2}),\nn\\
\Delta(E^-(z))&=&1\otimes 
E^-(zq^{e_0})+E^-(zq^{e_1})\otimes \Phi(zq^{e_1}),
\ea
where $a_i, b_i, c_i, d_i$ satisfy the following constraints:
\bac
b_1&=&a_1+1\otimes c,\quad e_0=a_2-c/2\otimes 1,\nn\\
d_0&=&a_1+1\otimes c/2,\quad d_2=a_2-3c/2\otimes 1,\nn\\
d_1&=&a_1+1\otimes 
c/2-c/2\otimes 1,\quad e_2=a_2-c/2\otimes 1-1\otimes c/2,\nn\\
e_1&=&a_1-1\otimes c/2,\quad b_2=a_2-c\otimes 1.
\ea

\subsection{\bf Representation theory of $\uq$}

Let us briefly recall the definitions of some
$\uq$ representations \cite{Idz93,JiMi94}.
For this, we still
need some notions from $\LS$ affine
algebra, which is generated by
$\{e_i,f_i,h_i,d; i=0,1\}$. We define on its Cartan
subalgebra $\hat{h}=\C h_0+\C h_1+\C d$ 
an invariant symmetric bilinear
form $(\quad,\quad)$ by
\be
(h_i,h_i)=2,\quad (h_i,h_{1-i})=-2, 
\quad (h_i,d)=\delta_{i,0},\quad
(d,d)=0,\quad i=0,1.
\ee
Let
$\hat{h}^*=\C \Lambda_0+\C \Lambda_1+\C \delta=
\C \al_0+\C \al_1+\C \Lambda_0$ be the dual space to
$h$ with
\be
<\Lambda_i,h_j>=\delta_{i,j},\quad
<\delta,d>=1,\quad
<\Lambda_i,d>=0,\quad
<\delta, h_i>=0,
\ee
where
\be
<\quad,\quad>:\quad \hat{h}^*\otimes \hat{h}\ra \C.
\ee
is the natural pairing,  $\Lambda_i$ are the
fundamental weights,  $\al_i$ are the positive
roots and $\delta=\al_0+\al_1$ is the null root. 
One can induce a symmetric bilinear form $(\quad,\quad)$ on 
$\hat{h}^*$ by
\bac
&&(\Lambda_i,\Lambda_j)={1\over 2}\delta_{i,1}\delta_{j,1},
\quad  (\Lambda_i,\delta)=1, \quad (\delta,\delta)=0,
\quad (\al_i,\al_i)=2,\nn\\
&&(\al_i,\al_{1-i})=2,\quad (\al_i,\Lambda_0)=
\delta_{i,0},\quad (\Lambda_0,\Lambda_0)=0,\quad i,j=0,1.
\ea
The weights $\lambda\in \hat{h}^*$ such that
\be
\lambda=n_{0}\Lambda_{0}+n_{1}\Lambda_{1},
\quad n_{0},n_{1}\in \N\backslash \{0\},
\ee
are called regular dominant integral weights, and $n_0+n_1=k$
is  the level that we have introduced previously.

Let us now turn to the representation theory of 
$\uq$ which is generated by 
$\{e_i,f_i,K^{\pm 1},\newline \gamma^{\pm 1},
q^{\pm d};\>i=0,1\}$.
Let $V$ be a $\uq$ representation and 
$\mu\in \hat{h}$, the subspace
$V_\mu\subset V$ defined by
\be
V_\mu =\{v\in V/ K^{\pm 1}v=q^{\pm 
<\mu,h_1>}v,\quad\gamma^{\pm 1}
v=q^{\pm k}v,\quad q^{\pm d}v=q^{\pm <\mu ,d>}v\},
\ee
is called a $\mu$-weight space, and any $v\in V_\mu$ is 
referred to
as a $\mu$-weight vector.
The representation $V$ becomes a weight
representation if it  is the direct sum of its weight spaces.
A $\uq$-highest weight vector  $v_\lambda$ in $V$ is a
$\lambda$-weight vector which satisfies the
additional condition
\be
e_iv_\lambda=0,\quad i=0,1.
\ee
The space $V$ is called a
$\uq$-highest weight representation if it is generated from a
$\lambda$-highest weight vector $v_\lambda$. In this case,
$v_\lambda$ is unique
(up to a multiplication by a scalar), and hence  $V$ 
is labelled by the
weight $\lambda$ as $V(\lambda)$.
This $V(\lambda)$ is called standard if it is
generated from a highest weight vector $v_\lambda$
with a dominant integral weight  $\lambda$ and such that
\be
f_i^{<\lambda,h_i>+1}v_\lambda=0,\quad i=0,1.
\ee
In  this case  it is also irreducible.
Let us now review the explicit constructions of the 
irreducible highest weight
representations (IHWR) in the case $k=1$, where they are also 
called basic representations. 

\subsection{Explicit construction of the $\uq$ basic
representations}

Let $U_q(\hat{h})$ be the
enveloping Heisenberg algebra 
embedded in $\uq$ and generated by $\{\alpha_n,
\g^{\pm 1/2}, q^{\pm d};\> n\in \Z^* \}$ with relations 
\cite{FrJi88}:
\bea
{[\al_n, \g^{\pm 1/2}]}&=&0,\quad \g^{1/2}\g^{-1/2}=\g^{-1/2}
\g^{1/2}=1,\\
{[\al_n,\al_m]}& =& \delta_{n+m,0}{[2n]\over n}.
\label{qH}\eea
Let 
$U_{q}(\hat{h^{-}})$ denote the Abelian 
subalgebra of $U_{q}(\hat{h})$ generated by
$\{\al_n; n<0\}$. 
Let $F$ be a Fock space constructed as 
$F=U_{q}(\hat{h^{-}})|0>$, 
where $|0>$ is a vacuum state with respect to the Heisenbeg 
algebra, i.e., 
\be 
\al_n|0>=0,\quad n\geq 0.
\ee
Now we define the action of $U_{q}(\hat{h})$ on $F$ as
\bac
\gamma^{\pm 1}&:&\quad x|0>\ra q^{\pm 1}x|0>,\nn\\
\al_n&:&\quad x|0>\ra \al_n x|0>,\quad n<0,\nn\\
\al_n&:&\quad x|0>\ra [\al_n,x]|0>,\quad n>0,\nn\\
q^{\pm d}&:&\quad x|0>\ra q^{\pm d}xq^{\mp d}|0>=q^{\pm d}x|0>, 
\quad  \forall x\in U_q(\hat{h^{-}}).
\ea
Then it is clear that $F$ is an irreducible representation
of $U_q(\hat{h})$. However,  since $F$ is 
annihilated by $\alpha_0$ it cannot be a $\uq$ IHWR, otherwise 
it would also be annihilated by $E^\pm(z)$, 
according to 
$q^{\alpha_0}E^\pm(z)q^{-\alpha_0}=q^{\pm 2} E^\pm(z)$.
This is obviously not the case. Therefore, we must extend
$F$ by some orthogonal space $W$ in order to make the
resulting space $G=F\otimes W$ an IHWR. In order to find
$W$ and the action of $\uq$ on $G$ let us introduce
two  vertex operators: 
\be
S_{\ep}^{\pm}(z)=\exp\{\pm \ep \sum_{n>0}{\al_{\pm n}\over
[n]}q^{-\ep n/2}z^{\mp n}\},\quad\ep=\pm,
\label{S}
\ee
which are viewed as  formal Laurent series in $z$ 
with coefficients acting on  $F$. 
Using  (\ref{qH}) we find the following commutation relations: 
\bac
S_{\ep}^{+}(z)S_{\epp}^{-}(w)&=&\displaystyle
(1-q^{-1-(\ep+\epp)/2}wz^{-1})^{\ep\epp}
(1-q^{1-(\ep+\epp)/2}wz^{-1})^{\ep\epp}
S_{\epp}^{-}(w)S_{\ep}^{+}(z),\nn\\
S_{\ep}^{\pm}(z)S_{\epp}^{\pm}(w)&=&
S_{\epp}^{\pm}(w)S_{\ep}^{\pm}(z),\nn\\
{[\al_n,S^+_{\ep}(z)]}&=&0,\nn\\
{[\al_{-n},S^+_{\ep}(z)]}&=&\displaystyle
-{\ep q^{-\ep n/2}z^{-n}[2n]\over n}S^+_{\ep}(z),\nn\\
{[\al_{-n},S^-_{\ep}(z)]}&=&0,\nn\\
{[\al_n,S^-_{\ep}(z)]}&=&\displaystyle
-{\ep q^{-\ep n/2}z^{n}[2n]\over n}S^-_{\ep}(z),
\label{SS}
\ea
where $n>0$ and $\ep=\pm$. 
To appreciate the usefulness of  $S^\pm_\epsilon(z)$ 
  in the
construction of $W$,
let us introduce the currents 
\be
Z^{\ep}(z)=S^-_{\ep}(z)E^{\ep}(z)S^+_{\ep}(z).
\label{ZOP}
\ee
Then, using 
\be
{[\al_n,E^{\ep}(z)]}={\ep q^{-\ep |n|/2}z^{n}[2n]
\over n}E^{\ep}(z),\quad \ep=\pm;\quad n\in 
\Z\backslash  \{0\},
\label{aln}\ee
which is equivalent to (\ref{eq5})
and (\ref{eq6}), we find
\bac
{[\al_n, Z^{\ep}(z)]}&=&0,\quad n\in \Z\backslash\{0\},\nn\\
{[\g^{\pm 1/2}, Z^{\ep}(z)]}&=&0,\quad \ep=\pm.
\label{aZ}
\ea
This means that the actions of $Z^\pm(z)$  on $G$ 
are orthogonal to 
that of $U_q(\hat{h^-})$  and therefore  
$Z^\pm(z)$ act nontrivially only on $W$. To be more specific
the actions of $U(\hat{h})$ and
$ Z^\pm(z)$  on $G=F\otimes W$ are given by
\bac
x&:& u\otimes v\ra xu\otimes v,\quad x\in U_q(\hat{h}),\nn\\
Z^\pm(z)&:& u\otimes v\ra u\otimes Z^\pm(z)v,
\quad u\in F,\quad v\in W.
\ea
Equations (\ref{eq5})
and (\ref{eq6})  imply
\bac
\Psi_{0}E^\pm(w)&=&
q^{\pm 2}E^\pm(w)\Psi_0,\nn\\
\Psi_0 Z^\pm(w)&=&
q^{\pm 2} Z^\pm(w)\Psi_0.
\label{PZ}
\ea
This means that $\Psi_0$ and $\Phi_0$ also act nontrivially only
on $W$, i.e.,
\bac
\Psi_0&:& u\otimes v\ra u\otimes \Psi_0 v,\nn\\
\Phi_0&:& u\otimes v\ra u\otimes \Phi_0 v, \quad u\in F,
\quad v\in W,
\ea
The currents $\Psi(z)$ and $\Phi(z)$ can also be written as
\bac
\Psi(z)&=&
S^{+}_{\ep}(zq^{-3\ep /2})S^{+}_{-\ep}(zq^{3\ep /2})
\otimes \Psi_{0},\nn\\
\Phi(z)&=&S^{-}_{\ep}(zq^{3\ep /2})S^{-}_{-\ep}
(zq^{-3\ep /2})\otimes \Phi_{0},\quad \ep=\pm.
\label{ps}\ea
For completeness let us also define the actions  
 of $q^{\pm d}$ and
$\gamma^{\pm 1}$ on $F\otimes W$ as
\bac
q^{\pm d}&:& u\otimes v\ra q^{\pm d}u\otimes q^{\pm d}v,\nn\\
\g^{\pm 1}&:& u\otimes v\ra q^{\pm 1}u\otimes v.
\ea
The homogeneous gradation of the $\uq$ currents 
is also inherited by $Z^\epsilon(z)$,
that is,
\be
q^{d} Z^{\epsilon}(z)q^{-d}= Z^{\epsilon}(zq^{-1}).
\label{hab}\ee

From the above analysis it is clear that the 
 action of the currents $E^{\ep}(z)$
on $F\otimes W$ 
decomposes  as  \be
E^{\ep}(z)=S^-_{-\ep}(zq^{-\ep })S^+_{-\ep }(zq^{\ep })
\otimes  Z^{\ep}(z),\quad \ep=\pm,
\label{XS}\ee
where we have used
\be
\bigl(S^{\pm}_{\ep}(z)\bigr)^{-1}=
S^{\pm}_{-\ep}(zq^{\pm\ep }),\quad \ep=\pm,
\label{PSS}
\ee
and (\ref{SS}).
Substituting $E^\ep(z)$ by the latter expression
 in (\ref{eq1})-(\ref{eq11}), 
we find the following algebra for  $Z^\ep(z)$:
\bac
&&{ Z^{\ep}(z) Z^{-\ep}(w)\over
(1-q^{-1}wz^{-1})
(1-qwz^{-1})}-{ Z^{-\ep}(w) Z^{\ep}(z)\over
(1-q^{-1}zw^{-1})
(1-qzw^{-1})}\nonumber\\
&&={\ep\over q-q^{-1}}
\left(\Psi_{0}\delta(zw^{-1}q^{-\ep })
-\Phi_{0}\delta(zw^{-1}q^{\ep })\right)\nonumber\\
&&={1\over q-q^{-1}}
\left(q^{\ep \al(0)}\delta(zw^{-1}q^{-1 })
-q^{-\ep\al(0)}\delta(zw^{-1}q)\right),
\label{equa1}\nn\\
&&w^{2} Z^{\ep}(z) Z^\ep(w)=z^{2}
 Z^\ep(w) Z^{\ep}(z).
\label{equa2}
\ea
These relations and (\ref{hab})  can in fact 
be solved explicitly for $Z^\epsilon(z)$. One finds  
\be
Z^{\ep}(z)=e^{\ep \al}z^{\ep\al_0+{1\over 2}(\al,\al)}=
e^{\ep \al}z^{\ep\al_0+1},\quad \ep=\pm,
\ee
where $\al$ is the $sl(2)$ positive simple root and 
$e^{\al}\in \C[P]$, $P=Q\cup(Q+\al/2)$ and $Q$ are the 
$sl(2)$ weight and root lattices, whereas $\C[P]$  
and $\C[Q]$ are  the corresponding Abelian group algebras 
respectively.
The elements $z^{\al_0}$ and $q^d$ act on $\C[P]$ as
\bac
z^{\al_0}e^\beta&=&z^{(\al,\beta)}e^\beta z^{\al_0},\nn\\
q^{d}e^{\beta}&=&e^{\beta}q^{d}q^{-\beta_0
-(\beta,\beta)/2},\quad \beta\in P.  
\ea

Finally, let us now summarize this construction.
The currents 
$E^{\ep}(z)$  act on $F\otimes \C[P]$ as
\be
E^{\ep}(z)=S^-_{-\ep}(zq^{-\ep })S^+_{-\ep }(zq^{\ep })\otimes
e^{\ep \al}z^{\ep\al_0+1}.
\label{cur}
\ee 
The subspaces $F\otimes \C[Q]$ and $F\otimes e^{\alpha\over
2}\C[Q]$ whose direct sum is $F\otimes\C[P]$ are 
$\uq$ invariant and
 irreducible. They are in fact isomorphic to the 
basic representations
$V(\Lambda_0)$ and $V(\Lambda_1)$ with highest weight states
 $1\otimes 1$ and $1\otimes e^{\al\over 2}$, 
respectively. 
This construction was first derived in
Ref. \cite{FrJi88}.
Unfortunately, in both cases the currents 
$Z^\pm(z)$ satisfy complicated generalized commutation 
relations which, contrary to the case $k=1$, cannot be solved in
terms  of just a group algebra \cite{BoVi96}. 
It turns out that one needs 
two more sets of Heisenberg generators in addition to  
the set $\{\al_n\}$. This solution seems to be necessary 
for the purpose of evaluating physical quantities which are
based on  explicit knowledge of the scalar product  
on $\uq$ representations. 

Although
explicit bases have been  
constructed  for  these 
representations in terms of the modes of 
$Z^{\epsilon}(z)$, the
scalar product on them has so far remained inaccessible.  
In the next section, we review the realization of $\uq$
in terms of three bosonic fields, whose modes satisfy 
simple Heisenberg algebras.
\newpage
\section{\bf Bosonization of higher level $\uq$ algebra}


Here, we extend the bosonic realization 
of $\uq$ given in (\ref{cur}) for $k=1$ to
arbitrary level $k$. 
For this purpose, 
it is most convenient to re-express (2.3) as a 
quantum current algebra generated by the set of currents
$\{E^\pm(z),\Psi(z),\Phi(z)\}$.  
This current algebra is given in terms of   
the following operator product expansions (OPE's):
\bac
\Psi(z)\Phi(w)&=&
{(z-wq^{2+k})(z-wq^{-2-k})\over (z-wq^{2-k})(z-wq^{-2+k})}
\Phi(w)\Psi(z),
\label{eeq1} \nn\\
\Psi(z)E^{\pm}(w)&=&
q^{\pm 2}{(z-wq^{\mp(2+k/2)})\over z-wq^{\pm (2-k/2)}}
E^\pm(w)\Psi(z), 
\label{eeq2}\nn\\
\Phi(z)E^{\pm}(w)&=&
q^{\pm 2}{(z-wq^{\mp(2-k/2)})\over z-wq^{\pm (2+k/2)}}
E^\pm(w)\Phi(z),
\label{eeq3}\nn\\
E^+(z)E^-(w)&\sim& {1
\over w(\q)}\left\{{\Psi(wq^{k/2})\over z-wq^k}-
{\Phi(wq^{-k/2})\over z-wq^{-k}}\right\},
\quad |z|>|wq^{\pm k}|,
\label{eeq4}\nn\\
E^-(z)E^+(w)&\sim& {1
\over w(\q)}\left\{{\Psi(wq^{k/2})\over z-wq^k}-
{\Phi(wq^{-k/2})\over z-wq^{-k}}\right\},
\quad |z|>|wq^{\pm k}|,
\label{eeq5}\nn\\
E^{\pm}(z)E^{\pm}(w)&=&
{(z q^{\pm 2}-w)\over z-w q^{\pm 2}}
E^{\pm}(w)E^{\pm}(z).
\label{op6}
\ea
Here, the symbol $\sim$  means that the regular terms  as $z$
approaches    $wq^{\pm k}$ are
being omitted. 
In the above OPE's the identification $\sum_{n\geq 0}z^n=
(1-z)^{-1}$ with $|z|<1$ is used.

The construction for general $k$ will be based on that
for $k=1$ \cite{FrJi88}. 
The natural generalization of (2.55) that is compatible with
the OPE's (\ref{eeq1}) is given by 
\bac
\!\!\!\!\Psi(z)&= &:V^+(zq^{k/2})V^-(zq^{-k/2}):\nn\\
&=&
q^{\sqrt{2k}\va_0}\exp\left(
\sqrt{2k}(\q)\sum_{n>0}\va_nz^{-n}\right),\nn\\
\!\!\!\!\Phi(z)&= &
:V^+(zq^{-k/2})
V^-(zq^{k/2}):\nn\\
&=&q^{-\sqrt{2k}\va_0}
\exp\left(-\sqrt{2k}(\q)\sum_{n<0}\va_nz^{-n}\right),\nn\\
\!\!\!\!E^{\pm}(z)&=&\sqrt{[k]}\psi^{\pm}(z)
V^{\pm}(z),
\label{ABE}
\ea
where 
\bac
V^\pm(z)&=&:e^{\pm i\varphi^\pm(z)}:,\nn\\
\va^{\pm}(z)&=&\va-i\va_0\ln{z}
+ik\sum_{n\neq 0}{q^{\mp |n|k/2}\over [nk]}\va_nz^{-n}.
\ea
Here $\psi^{\pm}(z)$ are deformed parafermionic
fields, $\va^{\pm}(z)$ are deformed bosonic fields. Their 
modes $\{\va, \va_n;\> n\in \Z\}$
obey the following Heisenberg  algebra:
\bac
{[\va_n,\va_m]}&=&{[2n][nk]\over 2kn}\delta_{n+m,0}\>,\nn\\
{[\va,\va_0]}&=&i\>.
\label{qhma}
\ea
Moreover, the bosonic normal ordering  symbol
:: indicates that in the product of fields between 
the two colons,
 the creation modes $\{\va_n, \va;n<0\}$ should
be moved to the left of the annihilation modes
$\{\va_n; n\geq 0\}$.
Henceforth, we will  use the standard convention that
operators defined at the same point are understood
to be normal ordered. 
Substituting $\Psi(z)$, $\Phi(z)$ and $E^\pm(z)$ in terms
of $V^\pm(z)$ and $\psi^\pm(z)$ in (\ref{eeq1}) 
we find the following
boson-parafermion algebra \cite{BoVi94}:
\bac
V^{\pm}(z)V^{\mp}(w)&=&
z^{-{2\over k}}{(q^{k+2}wz^{-1};q^{2k})_{\infty}\over
(q^{k-2}wz^{-1};q^{2k})_{\infty}}
:V^{\pm}(z)V^{\mp}(w):,\nn\\
V^{\pm}(z)V^{\pm}(w)&=&
z^{{2\over k}}{(q^{k\mp k-2}wz^{-1};q^{2k})_{\infty}\over
(q^{k\mp k+2}wz^{-1};q^{2k})_{\infty}}
:V^{\pm}(z)V^{\pm}(w):,\nn\\
\psi^{\pm}(z)\psi^{\mp}(w)&=&
z^{2\over k}
{(q^{k-2}wz^{-1};q^{2k})_{\infty}\over (q^{k+2}wz^{-1};q^{2k})_
{\infty}}\left({1\over (z-wq^{k})(z-wq^{-k})}
+{\rm regular\>as}\> z\ra wq^{\pm k}\right) ,\nn\\
\psi^{\pm}(z)\psi^{\pm}(w)&=&
{(wz^{-1})^{2\over k}(zq^{\pm 2}-w)\over
z-wq^{\pm 2}}
{(q^{k\mp k-2}zw^{-1};q^{2k})_{\infty}
(q^{k\mp k+2}wz^{-1};q^{2k})_{\infty}\over
(q^{k\mp k+2}zw^{-1};q^{2k})_{\infty}
(q^{k\mp k-2}wz^{-1};q^{2k})_{\infty}}
\psi^{\pm}(w)\psi^{\pm}(z),\nn\\
V^{\pm}(z)\psi^{\pm}(w)
&=&\psi^{\pm}(w)V^{\pm}(z)={\rm regular\>terms\>
as}\> z\ra w,\nn\\
V^{\pm}(z)\psi^{\mp}(w)
&=&\psi^{\mp}(w)V^{\pm}(z)={\rm regular\>terms\> 
as}\> z\ra w,
\label{OPEP}
\ea
where 
\be
(x;y)_\infty=\prod_{i=0}^\infty (1-xy^i).
\ee
When $k=1$ the latter algebra simplifies to
\bac
V^{\pm}(z)V^{\mp}(w)&=&
{1\over (z-wq)(z-wq^{-1})}
:V^{\pm}(z)V^{\mp}(w):,\nn\\
V^{\pm}(z)V^{\pm}(w)&=&
(z-wq^{1\mp 1})(z-wq^{-1\mp 1})
:V^{\pm}(z)V^{\pm}(w):,\nn\\
\psi^{\pm}(z)\psi^{\mp}(w)&=&
1+{\rm regular\>terms\>as}z\ra w,\nn\\
\psi^{\pm}(z)\psi^{\pm}(w)&=&
\psi^{\pm}(w)\psi^{\pm}(z),\nn\\
V^{\pm}(z)\psi^{\pm}(w)
&=&\psi^{\pm}(w)V^{\pm}(z)={\rm regular\> terms\>as}\>z\ra w,
\nn\\
V^{\pm}(z)\psi^{\mp}(w)
&=&\psi^{\mp}(w)V^{\pm}(z)={\rm regular\> terms\>as}\>z\ra w.
\ea
From these relations we see
that the parafermionic fields $\psi^{\pm}(z)$ can  be
identified with the identity operator. Therefore, we 
reproduce the construction of Ref. \cite{FrJi88} in the present
notations as
\bac
\!\!\!\!\Psi(z)&= &:V^+(zq^{1/2})
V^-(zq^{-1/2}):\nn\\
&=&
q^{\sqrt{2}\va_0}\exp\left(
\sqrt{2}(\q)\sum_{n>0}\va_nz^{-n}\right),\nn\\
\!\!\!\!\Phi(z)&= &
:V^+(zq^{-1/2})
V^-(zq^{1/2}):\nn\\
&=&q^{-\sqrt{2}\va_0}
\exp\left(-\sqrt{2}(\q)\sum_{n<0}\va_nz^{-n}\right),\nn\\
\!\!\!\!E^{\pm}(z)&=&
V^{\pm}(z).
\ea
To be more precise, we have $\alpha_n=\sqrt{2} \varphi_n$.

However, for $k>1$, $\psi^\pm(z)$ 
 cannot be identified with the   
 identity operator since they satisfy 
nontrivial OPE's. Unfortunately, these OPE's are 
generalized commutation relations. While one can still use
these  relations to construct  bases of the $\uq$
IHWR's,
one cannot easily use them for the 
purpose of evaluating scalar products  
of these IHWR's.
As we will see later however, 
it is precisely the scalar products
defined on the IHWR's that would lead to the physical
quantities, such as correlation functions and form 
factors. For this reason we need to
re-express the parafermionic fields $\psi^\pm(z)$ in terms
of two bosonic fields since we can 
always easily define  scalar products on  
Heisenberg representations.

\subsection{ Bosonization of higher level $\uq$}

We refer to the bosonic realization of $\uq$ for
arbitrary $k$,  
as the Wakimoto construction. 
For this we need  two extra Heisenberg 
algebras generated by $\{\va^2, \va^2_n\}$
and  $\{\va^{(3)},\va^{(3)}_n\}$ that bosonize $\psi^\pm(z)$. 
In the classical case, the  Wakimoto bosonization of  
$\psi^\pm(z)$ has the form 
\be
\psi^\pm(z)=i(a\partial\va^{(2)}(z)+b\partial\va^{(3)}(z))
e^{ic\va^{(2)}(z)+id\va^{(3)}(z)},
\ee
for some real parameters $a,b,c,d$ and free bosonic fields 
$\va^{(2)}(z)$ and $\va^{(3)}(z)$.
Therefore,  quantum   derivatives are naturally expected to
replace the classical ones.
Moreover, this bosonization  must be consistent with the
classical limit, i.e., $q\ra 1$. 
Finally, 
the OPE $\psi^+(z)\psi^-(w)$  is fixed to be 
that given in (\ref{OPEP}). 
Keeping track of the latter
requirements, we  define the general 
form of the bosonization of $\psi^\pm(z)$ as:
\be
[k]\psi^{\pm}(z)={\exp\{\pm i\sqs\va^{1,\pm}(z)\}\over z(\q)}
\left(\exp\{\pm i\sqs X_A^{(\pm)}(z)\}-
\exp\{\pm i\sqs X_B^{(\pm)}(z)\}   \right),
\ee
where  $\va^{1,\pm}(z)\equiv\va^\pm(z)$ 
and
\be
X_A^{(\pm)}(z)=
\va^{(2)}-i\va^{(2)}_0\ln{zq^{A^{(\pm)}_2}}-i\va^{(3)}_0
\ln{q^{A^{(\pm)}_3}}+
i\sum_{n\neq 0}\{A^{(\pm)}_2(n)
\va^{(2)}_n+A^{(\pm)}_3(n)
\va^{(3)}_n)\} {z^{-n}\over n}.
\label{XA}
\ee
The operators 
$X_B^{(\pm)}(z)$ is given by a similar expression 
to (\ref{XA}) with  $A$ being
replaced by $B$. The bosonic modes $\{\va^{(2)},  \va^{(2)}_n\}$
and $\{\va^{(3)},  \va^{(3)}_n\}$ satisfy the
following  Heisenberg algebras:
\bac
{[\va^{(j)}_n,\va^{(\ell)}_m]}&=&(-1)^{j-1}nI_j(n)\delta^{j,\ell}
\delta_{n+m,0},\nn\\
{[\va^{(j)},\va^{(\ell)}_0]}&=&(-1)^{j-1}i\delta^{j,\ell},
\qquad\qquad\qquad  j,\ell=2,3.
\ea
Here no sum with respect to $j$ is meant. 

First,  consistency with (\ref{eeq1})
requires the following  relations to  be satisfied:
\bac
\exp\{i\sqs X^{(+)}_B(z)\}\exp\{-i\sqs X^{(-)}_A(w)\}&=&
{z-wq^{k+2}\over q(z-wq^k)}\exp \left\{-{2\over k}\langle
\va^{1,+}(z)\va^{1,-}(w)\rangle\right\}\nn\\
&&\times :\exp\{i\sqs \left(X^+_B(z)-X^-_A(w)\right)\}:,
\\
X^{(+)}_B(wq^k)&=&X^{(-)}_A(w),\\
\exp\{i\sqs X^{(+)}_A(z)\}.\exp\{-i\sqs X_B^{(-)}(w)\}&=&
{q(z-wq^{-k-2})\over z-wq^{-k}}\exp \left\{-{2\over k}\langle
\va^{1,+}(z)
\va^{1,-}(w)\rangle\right\}\nn\\
&&\times:\exp\{i\sqs  \left(X_A^{(+)}(z)-
X_B^{(-)}(w)\right)\}:,\\
X_A^{(+)}(wq^{-k})&=&X_B^{(-)}(w),\\
\exp\{i\sqs X_A^{(+)}(z)\}.\exp\{-i\sqs X_A^{(-)}(w)\}&=&
q\exp \left\{-{2\over k}\langle
\va^{1,+}(z)\va^{1,-}(w)\rangle\right\}\\
&&\times:
\exp\{i\sqs \left(X_A^{(+)}(z)-
X_A^{(-)}(w)\right)\}:,\nonumber\\
\exp\{i\sqs X_B^{(+)}(z)\}.\exp\{-i\sqs X_B^{(-)}(w)\}&=&
q^{-1}\exp \left\{-{2\over k}\langle
\va^{1,+}(z)\va^{1,-}(w)\rangle\right\}\\
&&\times
:\exp\{i\sqs \left(X_B^{(+)}(z)- 
X_B^{(-)}(w)\right)\}:\nonumber.
\ea
Any solution to the above constraints 
for the  free parameters $A^{(\pm)}_2(n)$, 
$A^{(\pm)}_2$, $A^{(\pm)}_3(n)$,  
$A^{(\pm)}_3$, $B^{(\pm)}_2(n)$, $B^{(\pm)}_2$,
$B^{(\pm)}_3(n)$,  $B^{(\pm)}_3$, $I_2(n)$ and $I_3(n)$
leads to a deformation of the Wakimoto construction.
From these  relations we find
\bac
A^{(\pm)}_2(n)&=&q^{-nk}B^{(\mp)}_2(n),
\label{re1}\nn\\
A^\pm_3(n)&=&q^{-nk}B^{(\mp)}_3(n),
\label{re2}\nn\\
A^\pm_2&=&-B^{(\pm)}_2=k/2,
\label{re3}\nn\\
A^{(\pm)}_3&=&B^{(\mp)}_3=\pm{\sqrt{k(k+2)}\over 2},
\label{re4}
\ea
and
\bac
A^{(\pm)}_2(n)A^{(\pm)}_2(-n)I_2(n)-
A^{(\pm)}_3(n)A^{(\pm)}_3(-n)I_3(n)
&=&{k\over 2}\left({[n(k+2)]\over [nk]}-1\right), \quad n>0,\nn\\
A^{(\pm)}_2(n)A^{(\mp)}_2(-n)I_2(n)-
A^{(\pm)}_3(n)A^{(\mp)}_3(-n)I_3(n)&=&
{k\over 2}{[2n]\over [nk]},\quad n>0.
\label{mas}
\ea
Therefore, the only free parameters to be
determined are $A^{(\pm)}_2(n)$, $A^{(\pm)}_3(n)$, $I_2(n)$
and $I_3(n)$. 
These are restricted to 
satisfy the set of general  ``master"
equations (\ref{mas}). Each solution 
to these  equations yields
a  particular Wakimoto construction. 
In this manner, all  these constructions were re-derived 
in Ref. \cite{Bou93}.
For example, 
the first one is obtained by fixing the above parameters as 
\cite{Boual93}:
\bac
A^{(\pm)}_2(n)&=&q^{-nk/2},\nn\\
A^{(\pm)}_3(n)&=&\pm {1\over 2}
\sqrt{{k+2\over k}}(q^{-nk}-1),\nn\\
I_2(n)&=&{k\over 4}{[n(k+2)]+[2n]-[nk]\over [nk]},\nn\\
I_3(n)&=&{k^2\over k+2}{[n][n(k+2)/2]\over [nk][nk/2]}.
\ea
The second Wakimoto 
bosonization corresponds to the following values for
the parameters \cite{Mat92}:
\bac
A^{(\ep)}_2(n)&=&{nk\over [nk]}q^{-nk/2},\nn\\
A^{(-\ep)}_2(n)&=&{nkq^{-nk/2}\over [2n]}\left(
{[n(k+2)]\over [nk]}-1\right),\nn\\
A^{(\ep)}_3(n)&=&0,\nn\\
A^{(-\ep)}_3(n)&=&\ep\sqrt{k(k+2)}nq^{-nk/2}(q-q^{-1}){
[n]\over [2n]}, \nn\\
I_2(n)&=&{[2n][nk]\over 2kn^2},\nn\\
I_3(n)&=&{[2n][n(k+2)]\over 2n^2(k+2)},
\label{Mat1}
\ea
where $\ep$ is equal to $+$ or $-$ if
$n>0$ or $n<0$, respectively.
The third  bosonization corresponds to 
 the following  choices of parameters \cite{Bou93}:
\bac
A^{(\pm)}_2(n)&=&
\sqrt{{k+2\over 2}}{nk\over [nk]}q^{n(\pm   1-2)k/2},
\qquad\qquad n>0,\nn\\
A^{(\pm)}_2(n)&=&\sqrt{{k+2\over 2}}{nk\over [nk]}
q^{\pm nk/2},\qquad\qquad n<0,\nn\\
A^{(+)}_3(n)&=&-\sqrt{2k}{n\over [2n]}q^{n(f-k/2)},\nn\\
A^{(-)}_3(n)&=&-\sqrt{2k}{n\over [2n]}q^{n(f-2-3k/2)},\nn\\
I_2(n)&=&{[nk][n(k+2)]\over n^2k(k+2)}q^{|n|k},\nn\\
I_3(n)&=&{[2n]^2\over 4n^2}.
\ea
Here $f$ is a free parameter.

In Ref. \cite{Bou93}, it is shown that all these 
bosonizations can be obtained
from one another through linear 
combinations of 
their modes. Let us now show that another 
linear combination of these modes leads
to the bosonization of  Ref. \cite{Shi92}. We will   
 be using the latter in the
rest of  this paper for simplification reasons. 
Indeed, from the linear combinations
\bac
a_{X_1,n}&=&\sqrt{2k}q^{-2n-|n|}{[n(k+2)]\over [nk]}
\varphi^{(1)}_n+
\sqrt{k(2+k)}q^{-2n-(k+1)|n|} 
{[2n]\over [nk]}\varphi^{(2)}_n,\nn\\
a_{X_2,n}&=&-\sqrt{2k}q^{-n(2+k)-{k+2)|n|\over 2}} 
{[2n]\over [nk]}\varphi^{(1)}_n-
\sqrt{k(2+k)}q^{-n(2+k)-{(k+2)|n|\over 2}} 
{[2n]\over [nk]}\varphi^{(2)}_n,\nn\\
a_{X_3,n}&=&2\varphi^{(3)}_n,\nn\\
\sqrt{{k\over 
2}}a_{X_1,0}&=&(2+k)\varphi^{(1)}_0+\sqrt{2(2+k)}
\varphi^{(2)}_0,\nn\\
\sqrt{{k\over
2}}a_{X_2,0}&=&-2\varphi^{(1)}_0-\sqrt{2(2+k)}
\varphi^{(2)}_0,\nn\\
\sqrt{{k\over 2}}a_{X_3,0}&=&\sqrt{2k}\varphi^{(3)}_0,\nn\\
\sqrt{{k\over
2}}Q_{X_1}&=&i(2+k)\varphi^{(1)}+i\sqrt{2(2+k)}
\varphi^{(2)},\nn\\
\sqrt{{k\over 2}}Q_{X_2}&=&-2i
\varphi^{(1)}-i\sqrt{2(2+k)}\varphi^{(2)},\nn\\
\sqrt{{k\over 2}}Q_{X_3}&=&i\sqrt{2k}\varphi^{(3)},
\ea
we reproduce  the bosonization of Ref. \cite{Shi92} 
in terms of three fields $X_1,X_2,X_3$ as:
\bac
E^+(z)&=&{1\over z(q-q^{-1})}
\sum_{\ep=\pm 1}\ep E^+_\ep(z),\nn\\
E^-(z)&=&{1\over z(q-q^{-1})}
\sum_{\epsilon=\pm 1}\epsilon E^-_\epsilon(z),
\label{3fields}\ea
where
\bac
E^-_\epsilon(z)&=& \exp\{\partial X_1^{(\epsilon)}(q^{-2}z;
-{k+2\over 2})+X_2(2|q^{(\epsilon-1)(k+2)}z;-1)+
X_3(2|q^{(\epsilon-1)(k+1)-1}z;0)\},\nn\\
E^+_\ep(z)&=& \exp\{-X_2(2|q^{-k-2}z;1)-
X_3(2|q^{-k-2+\ep}z;0)\},
\ea
and
\be
\partial X_1^{(\epsilon)}(q^{-2}z;
-{k+2\over 2})=\epsilon \{(q-q^{-1})\sum_{n=1}^\infty 
a_{X_1,\epsilon n}z^{-\epsilon n}q^{(2\epsilon -{k+2\over 2})n}
+a_{X_1,0}\log(q)\}.
\ee
Here 
a deformed bosonic field is denoted by
\bac
X(L;M,N|z,\alpha)&=&-\sum_{n\neq 0}{[Ln]\over [Mn][Nn]}
a_{X,n}z^{-n}q^{|n|\alpha}+{L\over MN}a_{X,0}\ln(z)+{L\over
MN}Q_X,\nn\\
X(N|z,\alpha)&=&X(L;L,N|z,\alpha)=
-\sum_{n\neq 0}{a_{X,n}z^{-n}q^{|n|\alpha}\over [Nn]}
+{a_{X,0}\ln(z)\over N}+{Q_X\over N}.
\ea
It has the following two-point matrix elements:
\bac
<X(L;M,N|z,\alpha)X(L^\prime;M^\prime,N^\prime|w,
\alpha^\prime)>
&=&-\sum_{n>0}
{[Ln][L^\prime n][a_{X,n},a_{X,-n}]\over 
[Mn][M^\prime n][Nn] [N^\prime n]}
z^{-n}w^nq^{n(\alpha+\alpha^\prime)}\nn\\
&&+{LL^\prime [a_{X,0},Q_X]\over MM^\prime NN^\prime} \ln(z),
\nn\\
<X(N|z,\alpha)X(N^\prime|w,\alpha^\prime)>
&=&-\sum_{n>0}
{[a_{X,n},a_{X,-n}]\over 
[Nn] [N^\prime n]}
z^{-n}w^nq^{n(\alpha+\alpha^\prime)}\nn\\
&&+
{ [a_{X,0},Q_X]\over NN^\prime} \ln(z),
\ea
with
\bac
{[a_{X_1,n}, a_{X_1,m}]}&=&{[2n][(k+2)n]\over n}\delta_{n+m,0},
\nn\\
{[a_{X_2,n}, a_{X_2,m}]}&=&-{[2n]^2\over n}\delta_{n+m,0},\nn\\
{[a_{X_3,n}, a_{X_3,m}]}&=&{[2n]^2\over n}\delta_{n+m,0},\nn\\
{[a_{X_1,0}, Q_{X_1}]}&=&2(k+2),\nn\\
{[a_{X_2,0}, Q_{X_2}]}&=&-4,\nn\\
{[a_{X_3,0}, Q_{X_3}]}&=&4.
\ea

\subsection{\bf Screening currents and their bosonization} 

A screening current $S(z)$ has  OPE's with the $\uq$ currents
that are either regular or q-total derivatives, i.e.,
\bac
x(z)S(w)&=&\pm S(w)x(z) \sim {\rm regular\> terms\>as}\> 
z\ra w,\nn\\
x(z)S(w)&=&\pm S(w) x(z) \sim \pm _{\alpha}{\cal D}_w 
\left({Z(w)\over z-w}\right), 
 \quad x(z)=
\Psi(z),\Phi(z), E^\pm(z),
\label{scr}\ea
where $Z(w)$ is some normal ordered operator and the 
q-derivative is defined  by 
\be 
_{\alpha}{\cal D}_z(f(z))={f(zq^\alpha)-f(zq^{-\alpha})\over
z(q-q^{-1})}.
\ee 
The solution to the above relations (\ref{scr}) 
that is widely used in the literature is given by
\be
S(z)=-{1\over z(q-q^{-1})}\sum_{\delta=\pm
1}\delta S_\delta(z),
\ee
where
\be
S_\delta(z)= \exp\{-X_1(k+2|q^{-2}z;
-{k+2\over 2})-X_2(2|q^{-k-2}z;-1)-
X_3(2|q^{-k-2+\delta}z;0)\}.
\ee 
Let us note that in Ref. \cite{BoWe95} a general  
theory for constructing
screening currents of "pure exponential type" was 
presented and
it allows for the derivation of two more screening currents.
This construction is quite complicated for
$U_q(\widehat{sl(n)})$, $n\geq 3$ but in the case of 
$\uq$ it can simply be summarised as follows: 
A screening current  can be automatically 
derived from each step current, i.e.,  
$E^+(z)$ and $E^-(z)$. All it takes is to write $E^\pm(z)$ in
the form
\be
E^\pm(z)=e^{X^\pm(z)}{\cal D}_z e^{Y^\pm(z)},
\ee
then up to coefficients $\ep_+,\ep_-=\pm 1$ one can   
read off automatically two screening currents 
$S^\pm(z)$ as
\be
S^\pm(z)=e^{\ep_\pm Y^\pm(z)}.
\ee
The coefficients $\ep_\pm$ are determined from the OPE's
$E^\pm(z)S^\pm(w)$ which are required to be q-total
derivatives.
However, there is some redundancy, in that, we need only
one of them, say $S^+(z)$ 
which we denote by $\eta(z)$ (to follow
the usual notations in the literature). From the
simple bosonization 
 of $E^+(z)$, 
that of $\eta(z)$ is then found as
\be
\eta(z)=e^{X_3(2|q^{-k-2}z;0)}.
\label{eta}\ee  

\newpage

\section{ Diagonalization of the higher spin Heisenberg model}

\subsection{Vertex operators}

The objects of main interest that bridge the 
mathematical aspects 
of the quantum affine algebra $\uq$ and the physical aspects
of the Heisenberg model are called vertex operators. 
More 
specifically,  the matrix elements
of local operators of this model (i.e., form factors) 
are expressed as traces
of products of vertex operators. Since the latter operators 
 act as intertwiners
of $\uq$ representations let us  briefly recall 
certain features of 
these  representations that are relevant for our purposes.

As described in section 2.2 the $\uq$  highest weight    
representations are denoted by $V(\lambda)$  
where $\{\lambda=\lambda_{m}=(k-m)
\Lambda_0+m\Lambda_1,~ m=0,\dots,k\}$  
and $\{\Lambda_0,\Lambda_1\}$ are the sets of 
$\uq$ dominant highest   
weights and fundamental weights, respectively. 
A $\uq$ evaluation representation 
$V^{(\ell)}(z)$ $(0\leq \ell\leq k)$  
is the affinitization of the spin $\ell/2$ representation.  
It is isomorphic to 
$V^{(\ell)}\otimes {\bf C}[z,z^{-1}]$, where 
$V^{(\ell)}$ is the   
$\uq$ $(\ell+1)$-dimensional representation,   
 with a basis $\{v_m^{(\ell)},~ 0\leq m\leq \ell\}$ 
such that: 
\bac
e_1v^{(\ell)}_m=[m]v^{(\ell)}_{m-1},\quad 
f_1v^{(\ell)}_m=[\ell-m]v^{(\ell)}_{m+1},\quad 
t_1v^{(\ell)}_m=q^{\ell-2m}v^{(\ell)}_m,\nn\\ 
e_0v^{(\ell)}_m=[\ell-m]v^{(\ell)}_{m+1},\quad 
f_0v^{(\ell)}_m=[m]v^{(\ell)}_{m-1},\quad 
t_0v^{(\ell)}_m=q^{2m-\ell}v^{(\ell)}_m, 
\ea
where it is understood that $v^{(\ell)}_m=0$ if $m>\ell$ or 
$m<0$.   
$V^{(\ell)}(z)$ is equipped with the following 
$\uq$ module  
structure  \cite{Idzal93}:  
\bac
e_1v^{(\ell)}_m\otimes z^n&=&[m]
v^{(\ell)}_{m-1}\otimes z^n,\nn\\ 
e_0v^{(\ell)}_m\otimes z^n&=&
[\ell-m]v^{(\ell)}_{m+1}\otimes z^{n+1},\nn\\ 
f_1v^{(\ell)}_m\otimes z^n&=&[\ell-m]v^{(\ell)}_{m+1}
\otimes z^n,\nn\\ 
f_0v^{(\ell)}_m\otimes z^n&=&[m]v^{(\ell)}_{m-1}\otimes 
z^{n-1},\nn\\ 
t_1v^{(\ell)}_m\otimes z^n&=&q^{\ell-2m}v^{(\ell)}_m\otimes 
z^n,\nn\\ 
t_0v^{(\ell)}_m\otimes z^n&=&
q^{2m-\ell}v^{(\ell)}_m\otimes z^n. 
\ea 
In terms of the Drinfeld realization this becomes 
\bac
\gamma^{\pm 1/2}v^{(\ell)}_m\otimes z^s&=&v^{(\ell)}_m\otimes 
z^s,\nn\\  
K v_m^{(\ell)}\otimes z^s&=&q^{\ell-2m}v^{(\ell)}_m\otimes 
z^s,\nn\\  
E^+_nv^{(\ell)}_m\otimes z^s&=&
q^{n(\ell+2-2m)}[m]v^{(\ell)}_{m-1}\otimes z^{s+n},\nn\\  
E^-_nv^{(\ell)}_m\otimes
z^s&=&q^{n(\ell-2m)}[\ell-m]v^{(\ell)}_{m+1}
\otimes z^{s+n},\nn\\  
\alpha_nv^{(\ell)}_m\otimes z^s&=&{1\over
n}\{[n\ell]-q^{n(\ell-m+1)}(q^n+q^{-n}) 
[nm]\}v^{(\ell)}_m\otimes z^{s+n}.  
\ea  
Let $V^{(\ell)*a^{\pm 1}}(z)$ be two evaluation 
representations dual to 
$V^{(\ell)}(z)$  through the 
action of the  
antipode $a$ and its inverse $a^{-1}$, respectively.
Then we have the following isomorphisms  \cite{Idzal93}: 
\be 
C^{(\ell)}_\pm:\, V^{(\ell)}(zq^{\mp 2})\curlra
V^{(\ell)*a^{\pm 1}}(z), 
\label{vo1}
\ee 
where  
\be\br{rcl} 
C^{(\ell)}_\pm(v^{(\ell)}_m\otimes (zq^{-2})^n)&=&{\cal
C}_{\pm, m}^{(\ell)}{v^{(\ell)*}_{\ell-m}}\otimes z^n 
,\\ 
{\cal
C}_{\pm m}^\ell&=&
(-1)^{m}q^{m^2-m(\ell\mp 1)}{\qbinom{\ell}{m}}^{-1}
,\quad 0\leq m\leq \ell. \er\ee
Here $\{v^{(\ell)*}_m,\, 0\leq m\leq \ell\}$ is the basis of 
$V^{(\ell)*}$ dual 
to $V^{(\ell)}$, and  the notation  
$\qbinom{x}{y}$ defines the q-analogue of the 
binomial coefficient as 
\bac 
\qbinom{x}{y}&=&{[x]!\over [y]![x-y]!},\nn\\ 
{[x]}!&=&[x][x-1]\dots [1]. 
\ea

Following the terminology of  
 Refs.  \cite{Daval92,Idzal93} there are two types of vertex
operators  referred to as type I vertex operators 
$\vt^{\mu,V^{(\ell)}}_{\lambda}(z)$ and 
type II vertex
operators $\~\Psi^{V^{(\ell)},\mu}_{\lambda}(z)$ respectively.   
They are defined as maps between $\uq$ representations  in  
the following way:  
\bac
&& \vt^{\mu,V^{(\ell)}}_{\lambda}(z):V(\lambda) 
\ra V(\mu) \otimes V^{(\ell)}(z),\nn\\
&&\~\Psi^{V^{(\ell)},\mu}_{\lambda}(z):V(\lambda) 
\ra V^{(\ell)}(z) \otimes V(\mu).  
\label{vo2}
\ea

Their main feature is that they satisfy the
following intertwining properties:  
\bac
&&\vt^{\mu, V^{(\ell)}}_\lambda(z) \circ x = \Delta(x)\circ  
\vt^{\mu, V^{(\ell)}}_\lambda(z),\nn\\
&&\~\Psi^{V^{(\ell)},\mu}_\lambda(z) \circ x = \Delta(x)\circ  
\~\Psi^{V^{(\ell)},\mu}_\lambda(z), ~~~~~\forall~ x\in \uq .  
\label{vo8} 
\ea  
This simply means that they commute with $\uq$ when they
act as above. 
Let us also note that all vertex operators  conserve 
the homogeneous gradation of $\uq$ highest weight modules, i.e.,
\bac
(d\otimes {\rm id})\vt^{\mu, V^{(\ell)}}_\lambda(z)-
\vt^{\mu, V^{(\ell)}}_\lambda(z)d&=&-z{d\over dz}
\vt^{\mu, V^{(\ell)}}_\lambda(z),\nn\\  
(d\otimes {\rm id})\~\Psi^{ V^{(\ell)},\mu}_\lambda(z)-
\~\Psi^{\mu, V^{(\ell)}}_\lambda(z)d&=&-z{d\over dz}
\~\Psi^{ V^{(\ell)},\mu}_\lambda(z).
\label{vvo11}
\ea
It is also convenient to define components of 
the above vertex operators    
\bac
\vt^{\mu, V^{(\ell)}}_\lambda(z) &=&   
\sum\limits_{j=0}^{\ell}  \vt^{\mu,
V^{(\ell)}}_{\lambda,j}(z)\otimes v^{(\ell)}_j\nn\\
\~\Psi^{V^{(\ell)},\mu}_\lambda(z) &=&   
\sum\limits_{j=0}^{\ell} v^{(\ell)}_j \otimes
\~\Psi^{V^{(\ell),\mu}}_{\lambda,j}(z).
\ea  

From a physical point of view only type I
vertex operators with $\ell=k$ are relevant since they
realize the action of  local operators on the eigenspace  
of the spin $s=k/2$
Heisenberg model. Moreover, only type II vertex operators with
$\ell=1$ are important for physical purposes because the
elementary excitations which are created by these operators
are always spin-$1/2$ spinons, regardless of the value $s=k/2$ 
of
the   local spin  variables. In these physical situations,
the vertex operators are normalized as follows:
\bac
\vt^{\sigma(\lambda), V^{(k)}}_\lambda(z)|\lambda>&=& 
|\sigma(\lambda)>\otimes v^{(k)}_m+\cdots,\quad
\lambda=m\Lambda_0+n\Lambda_1,\nn\\
\~\Psi^{ V^{(1)},\lambda_\pm}_\lambda(z)|\lambda>&=&
v_{\mp}\otimes |\lambda_\pm>+\cdots,\quad \lambda_\pm=
(m\mp 1)\Lambda_0+(n\pm 1)\Lambda_1.
\label{vo12}
\ea
Moreover, since the latter operators will be
interpreted as creation operators of local spin states
and annihilation operators of  eigenstates, respectively, one
needs their conjugate vertex operators which play opposite
roles. They are obtained through the isomorphisms 
(\ref{vo1}) and the
substitution of $V(z)$ by
of $V^{*a^{\pm 1}}(z)$ in (\ref{vo2}), that is,
\bac
\vt_{\lambda}^{\sigma(\lambda),V^{(k)*a^{\pm 1}}}(z)&=&
\alpha_{\lambda,\pm}^{\sigma(\lambda)} 
({\rm id}\otimes  C^{(k)}_\pm)
\vt^{\sigma(\lambda), V^{(k)}}_{\lambda}(zq^{\mp 2}),\nn\\
\~\Psi_{\lambda}^{V^{(1)*a^\pm},\lambda_\epsilon}(z)&=&
\alpha_{\lambda,\pm}^{\lambda_\epsilon} 
(C^{(1)}_\pm \otimes {\rm id})
\~\Psi^{ V^{(1)}, \lambda_\epsilon}_{\lambda}
(zq^{\mp 2}), \quad \epsilon=\pm,    
\ea
where the normalization constants are fixed by 
\cite{Idzal93}
\bac
\vt_{\lambda}^{\sigma(\lambda),V^{(k)*a^\pm}}(z)
|\lambda>&=&
|\sigma(\lambda)>\otimes v_n^{(k)*}+\cdots,\nn\\ 
\~\Psi_{\lambda}^{V^{(1)*a^\pm },
\lambda_\epsilon}(z)|\lambda>&=&
v^*_\epsilon \otimes |\lambda_\epsilon>+\cdots.
\ea
Here for notational convenience 
$v^*_0=v^*_+$ and $v^*_1=v^*_-$.
One then finds
\bac
\alpha^{\sigma(\lambda)}_{\lambda,\pm}&=&
{1\over {\cal C}_{\pm, k-n}^{(k)}},\nn\\
\alpha_{\lambda,\pm}^{\lambda_\epsilon}&=&
{1\over {\cal C}_{\pm, {1+\epsilon \over 2}}^{(1)}}. 
\ea

For the purpose of translating 
the action of a local operator on the
Hilbert space 
${\cal F}^{(m)}=
{\rm End}(V(\lambda_m))$, 
inverse vertex operators have also been introduced
in Ref. \cite{Idzal93}.  They  
intertwine  
$\uq$ representations in the following order: 
\be 
\vt_{\sigma(\lambda_m),V^{(k)}}^{\lambda_m}(z):
\quad V(\sigma(\lambda_m))\otimes  
V^{(k)}(z) 
\rightarrow V(\lambda_m), 
\ee 
where $\lambda_{m}=(k-m)\Lambda_0+m\Lambda_1$ and  
$\sigma(\lambda_m)=m\Lambda_0+(k-m)\Lambda_1$.  
They have $k+1$ components 
$ \vt_{\sigma(\lambda_m),V^{(k)},j}^{\lambda_m}(z)$ 
defined by 
\be 
\vt_{\sigma(\lambda_m),V^{(k)},j}^{\lambda_i}(z)|u>= 
\vt_{\sigma(\lambda_m),V^{(k)}}^{\lambda_m}(z)(|u>\otimes v_j)= 
({\rm id}_{V(\lambda_m)}\otimes <v_j,.>)
\vt^{\lambda_m, V^{(k)*a}}_{\sigma(\lambda_m)} 
(z)|u>, 
\ee 
for any $|u>\in V(\sigma(\lambda_m))$. 
Using the above normalizations one then finds 
\be 
\vt_{\sigma(\lambda_m),V^{(k)},j}^{\lambda_m}(z)=
{{\cal C}_{+,k-j}^{(k)}\over  
{\cal C}_{+,m}^{(k)}} 
\vt^{\lambda_m, V^{(k)}}_{\sigma(\lambda_m),k-j}(zq^{-2}),
\qquad 0\leq m\leq k. 
\label{vvo3}
\ee 
Furthermore, from the above definitions it has been shown  
in Ref. \cite{Idzal93} that   
\bac
\vt_{\lambda_{k-m},V^{(k)}}^{\lambda_m}(z)\circ 
\vt^{\lambda_{k-m}, 
V^{(k)}}_{\lambda_m}(z)&=& 
g_{\lambda_m}{\rm id}_{V(\lambda_m)},\nn\\
\vt_{\lambda_{m}}^{\lambda_{k-m},V^{(k)}}(z)\circ 
\vt^{\lambda_m}_{\lambda_{k-m},V^{(k)}}(z)&=& 
g_{\lambda_m}{\rm id}_{V(\lambda_{k-m})}, 
\label{}
\ea 
where $g_{\lambda_m}$ are scalar functions given by
\be
 g_{\lambda_m}=q^{(k-m)m} 
{\qbinom{k}{m}}
{(q^{2(k+1)};q^4)_\infty \over
(q^2;q^4)_\infty }.
\ee

\subsection{\bf Commutation relations}

In this section, we briefly summarize 
the commutation relations
of the vertex operators following Ref.  \cite{Idzal93}.
There,  they have been derived from 
the q-KZ equation. For this purpose we need to introduce
an $R$ matrix whose entries are the Boltzmann weights of 
 the classical two dimensional lattice model, 
and which is equivalent
to the higher spin Heisenberg model. The mathematical 
interpretation of this  matrix is that it intertwines 
the action of  $\uq$ 
 evaluation representations, that is,
\be
R(z_1,z_2)\Delta(x)=\Delta^\prime(x)R(z_1,z_2), \quad 
R(z_1,z_2)\in {\rm End}(V_{z_1}^{(m)}\otimes
V_{z_2}^{(n)}),\quad \forall x\in\uq,
\label{vo4}
\ee 
where $\Delta$ is the comultiplication of $\uq$, 
$\Delta^\prime =\Delta P$, with $P$ being the transposition
operator, i.e., $P(u\otimes v)=v\otimes u$. 
This $R$ matrix, denoted by $\bar R_{mn}(z_1/z_2)$, 
is uniquely defined through the above
intertwining relations and the following normalization:
\be
\bar R_{mn}(z_1/z_2)(v_0^{(m)}\otimes v_0^{(n)})=
v_0^{(m)}\otimes v_0^{(n)}.
\ee  
In the case of most physical interest, i.e., $m=n=k$, one finds
explicitly \cite{Idz93}
\be
\bar R_{kk}(z)P=\sum_{l=1}^k \prod_{r=1}^l \left(
{1-z q^{2k-2r+2} \over z-q^{2k-2r+2}}\right)\> P_l,
\ee
where $P_l$ is the projector defined as:
\bac
&&V^{(k)}\otimes V^{(k)}=V^{(2k)}\oplus V^{(2k-2)}\oplus \cdots
\oplus V^{(2)}\oplus V^{(0)},\nn\\
&&P_l:\quad V^{(k)}\otimes V^{(k)}\rightarrow V^{(2k-2l)}.
\ea
In Ref. \cite{Idzal93} other $R$ matrices, 
which are different from
$\bar R$ only by an overall scale factor, are also introduced as
\bac
R_{kk}(z)&=&r_{kk}\bar R_{kk}(z),\nn\\
R_{VV}(z)&=&-R_{11}(z),\nn\\
R_{VV^*}(z)&=&({\rm id}\otimes C_-)R_{VV}(zq^{-2})  
({\rm id}\otimes C_-)^{-1},\nn\\
R_{{V^*}V^*}(z)&=&(C_-\otimes C_-)R_{VV}(z)  
(C_-\otimes C_-)^{-1},\nn\\
r_{kk}(z)&=&z^{-k/2}{(q^2z;q^4)_\infty
(q^{2k+2}z^{-1};q^4)_\infty \over 
(q^2z^{-1};q^4)_\infty (q^{2k+2}z;q^4)_\infty}.
\label{vo3}
\ea  
From the solution to  the q-KZ equation for the 
two-point vacuum to vacuum
matrix elements of vertex operators, the following 
key commutation relations have been found \cite{Idzal93}:
\bac
\vt_{\sigma(\lambda)}^{\lambda,V_2^{(k)}}(z_2)
\vt_{\lambda}^{\sigma(\lambda),V_1^{(k)}}(z_1)
&=&
R_{kk}(z_1/z_2)
\vt_{\sigma(\lambda)}^{\lambda,V_1^{(k)}}(z_1)
\vt_{\lambda}^{\sigma(\lambda),V_2^{(k)}}(z_2),\nn\\
\~\Psi_{\sigma(\lambda)}^{V_1^{(1)},\sigma(\lambda_\pm)}(z_1)
\vt_{\lambda}^{\sigma(\lambda),V_2^{(k)}}(z_2)
&=&
\tau(z_1/z_2)^{-1}
\vt_{\lambda_\pm}^{\sigma(\lambda_\pm),V_2^{(k)}}(z_2)
\~\Psi_{\lambda}^{V_1^{(1)},\lambda_\pm}(z_1),\nn\\
\~\Psi_{\sigma(\lambda)}^{V_1^{(1)*a^{-1}},
\sigma(\lambda_\pm)}(z_1)
\vt_{\lambda}^{\sigma(\lambda),V_2^{(k)}}(z_2)
&=&
\tau(z_1/z_2)
\vt_{\lambda_\pm}^{\sigma(\lambda_\pm),V_2^{(k)}}(z_2)
\~\Psi_{\lambda}^{V_1^{(1)*a^{-1}},\lambda_\pm}(z_1),\nn\\
\~\Psi_{\mu}^{V_1^{(1)},\nu}(z_1)
\~\Psi_{\lambda}^{V_2^{(1)*a^{-1}},\mu}(z_2)
&=&
(-q)^{\pm \delta_{\lambda,\mu}\delta_{\mu,\mu^\prime}}
R_{VV^*}(z_1/z_2)\sum_{\mu^\prime}
\~\Psi_{\mu^\prime}^{V_2^{(1)*a^{-1}},\nu}(z_2)
\~\Psi_{\lambda}^{V_1^{(1)},\mu^\prime}(z_1)\nn\\
&&\times W\left.\left(\matrix{\lambda&\mu\cr 
\mu^\prime&\nu\cr}\right|q^{-2}z_1/z_2\right),\nn\\
\~\Psi_{\mu}^{V_1^{(1)},\nu}(z_1)
\~\Psi_{\lambda}^{V_2^{(1)},\mu}(z_2)
&=&
R_{VV}(z_1/z_2)\sum_{\mu^\prime}
\~\Psi_{\mu^\prime}^{V_2^{(1)},\nu}(z_2)
\~\Psi_{\lambda}^{V_1^{(1)},\mu^\prime}(z_1)
W\left.\left(\matrix{\lambda&\mu\cr 
\mu^\prime&\nu\cr}\right|z_1/z_2\right),\nn\\
\label{vo7}
\ea
where $\mu=\lambda_\pm$. The last commutation relation  
still holds if 
 $\~\Psi_{\lambda}^{V^{(1)},\mu}(z)$ and 
$R_{VV}(z)$ are replaced by 
$\~\Psi_{\lambda}^{V^{(1)*a^{-1}},\mu}(z)$ and 
$R_{V^{*}V^{*}}(z)$, respectively.
The $W$ factors are the face type Boltzmann weights and
are given in the appendix, and 
\bac
\tau(z)&=&z^{-1/2}{\theta_{q^4}(qz)\over
\theta_{q^4}(qz^{-1})},\nn\\
\theta_p(z)&=&(z;p)_\infty (pz^{-1};p)_\infty (p;p)_\infty.
\ea

Moreover, it is worthnoting that in the limit
$|z_1/z_2|\rightarrow 1$,  the products
$\~\Psi_{\mu}^{V_1^{(1)},\nu}(z_1)
\~\Psi_{\lambda}^{V_2^{(1)},\mu}(z_2)$ and 
$\~\Psi_{\mu}^{V_1^{(1)*a^{-1}},\nu}(z_1)
\~\Psi_{\lambda}^{V_2^{(1)*a^{-1}},\mu}(z_2)$ are
holomorphic, whereas 
$\~\Psi_{\lambda_\pm}^{V_1^{(1)},\nu}(z_1)
\~\Psi_{\lambda}^{V_2^{(1)*a^{-1}},\lambda_\pm}(z_2)$ have   
 simple poles   
when $\lambda_\pm=\nu$, and whose residues are given by
\be
{\rm Res}_{z_1=z_2}
\~\Psi_{\lambda_\pm}^{V_1^{(1)},\nu}(z_1)
\~\Psi_{\lambda}^{V_2^{(1)*a^{-1}},\lambda_\pm}(z_2)
d({z_1\over z_2})=(v_+\otimes v^*_++v_-\otimes v_-^*)
\otimes g^\pm_\lambda {\rm id}_{V(\lambda)},
\ee
where
\bac
g^+_\lambda&=&q^{-1}\xi(q^2;1,q^4){(q^2;p)_\infty 
(q^{2n+4};p)_\infty \over
(p;p)_\infty (q^{2n+2};p)_\infty},\nn\\
g^-_\lambda&=&\xi(q^2;1,q^4){(q^2;p)_\infty 
(pq^{-2n};p)_\infty \over
(p;p)_\infty (pq^{-2n-2};p)_\infty},\nn\\
\xi(z;x,y)&=&{(xz;p,q^4)_\infty (x^{-1}yz;p,q^4)_\infty \over
(q^2xz;p,q^4)_\infty (q^{-2}x^{-1}yz;p,q^4)_\infty},\nn\\
(z;x,y)_\infty&=&\prod_{m,n=0}^\infty (1-zx^{m}y^n).
\ea

\subsection{\bf Diagonalization of the higher spin Heisenberg 
model}

In this section, we  briefly summarize 
the diagonalization of
the higher spin Heisenberg model as performed in Ref. 
\cite{Idzal93}.
We do not go  into the details of the proofs which are  
already well presented in this reference.

In this framework, it is more
convenient  to diagonalize the transfer matrix 
instead. This is of course equivalent to the diagonalization 
of the Hamiltonian. For this crucial purpose, 
let us describe the main steps and results of 
Ref. \cite{Idzal93}:

1.$\>$ Given the $R$ matrix  $R_{kk}(z)$ (\ref{vo3}) 
the higher spin Hamiltonian is obtained from it as
\bac
PR_{kk}(z_1/z_2)&=&{\rm const}(1+uh_{l+1,l}+\cdots ),\quad 
u\rightarrow
0,\quad z_1/z_2=e^{u},\nn\\{\cal H}&=&
\sum_{l\in \Z}h_{l+1,l}=
{\rm const}\> z {d\over dz} \log(T(z))|_{z=1},
\label{vo5}
\ea   
where $T(z)$ is the transfer matrix.

2.$\>$ Due to the intertwining relations (\ref{vo4}) and the
thermodynamic limit (i.e. $\sum_{l\in\Z}$), ${\cal H}$ is
invariant under the quantum affine symmetry $\uq$
\be
{[\Delta^{(\infty)}(x),{\cal H}]}=0,\quad \forall x\in \uq,
\ee 
where $\Delta^{(\infty)}$ is the infinite comultiplication. This
is defined from  $\Delta$ in (\ref{coco}) as
\be
\Delta^{(\infty)}=\lim_{n\rightarrow \infty} 
\Delta^{(n)},\quad \Delta^{(n)}=(\Delta\otimes {\rm id})
\Delta^{(n-1)},\quad \Delta^{(1)}=\Delta.
\ee

3.$\>$ The Hilbert space ${\cal F}$ of ${\cal H}$  is given by
\bac
{\cal F}&=&\sum_{\lambda,\mu}{\cal F}_{\lambda\mu},\nn\\
{\cal F}_{\lambda\mu}&=&\sum_{\lambda,\mu}V(\lambda)\otimes 
V(\mu)^{*a}=\sum_{\lambda,\mu}{\rm Hom}(V(\mu),V(\lambda)), 
\ea
where $V(\mu)^{*a}$ is the dual representation obtained from 
$V(\mu)$ through
the antipode $a$. The second  definition of 
${\cal F}_{\lambda\mu}$ is more convenient in defining the 
action of $\uq$ on it. To see this, 
let $f\in  {\rm Hom}(V(\mu),V(\lambda))$
then the left action of $\uq$ is defined by
\be
xf=\sum x_{(1)}\circ f \circ a(x_{(2)}),\quad 
\Delta(x)=\sum x_{(1)}\otimes x_{(2)},\quad \forall x\in\uq.
\ee
Moreover, let ${\cal F}^r_{\lambda\mu}$ denote 
${\cal F}_{\lambda\mu}$ as a $\uq$ right representation, 
then we have 
\be
fx=\sum a^{-1}(x_{(2)})\circ f \circ x_{(1)},\quad f\in 
{\cal F}^r_{\lambda\mu}.
\ee
There is a unique vacuum state $|vac>_{\lambda}\in 
{\cal F}_{\lambda\lambda}$, and which  is identified with 
the  identity element ${\rm id}_{V(\lambda)}$. This also holds 
for $_{\lambda}<vac|\in {\cal F}^r_{\lambda\lambda}$. 
There is a natural scalar product 
\be
<f|g>={tr_{V(\lambda)}(q^{-2\rho} fg)\over tr_{V(\lambda)}
(q^{-2\rho})}, \quad f\in {\cal F}^r_{\lambda\mu},\quad
g\in {\cal F}_{\mu\lambda}
\ee
with the property
\be
<fx|g>=<f|xg>,\quad \forall x\in\uq.
\ee 
Let us mention that as defined above, ${\cal F}$ is free from
the potential divergences, which plug the formal space
$\cdots \otimes V^{(k)}\otimes V^{(k)}\otimes \cdots$ on
which ${\cal H}$  acts naturally.

4.$\>$ The eigenspace ${\cal F}$ has a Fock space 
structure (i.e.,
a particle or spinon picture). 
This means that all its states can be
obtained  through successive actions of
a set of creation operators on the vacuum state. They can also
be annihilated through a set of annihilation operators. 
These operators are defined through type II vertex operators 
with $\ell=1$.
More specifically, 
creation operators $\varphi^{*\lambda^\prime}_{\lambda,
\epsilon}(z)$ acting on ${\cal F}_{\lambda\mu}=
{\rm Hom}(V(\mu),V(\lambda))$ are defined by
\be
 \varphi^{*\lambda^\prime}_{\lambda,
\epsilon}(z): {\cal F}_{\lambda\mu}\rightarrow 
{\cal F}_{\lambda^\prime\mu},\quad f\rightarrow 
 \~\Psi^{*\lambda^\prime}_{\lambda,\epsilon}(z)\circ f,\quad 
f\in {\cal F}_{\lambda\mu}.
\ee
Similarly, the annihilation operators 
$\varphi^{\mu^\prime}_{\mu,
\epsilon}(z)$ acting  on ${\cal F}^r_{\lambda\mu^\prime}$ 
are defined by
\be
 \varphi^{\mu^\prime}_{\mu,
\epsilon}(z): {\cal F}^r_{\lambda\mu^\prime}\rightarrow 
{\cal F}^r_{\lambda\mu},\quad f\rightarrow f\circ 
 \~\Psi^{\mu^\prime}_{\mu,\epsilon}(z),\quad 
f\in {\cal F}^r_{\lambda\mu^\prime}.
\ee
With these definitions, it has been shown in Ref. 
\cite{Idzal93} that the
following commutation relations hold:
\bac
\varphi_{\mu,\epsilon_1}^{\nu}(z_1)
\varphi_{\lambda,\epsilon_2}^{\mu}(z_2)
&=&
\sum_{\mu^\prime,\epsilon^\prime_1,\epsilon^\prime_2}
R_{VV}(z_1/z_2)^{\epsilon^\prime_1,\epsilon^\prime_2}_
{\epsilon_1,\epsilon_2}
\varphi_{\mu^\prime,\epsilon^\prime_2}^{\nu}(z_2)
\varphi_{\lambda,\epsilon^\prime_1}^{\mu^\prime}(z_1)
W\left.\left(\matrix{\lambda&\mu\cr 
\mu^\prime&\nu\cr}\right| z_1/z_2\right),\nn\\
\varphi_{\mu,\epsilon_1}^{*\nu}(z_1)
\varphi_{\lambda,\epsilon_2}^{*\mu}(z_2)
&=&
\sum_{\mu^\prime,\epsilon^\prime_1,\epsilon^\prime_2}
R_{{V^*}V^*}(z_1/z_2)^{\epsilon^\prime_1,\epsilon^\prime_2}_
{\epsilon_1,\epsilon_2}
\varphi_{\mu^\prime,\epsilon^\prime_2}^{*\nu}(z_2)
\varphi_{\lambda,\epsilon^\prime_1}^{*\mu^\prime}(z_1)
W\left.\left(\matrix{\lambda&\mu\cr 
\mu^\prime&\nu\cr}\right|z_1/z_2\right),\nn\\
\varphi_{\mu,\epsilon_1}^{\nu}(z_1)
\varphi_{\lambda,\epsilon_2}^{*\mu}(z_2)
&=&
\sum_{\mu^\prime,\epsilon^\prime_1,\epsilon^\prime_2}
R_{VV^*}(z_1/z_2)^{\epsilon^\prime_1,\epsilon^\prime_2}_
{\epsilon_1,\epsilon_2}
\varphi_{\mu^\prime,\epsilon^\prime_2}^{*\nu}(z_2)
\varphi_{\lambda,\epsilon^\prime_1}^{\mu^\prime}(z_1)
W\left.\left(\matrix{\lambda&\mu\cr 
\mu^\prime&\nu\cr}\right|z_1/z_2\right),\nn\\
&&\times (-q)^{\pm \delta_{\lambda,\nu} \delta_{\mu,\mu^\prime}}
W\left.\left(\matrix{\lambda&\mu\cr \mu^\prime&\nu\cr}
\right|
z_1q^{-2}/z_2\right)\nn\\
&&+g_\lambda^\pm
\delta_{\delta,\nu}\delta_{\epsilon_1,\epsilon_2}\delta(z_1/z_2).
\ea
Here $\mu=\lambda_\pm$. These commutations relations
are typical for creation and annihilation operators in
integrable models of quantum field theory. 

Therefore, a typical 
n-particle  state in ${\cal F}_{\lambda_{n}\lambda}$ is
given by
\be
|z_n,\cdots,z_1>_{\epsilon_n,\cdots,\epsilon_1}^{\lambda_n,\cdot
,\lambda_1,\lambda}=  
\varphi^{*\lambda_n}_{\lambda_{n-1},\epsilon_n}(z_n)
\cdots 
\varphi^{*\lambda_1}_{\lambda,\epsilon_1}(z_1)|vac>_\lambda =
\~\Psi^{*\lambda_n}_{\lambda_{n-1},\epsilon_n}(z_n)\cdots
\~\Psi^{*\lambda_1}_{\lambda,\epsilon_1}(z_1).
\ee
Similarly, a typical state in 
${\cal F}^r_{\lambda\lambda^\prime_{m}}$ is represented as
\be
_{\epsilon^\prime_1,\cdots,\epsilon^\prime_m}^{\lambda,
\lambda^\prime_1,\cdot
,\lambda^\prime_m}
<z^\prime_1,\cdots,z^\prime_m|=_{\lambda}<vac|
\varphi^{\lambda}_{\lambda^\prime_1,
\epsilon^\prime_1}(z^\prime_1)\cdots
\varphi^{\lambda^\prime_{m-1}}_{\lambda_m^\prime,
\epsilon^\prime_1}(z^\prime_m)=
\~\Psi^{\lambda}_{\lambda^\prime_1,
\epsilon^\prime_1}(z^\prime_1)\cdots
\~\Psi^{\lambda_{m-1}^\prime}_{\lambda^\prime_m,
\epsilon^\prime_1}(z^\prime_m).
\ee
The scalar product of the latter two states  then
reads
\bac
&&
_{\epsilon^\prime_1,\cdots,\epsilon^\prime_m}^{\lambda,
\lambda^\prime_1,\cdot
,\lambda^\prime_m}
<z^\prime_1,\cdots,z^\prime_m
|z_n,\cdots,z_1>_{\epsilon_n,\cdots,\epsilon_1}^{\lambda_n,\cdot
,\lambda_1,\lambda}\nn\\
&&= {Tr_{V(\lambda)}\left(q^{-2\rho} 
\~\Psi^{\lambda}_{\lambda^\prime_1,
\epsilon^\prime_1}(z^\prime_1)\cdots
\~\Psi^{\lambda_{m-1}^\prime}_{\lambda^\prime_m,
\epsilon^\prime_1}(z^\prime_m)
\~\Psi^{*\lambda_n}_{\lambda_{n-1},\epsilon_n}(z_n)\cdots
\~\Psi^{*\lambda_1}_{\lambda,\epsilon_1}(z_1)\right)\over
Tr_{V(\lambda)}\left(q^{-2\rho}\right)}.
\ea

5.$\>$ To see that the Fock space of the above creation and
annihilation operators is the eigenspace of the Hamiltonian let
us recall the definition of the transfer matrix $T(z)$ 
in this framework. To simplify the notations let us write
\bac
\~\Phi^{\sigma(\lambda) V}_{\lambda}(z)&=&
\sum_{j}\Phi_j(z)\otimes v_j^{(k)},\nn\\
\~\Phi^{\mu V^{*a^{-1}}}_{\sigma(\mu)}(z)&=&
\sum_{j}\Phi_j^*(z)\otimes v_j^{(k)*}, 
\ea
and let $\Phi^{*t}(z)$ denote the transpose of 
$\~\Phi^{\mu V^{*a^{-1}}}_{\sigma(\mu)}(z)$. 
Then, $T(z)$ is  defined
by
\bac
&&T(z)=T^{\sigma(\lambda)\sigma(\mu)}_{\lambda\mu}(z):
{\cal F}_{\lambda\mu}=V(\lambda)\otimes V(\mu)^{*a}
\stackrel{\Phi(z)\otimes {\rm id}}\rightarrow 
V(\sigma(\lambda))\otimes V_z^{(k)}\otimes
V(\mu)^{*a}\nn\\
&&\stackrel{{\rm id}\otimes \Phi^{*t}(z)}\rightarrow
V(\sigma(\lambda))\otimes V(\sigma(\mu))^{*a}=
{\cal F}_{\sigma(\lambda)\sigma(\mu)}.
\ea
This means 
\be
T(z)(f)=\sum_{j} \Phi_j(z)\circ f\circ \Phi_j^*(z),\quad
\forall f\in{\cal F}_{\lambda\mu}.
\ee
In particular, $g_\lambda^{-1}T(1)$ is identified with  
the translation operator. Using the last equation, the
commutation relations of vertex operators (\ref{vo7}), and 
the identity
\be
\sum_{j} \Phi_j(z)\circ \Phi_j^*(z)=g_\lambda {\rm id},
\ee
one finds the following
important relations \cite{Idzal93}: 
\bac
&&T^{\sigma(\lambda)\sigma(\lambda)}_{\lambda\lambda}(z)
|vac>_\lambda=g_\lambda |vac>_{\sigma(\lambda)},\nn\\
&&T^{\sigma(\lambda^\prime)
\sigma(\lambda)}_{\lambda^\prime\lambda}(z_2)
 \varphi^{\lambda^\prime}_{\lambda,\epsilon}(z_1)=
\tau(z_1/z_2)
\varphi^{\sigma(\lambda^\prime)}_{\sigma(\lambda),\epsilon}(z_1)
T^{\sigma(\lambda)
\sigma(\mu)}_{\lambda\mu}(z_2),\nn\\
&&T^{\sigma(\lambda^\prime)
\sigma(\lambda)}_{\lambda^\prime\lambda}(z_2)
 \varphi^{*\lambda^\prime}_{\lambda,\epsilon}(z_1)=
\tau(z_1/z_2)^{-1} 
\varphi^{*\sigma(\lambda^\prime)}_{\sigma(\lambda),\epsilon}(z_1)
T^{\sigma(\lambda)
\sigma(\mu)}_{\lambda\mu}(z_2),\nn\\
&&\~T(z_2) \varphi^{\mu}_{\lambda,\epsilon}(z_1)
\~T(z_2)^{-1}=\tau(z_1/z_2)^2 
\varphi^{\mu}_{\lambda,\epsilon}(z_1),\nn\\
&&\~T(z_2) \varphi^{*\mu}_{\lambda,\epsilon}(z_1)
\~T(z_2)^{-1}=\tau(z_1/z_2)^{-2} 
\varphi^{*\mu}_{\lambda,\epsilon}(z_1),
\ea
where 
$\~T(z)=T_{\sigma(\lambda)\sigma(\mu)}^{\lambda\mu}(z) 
T^{\sigma(\lambda)\sigma(\mu)}_{\lambda\mu}(z)$.
Consequently, (\ref{vo5}) implies that
\bac
{[{\cal H}, \varphi^{\mu}_{\lambda,\epsilon}(z)]}&=&
-e(z)\varphi^{\mu}_{\lambda,\epsilon}(z),\nn\\
{[{\cal H}, \varphi^{*\mu}_{\lambda,\epsilon}(z)]}&=&
e(z)\varphi^{*\mu}_{\lambda,\epsilon}(z),
\ea
with 
\be
e(z)={\rm const}\> z{d\over dz} \log(\tau(z))={\rm const}\>
\sqrt{1-k^2\cos^2(p(z))},
\label{vo6}
\ee
where the momentum $p(z)$ of an elementary excitation
 (spinon) is as usual defined by
\be
\tau(z)=e^{-ip(z)},
\ee 
and $k$ should not be confused with the level of $\uq$ or the
spin $k/2$ of the Heisenberg model. It is the 
elliptic module associated with  a nome $-q$.  
Obviously, $e(z)$ is nothing but the energy of a 
spinon, which is
a spin 1/2 excitation. 
Equation (\ref{vo6}) is then the dispersion 
relation of a spinon and
which is indepedent of the value $k/2$ of local spin variables.

6.$\>$ Let $L$ be a local operator acting on the subspace  
$V^{(k)\otimes n}$, i.e., $L\in {\rm End}(V^{(k)\otimes n})$, 
then its action on ${\cal F}_{\lambda\mu}$ reads
\be 
L(f)={\cal L}_\lambda \circ f,\quad f\in {\cal F}_{\lambda\mu},
\ee
and where ${\cal L}_\lambda$, in turn, 
is defined in terms of type 
I vertex operators as 
\bac
{\cal L}_{\lambda}&=& (g_{\lambda} g_{\lambda^{(1)}}
\cdots g_{\lambda^{(n-1)}})^{-1}
\vt_{\lambda^{(1)},V}^{\lambda}(z_1)\circ 
(\vt_{\lambda^{(2)},V}^{\lambda^{(1)}}(z_2)\otimes 
{\rm id}_V)\circ\cdots\circ 
(\vt_{\lambda^{(n)},V}^{\lambda^{(n-1)}}(z_n)\otimes 
{\rm id}_{V^{\otimes (n-1)}})\nn\\
&&\circ ({\rm id}_{V(\lambda^{(n)})}\otimes L)\nn\\
&&\circ 
(\vt_{\lambda^{(n-1)}}^{\lambda^{(n)},V}(z_n)\otimes 
{\rm id}_{V^{\otimes (n-1)}})\circ \cdots \circ 
(\vt_{\lambda^{(1)}}^{\lambda^{(2)},V}(z_2)
\otimes {\rm id}_{V})\circ 
\vt_{\lambda}^{\lambda^{(1)},V}(z_1).
\label{vvo1}\ea
where $\lambda^{(l)}=\sigma(\lambda^{(l-1)})$ and 
$V=V^{(k)}$. In this sense we can think of type I vertex
operators as the creation of local spin states belonging to the
subspace $V^{(k)\otimes n}$. 
Furthermore, from the commutation relations (\ref{vo7}) 
it can
be shown that local operators always commute with the creation
operators, that is,
\be
{\cal L}_\mu \~\Psi^{*\mu}_{\lambda,\epsilon}(z)=
\~\Psi^{*\mu}_{\lambda,\epsilon}(z){\cal L}_\lambda.
\ee
A matrix element   
of $L$ (form factor)
is then  given by
\be
_{\lambda}<vac|L 
|z_n,\cdots,z_1>_{\epsilon_n,\cdots,\epsilon_1}^{\lambda_n,\cdot
,\lambda_1,\lambda}
={Tr_{V(\lambda)}\left(q^{-2\rho}{\cal L}_{\lambda} 
\~\Psi^{*\lambda_n}_{\lambda_{n-1},\epsilon_n}(z_n)\cdots
\~\Psi^{*\lambda_1}_{\lambda,\epsilon_1}(z_1)\right)
\over
Tr_{V(\lambda)}(q^{-2\rho})}.
\label{vo11}
\ee
In particular, the correlation functions of $L$ (its 
vacuum to vacuum 
expectation values) simplify to 
\be
_{\lambda}<vac|L 
|vac>_{\lambda}
={Tr_{V(\lambda)}\left(q^{-2\rho}{\cal L}_{\lambda}\right) 
\over
Tr_{V(\lambda)}(q^{-2\rho})}.
\label{vvo2}
\ee

Let us however note that 
the evaluation of the above
traces is a nontrivial task because the representations 
are infinite-dimensional and so the vertex operators. 
For this purpose we use the powerful bosonization 
technique, that is,  the representation of  $\uq$ currents,
screening currents, 
irreducible highest weight modules and  vertex operators
in terms of modes satisfying a simple Heisenberg algebra.  
We have already presented the bosonization of the currents
and screening currents and therefore  
in the next section we  focus just on the bosonization of the 
vertex operators and $\uq$  irreducible highest weight 
representations.

\newpage

\section{\bf Bosonization of  vertex operators and IHWM's}

\subsection{\bf Bosonization of vertex operators}

Here we will be using the bosonization of $\uq$ 
as given in (\ref{3fields}). 
The results of this section is spread throughout many
papers in the literature but we will be following 
closely  Ref. \cite{Kon94}, especially as far as the 
cohomology analysis is concerned.
Due to the reducibility of the Fock 
spaces as $\uq$ modules
the bosonization of the  vertex operators is a
little trickier. To see this, let us construct 
Fock spaces $F_{\ell,s,t}$ from the successive actions
of the creation modes $\{a_{X_1,-n_1},
a_{X_2,-n_2}, 
a_{X_3,-n_3},\> n_1,n_2,n_3 >0\}$ on the vector 
$|\ell;s,t>=e^{{\ell\over
2(k+2)}Q_{X_1}+{s\over 2}Q_{X_2}+{t\over 2}Q_{X_3}}|0;0,0>$.
Because the currents do not depend explicitly on $Q_{X_1}$ 
they act naturally on the larger Fock spaces $F_\ell=
\oplus_{s,t\in \Z}F_{\ell,s,t}$. However, as can easily be
verified from their characters, these Fock spaces are larger
than the $\uq$ Verma modules
$L(\lambda_\ell)$ with highest weights $\lambda_\ell=
(k-\ell)\Lambda_0+\ell\Lambda_1$.  Note that the 
particular state 
\be
|\ell>=|\ell;0,0>=|vac>_{\lambda_\ell}= |\lambda_\ell>
\ee 
is a common highest weight state to $F_\ell$, 
 $L(\lambda_\ell)$, and to the 
irreducible highest weight module $V(\lambda_\ell)$. 
Therefore 
a projector is needed to  project the larger space
$F_\ell$ on  $L(\lambda_\ell)$.
Such a projector
is easily constructed from the 
screening current $\eta(z)$ (\ref{eta})
and also the ``screening operator" $a_{X_2}+a_{X_3}$ 
acting on  $F_\ell$. Both these
operators commute up to total derivatives with all 
$\uq$ currents. The bosonizations
of  $\eta(z)$ (\ref{eta})
and its ``conjugate" $\xi(z)$ are given by
\bac
\eta(z)&=&\exp({X_3(2|q^{-k-2}z;0))},\nn\\
\xi(z)&=&\exp{(-X_3(2|q^{-k-2}z;0))}.
\ea
One can easily check that these two fields are conjugate to each
other in the following sense:
\be
\xi(z)\eta(w)=-\eta(w)\xi(z)={q^3\over z-w}+:\xi(z)\eta(w):.
\ee
This means that their zero modes in 
$\eta(z)=\sum_{n\in \Z} z^{-n-1}\eta_n$ and 
$\xi(z)=\sum_{n\in \Z}z^{-n}\xi_n$ 
satisfy the anticommutation relation
\be
\{\xi_0,\eta_0\}=1.
\ee
Therefore, $L(\lambda_\ell)$ is isomorphic to the following
restricted space $\~F_\ell$:
\be
V(\lambda_\ell)\sim  \~F_\ell=\oplus_{s,t\in\Z}{\rm Ker}_{\eta_0}
\left({\rm
Ker}_{a_{X_2}+a_{X_3}}(F_{\ell,s,t})\right)=\oplus_{s\in\Z}
{\rm Ker}_{\eta_0}
(F_{\ell,s,s}).
\ee

Let us now define the actions of the vertex operators so that
they intertwine  the Fock  
spaces $\~F_\ell$. Let the bosonization of type I vertex
operator be 
\be
\~\Phi_{\lambda_{\ell_1}}^{\lambda_{\ell_2},V^{(\ell)}}(z)=
g_{\lambda_{\ell_1}}^{\lambda_{\ell_2},V^{(\ell)}}(z)
\sum_{m=0}^\ell
\~\Phi^{(r)}_{\ell,m}(z)\otimes 
v_m^{(\ell)},\quad 2r=\ell+\ell_1-\ell_2, 
\label{vvo4}
\ee
where   
the screened vertex operator components 
$\~\Phi^{(r)}_{\ell,m}(z)$ are given
in terms of the bare ones 
$\~\Phi^{(\ell)}_{m}(z)$ by
\be
\~\Phi^{(r)}_{\ell,m}(z)=
\int_{q^\ell pz}^{q^\ell z\infty} d_{p}t_{1} \dots  
\int_{q^\ell pz}^{q^\ell z\infty} d_{p}t_{r} S(t_1)\dots S(t_r)
\~\Phi^{(\ell)}_{m}(z).
\label{vo13}
\ee
The bosonization of the latter bare components  
$\~\Phi^{(\ell)}_{m}(z)$ is in turn
derived from  that of $\~\Phi^{(\ell)}_{\ell}(z)$ through the
successive applications of the intertwining relations
(\ref{vo8}).
More specifically,
using these relations (\ref{vo8}), the   
comultiplication (\ref{comult}), and the fact that   
$N_+v^{(\ell)}_{0}=N_-v^{(\ell)}_\ell=0$,
$N_{\pm}v^{(\ell)}_m\in   
F[z,z^{-1}]  
v^{(\ell)}_{m\mp 1}$, one arrives at:  
\bac                    
{[E^{+}(w),\~\Phi^{(\ell)}_{\ell}(z)]}&=&0,\nn\\         
{[H_n, \~\Phi^{(\ell)}_{\ell}(z)]}&=&q^{2n}z^n{[n\ell]\over 
n} q^{k(n+|n|/2)}\~\Phi^{(\ell)}_{\ell}(z),\quad n\neq 0,\nn\\  
K\~\Phi^{(\ell)}_{\ell}(z)K^{-1}&=&q^{\ell} 
\~\Phi^{(\ell)}_{\ell}(z),\nn\\
\~\Phi^{(\ell)}_{m}(z)&=&{1\over {[\ell-m]!}}[\dots[[
\~\Phi^{(\ell)}_{\ell}(z),
E^{-}_0]_{q^{}},E^{-}_0]_{q^{2}},\dots,  
E^-_0]_{q^{\ell-m}}.  
\label{vo10}
\ea
Here the quantum commutator $[A,B]_{q^{x}}$ is defined by  
\be  
[A,B]_{q^x}=AB-q^xBA.  
\ee  
These relations can easily be solved for the  bosonization of  
$\~\Phi^{(\ell)}_{\ell}(z)$, and one finds
\be
\~\Phi^{(\ell)}_{\ell}(z)=\exp{(X_1(\ell;2,k+2|q^{k}z;{k+2\over
2}))}.
\label{vo9}
\ee 

Finally, due to the commutation relations
\bac
{[\~\Phi^{(\ell)}_{\ell}(z),\eta_0]}&=&0,\nn\\
{[\~\Phi^{(\ell)}_{\ell}(z),a_{X_2,0}+a_{X_3,0}]}&=&0,
\ea
one can easily see that the screened vertex operators
intertwine the restricted Fock spaces
as: 
\be
\~\Phi^{(r)}_{\ell,m}(z): \~F_{\ell_1}\rightarrow   
\~F_{\ell_2}\otimes V^{(\ell)}(z)
\ee
The normalization functions 
$g_{\lambda_{\ell_1}}^{\lambda_{\ell_2},V^{(\ell)}}(z)$ are 
 fixed by.
\be
\~\Phi_{\lambda_{\ell_1}}^{\lambda_{\ell_2},V^{(\ell)}}(z)
|\ell_1;0,0>=|\ell_2;0,0>\otimes v+\cdots.
\ee

In fact, unlike  the case  $k=1$ with the simple bosonization
(\ref{cur}), there is a second bosonization  here of the 
vertex operator satisfying the intertwining 
relations. To distinguish it from the first one given in
(\ref{vo9}),
let us denote it by  
\be
\hat \Phi_{\lambda_{\ell_1}}^{\lambda_{\ell_2},V^{(\ell)}}(z)=
\hat g_{\lambda_{\ell_1}}^{\lambda_{\ell_2},V^{(\ell)}}(z)
\sum_{m=0}^\ell
\hat \Phi^{(r)}_{\ell,m}(z)\otimes 
v_m^{(\ell)}, 
\ee
with the bare components $\hat \Phi^{(r)}_{\ell,m}(z)$ 
being related to 
$\hat \Phi^{(\ell)}_{\ell}(z)$ in the same manner as described by 
(\ref{vo10}). Moreover, the bosonization of 
$\hat \Phi^{(\ell)}_{\ell}(z)$  is given by
 \be
\hat \Phi^{(\ell)}_{\ell}(z)=\exp(X_1(k+1-\ell;2,{k+2\over
2}|q^{k}z;0)+X_2(k+1-\ell;1,2|z;-1)+X_3(k-\ell;1,2|z;0)).
\ee  
It is typical in the context of conformal field theory to 
refer to this second bosonized vertex operator as a conjugate
vertex operator.

With obvious modifications, a type II vertex operator  
intertwines also the restricted Fock spaces when 
they are screened as
\be
\~\Psi_{\lambda_{\ell_1}}^{V^{(\ell)},\lambda_{\ell_2}}(z)=
h_{\lambda_{\ell_1}}^{V^{(\ell)},\lambda_{\ell_2}}(z)
\sum_{m=0}^\ell v_m^{(\ell)}\otimes
\~\Psi^{(r)}_{\ell,m}(z),
\ee
with
\be
\~\Psi^{(r)}_{\ell,m}(z)=
\int_{0}^{q^{k-\ell}z} d_{p}t_{1} \dots  
\int_{0}^{q^{k-\ell}z} d_{p}t_{r}  
\~\Psi^{(\ell)}_{m}(z) S(t_1)\dots S(t_r).
\ee
Similarly, the bosonizations of the bare components 
$\~\Psi^{(\ell)}_{m}(z)$ are derived
recursively from that of  $\~\Psi^{(\ell)}_{0}(z)$ as shown
by the following explicit form of the intertwining relations 
(\ref{vo8}):
\bac                    
{[E^{-}(w),\~\Psi^{(\ell)}_{0}(z)]}&=&0,\nn\\                     
{[H_n, \~\Psi^{(\ell)}_{0}(z)]}&=&-z^n{[n\ell]\over 
n} q^{k(n-|n|/2)}\~\Psi^{(\ell)}_{0}(z),\quad n\neq 0,\nn\\  
K\~\Psi^{(\ell)}_{0}(z)K^{-1}&=&q^{-\ell} 
\~\Psi^{(\ell)}_{0}(z),\nn\\
\~\Psi^{(\ell)}_{m}(z)&=&{1\over {[m]!}}[\dots[[
\~\Psi^{(\ell)}_{0}(z),
E^{+}_0]_{q^{\ell}},E^{+}_0]_{q^{\ell-2}},\dots,  
E^+_0]_{q^{\ell-2(m-1)}}.  
\label{vvo9}
\ea
From the first three commutation relations
one finds the bosonization of $\~\Psi^{(\ell)}_{0}(z)$ 
\be
\~\Psi^{(\ell)}_{0}(z)=
\exp({X_1(\ell;2,k+2|q^{k-2}z;-{k+2\over 2})+
X_2(\ell;2,1|q^{-2}z;0)+X_3(\ell;2,1|q^{-2}z;0)}).
\ee 
Moreover, one can easily check 
\bac
{[\~\Psi^{(\ell)}_{0}(z),\eta_0]}&=&0,\nn\\
{[\~\Phi^{(\ell)}_{0}(z),a_{X_2,0}+a_{X_3,0}]}&=&0,
\ea
and hence the screened components intertwine the restriced 
Fock spaces as
\be
\~\Psi^{(r)}_{\ell,m}(z): \~F_{\ell_1}\rightarrow   
V^{(\ell)}(z)\otimes \~F_{\ell_2}.
\ee 
Here also there exists a second bosonization given by
\be
\hat \Psi_{\lambda_{\ell_1}}^{V^{(\ell)},\lambda_{\ell_2}}(z)
=
\hat h_{\lambda_{\ell_1}}^{V^{(\ell)},\lambda_{\ell_2}}(z)
\sum_{m=0}^\ell v_m^{(\ell)}\otimes
\hat \Psi^{(r)}_{\ell,m}(z),
\ee
with 
\be
{\hat \Psi}^{(\ell)}_{0}(z)=
\exp(-{X_1(\ell+2;2,k+2|q^{k-2}z;-{k+2\over 2})-
X_2(2|q^{-2}z;0)}).
\ee 
For the same reason as in the case of 
type I vertex operators, we refer
to ${\hat \Psi}^{(\ell)}_{0}(z)$ as a type II conjugate vertex
operator. 

Let us however recall that in conformal field theroy
conjugate vertex operators play
an important role in simplifying the calculation of 
correlation functions of
vertex operators only on a genus zero manifold, i.e., the 
sphere.  
 They do not enter in the calculation of higher genus correlation
functions, and in particular in the case of the correlation
functions on a torus. As we have  seen in (\ref{vo11}) 
the physical
correlation functions and form factors of the Heisenberg model 
 are  expressed as 
traces of vertex operators
in a similar manner as genus one correlation functions in
conformal field theory,  conjugate vertex operators 
 are then not particularly useful for our purposes. 
Therefore, we will not consider them any further.  

Let us now briefly describe how the above contours
of the Jackson integrals were defined in Ref. \cite{Kon94}. 
The basic requirement is the following commutation relation:
\be
{[\~\Psi^{(\ell)}_{\ell,0}(z), E^-(w)]}=0,
\ee    
which is one of the intertwining relations. 
Here, a q-difference appears in  the  OPE
\be
:\~\Psi^{(\ell)}_{\ell,0}(z)::S(t)::E^{-}(w):\sim
_{k+2}\partial_{t}\left(
{(q^{\ell+2}tz^{-1};p)_\infty \over (w-t)(q^{-\ell+2} 
tz^{-1};p)_\infty}
:\~\Psi^{(\ell)}_{\ell,0}(z)
e^{-X_1(k+2|q^{-2}t;{k+2\over 2})}:\right).
\ee
With this OPE one can easily check using the definition
of a Jackson integral (\ref{ap12}) that the 
following integral of a
total q-difference  
\be
\int_{0}^{a} d_pt _{k+2}\partial_{t}\left(
{(q^{\ell+2}tz^{-1};p)_\infty \over (w-t)(q^{-\ell+2} 
tz^{-1};p)_\infty}
:\~\Psi^{(\ell)}_{\ell,0}(z)e^{-X_1(k+2|q^{-2}t;{k+2\over 2})}:
\right)
\ee
vanishes if 
\be
a=q^{k-\ell}p^{-j}z,\quad j=0,1,\cdots
\ee
In the above contours  the first 
choice $j=0$ has been picked. 
The normalization functions are obtained 
through the conventions made in (\ref{vo12}) 
and the various OPE's.
One finds explicitly
\bac
h_{\lambda_\ell}^{V^{(1)},\lambda_{\ell+1}}(z)&=&
q^{-1-2\ell(k-2)s}z^{-2\ell s},\nn\\
h_{\lambda_\ell}^{V^{(1)},\lambda_{\ell-1}}(z)&=&
-q^{2+{\ell\over
2}-2(3\ell+4)s}z^{(\ell+2)s}B_q(1-2\ell s,-2s)^{-1},
\ea
with $s=1/2(k+2)$.

\subsection{\bf Irreducible highest weight representations}

Since correlation functions and form factors are
given
as traces over irreducible highest weight representations 
$V(\lambda_\ell)$, which are
embedded in the Verma modules $L(\lambda_\ell)$, 
one needs then another projection from the latter to the former.
This is provided through a BRST cohomology analysis. As in the
classical case a BRST operator $Q_n$ can be defined as
\be
Q_n=\int_0^\infty d_{p}t J_n(t),
\ee
with the BRST current $J_n(t)$ being given by
\be
J_n(t)={1-{\cal A}_s^n\over 1-{\cal A}_s}
\int_{0}^{q^2t} d_{p}t_{2}\cdots 
 \int_{0}^{q^2t} d_{p}t_{n} S(t)S(t_2)\cdots S(t_n),
\ee
and being single valued on the restricted Fock spaces
$\~F_{\ell_{n,\bar n}}\equiv \~F_{n,\bar n}$. Here, 
${\cal A}_s$ is defined by 
\bac
S(t_1)S(t_2)&=&{\cal A}_{s}(t_1/t_2)
S(t_2)S(t_1),\nn\\
{\cal A}_{s}(t)&=&t^{-2\over k+2} 
{\theta_p(pq^{-2}t)\over \theta_p(pq^{2}t)},\nn\\
{\cal A}_{s}&=&
{\cal A}_{s}(q^{-2-\epsilon}),
\ea
where,  
$\epsilon$ is a
regularization parameter and   
${\cal A}_{s}(t)$ is a pseudo constant, that is,
${\cal A}_{s}(tp^n)=
{\cal A}_{s}(t)$, $n\in 
\Z$.

Following
the analysis of the general case we consider
$\ell_{n,\bar n}=
n-\bar n {P\over P^\prime}-1$ with $n$ and $\bar n$
being integers, and  $P$ and 
$P^\prime$ coprime integers satisfying 
${P\over P^\prime}=k+2$. However, the highest weights 
$\lambda_{\ell_{n,\bar n}}=(k-\ell_{n,\bar n})\Lambda_0 +
\ell_{n,\bar n}\Lambda_1$  of $V(\lambda_\ell)$ 
are such that 
$1\leq n \leq P$ and $0\leq \bar n \leq P^\prime$ (our present
case corresponds to $P=k+2$ and $P^\prime=1$). Given this, 
it has been
partially shown  in Ref. \cite{Kon94} that 
\be
Q_nQ_{P-n}=Q_{P-n}Q_n=0,
\ee
and that the following sequence is a complex:\\
$$\cdots \stackrel{Q_{n}}\rightarrow \~F_{-n+2P,\bar n}
\stackrel{Q_{P-n}}\rightarrow 
\~F_{n,\bar n} \stackrel{Q_{n}}\rightarrow \~F_{-n,\bar n} 
\stackrel{Q_{P-n}}\rightarrow
\~F_{n-2P,\bar n}\stackrel{Q_{n}}\rightarrow\cdots$$ 
As a consequence, 
 it was claimed that 
the $\uq$ irreducible highest
weight representation 
$V(\lambda_{\ell_{n,\bar n}})$ with 
$1\leq n\leq P$, $0\leq \bar n \leq P^\prime$ 
is isomorphic to the single 
nonvanishing BRST cohomology group as follows:
\bac
&&{\rm Ker}Q^{[s]}_n/{\rm Im}Q^{[s-1]}_n\sim 
V(\lambda_{\ell_{n,\bar n}})
\quad {\rm if}\quad s\in \Z \backslash \{0\},\nn\\
&&{\rm Ker}Q^{[s]}_n/{\rm Im}Q^{[s-1]}_n\sim 0\quad {\rm
if}\quad s\neq 0,
\ea
where 
\bac
Q^{[2s]}_n&=&Q_n:\~F_{n-2sP,\bar n}\rightarrow 
\~F_{-n-2sP,\bar n},\quad s\in \Z,\nn\\
Q^{[2s+1]}_n&=&Q_{P-n}:\~F_{-n-2sP,\bar n}\rightarrow 
\~F_{-n-2(s+1)P,\bar n}, \quad s\in \Z.
\ea
This is a natural generalization of the classical (i.e., $q=1$)
theorem. 
From the Lefschetz formula  the trace of 
a BRST invariant operator ${\cal O}$ over
$V(\lambda_{\ell_{n,\bar n}})$ follows as:
\be
Tr_{V(\lambda_{\ell_{n,\bar n}})}{\cal O}=
\sum_{s\in \Z}(-1)^s Tr_{\~F^{[s]}_{n,\bar n}}{\cal O}^{[s]}
\label{vvo5}\ee
where 
\bac
\~F^{[s]}_{n,\bar n}&=&\~F_{n-Ps,\bar n},\quad s\in 2\Z,\nn\\
\~F^{[s]}_{n,\bar n}&=&\~F_{-n-P(s-1),\bar n},\quad 
s\in 2\Z+1,
\label{vvo6}\ea
 and the graded BRST invariant 
operators ${\cal 
O}^{[s]}$ are obtained recursively from ${\cal
O}$ through the relations
\be
Q^{[s]}_n {\cal O}^{[s]}={\cal O}^{[s+1]}
Q^{[s]}_n,\quad {\cal O}^{[0]}={\cal O}.
\label{vvo7}
\ee
In particular, using (\ref{vo13}) and the fact 
\be
\mu^{-d}S_{\delta,0}=S_{\delta,0}\mu^{-d},
\ee
one can show that 
the screened vertex operators are BRST invariant 
and satisfy
\bac
\mu^{-d}\nu^{-\alpha_0}Q_{n_2}\~\Phi^{(r)}_{\ell,m}(z)&=&{\cal
A}_{\~\Phi,\ell}^{n_2}
\mu^{-d}\nu^{-\alpha_0}\~\Phi^{(r+n_3-n_1)}_{\ell,m}(z)Q_{n_1},
\nn\\
Q_{P-n_2}\mu^{-d}\nu^{-\alpha_0}\~\Phi^{(r+n_3-n_1)}_{\ell,m}(z)
&=&{\cal A}_{\~\Phi,\ell}^{P-n_2}
\mu^{-d}\nu^{-\alpha_0}\~\Phi^{(r)}_{\ell,m}(z)Q_{P-n_1},
\label{vvo8}
\ea
where ${\cal A}_{\~\Phi,\ell}$
is  defined by
\bac
S(t)\~\Phi^{(\ell)}_{\ell}(z)&=&{\cal A}_{\~\Phi,\ell}(z/t)
\~\Phi^{(\ell)}_{\ell}(z)S(t),\nn\\
{\cal A}_{\~\Phi,\ell}(z)&=&(q^{k+2}z)^{\ell\over k+2} 
{\theta_p(q^{-\ell}z)\over \theta_p(q^{\ell}z)},\nn\\
{\cal A}_{\~\Phi,\ell}&=&
{\cal A}_{\Phi,\ell}(q^{\ell+\epsilon}).
\ea
Here,  
$\epsilon$ is a
regularization parameter and  
${\cal A}_{\~\Phi,\ell}(z)$ a pseudo constant, that is,
${\cal A}_{\~\Phi,\ell}(zp^n)=
{\cal A}_{\~\Phi,\ell}(z)$, $n\in 
\Z$. 
Similarly,  physical type II vertex operators are 
 BRST invariant and satisfy
\bac
\mu^{-d}\nu^{-\alpha_0}Q_{n_2}
\~\Psi^{(r)}_{1,m}(z)&=&{\cal
A}_{\~\Psi,1}^{n_2}
\mu^{-d}\nu^{-\alpha_0}\~\Psi^{(1-r)}_{1,m}(z)Q_{n_1},\quad 
n_2-n_1=1-2r,\quad r=0,1;\nn\\
Q_{P-n_2}\mu^{-d}\nu^{-\alpha_0}
\~\Psi^{(1-r)}_{1,m}(z)
&=&{\cal A}_{\~\Psi,1}^{P-n_2}
\mu^{-d}\nu^{-\alpha_0}
\~\Psi^{(r)}_{1,m}(z)Q_{P-n_1},\quad n_2-n_1=1-2r,
\quad r=0,1;\nn\\ 
\label{vvo10}
\ea
where 
${\cal A}_{\~\Psi,1}$
is  given by the following relations:
\bac
S(t)\~\Psi^{(1)}_{1}(z)&=&{\cal A}_{\~\Psi,1}(z/t)
\~\Psi^{(1)}_{1}(z)S(t),\nn\\
{\cal A}_{\~\Psi,1}(z)&=&(q^{-2}z)^{1\over k+2} 
{\theta_p(q^{k+1}z)\over \theta_p(q^{k-1}z)},\nn\\
{\cal
A}_{\~\Psi,1}&=&{\cal
A}_{\~\Psi,1}(q^{-k-1+\epsilon}). 
\ea
Here also, $\epsilon$ is a regularization parameter and 
${\cal A}_{\~\Psi,1}(z)$ is a pseudo constant, i.e., 
${\cal A}_{\~\Psi,1}(zp^n)={\cal A}_{\~\Psi,1}(z)$, $n\in \Z$. 

Now let us make some remarks about 
the explicit  evaluation of the traces over
$\~F^{[s]}_{n,\bar n}$. The bosonization of 
the currents is independent of the zero mode $\xi_0$ since
\be
E^+(z)=_{1}\partial_{z}\xi(z)\exp{-X_2(2|q^{-k-2}z;1)}.
\ee
Therefore, states with $\xi_0$ must not be created 
by  vertex operators and propagate. This means that the
traces over the restricted Fock spaces 
of the BRST invariant operators can be evaluated in terms of
the simpler traces over the unrestricted Fock spaces as follows:
\be
Tr_{\~F^{[s]}_{n,\bar n}} {\cal O}^{[s]}=
Tr_{F^{[s]}_{n,\bar n}}\left({\cal O}^{[s]} \oint
{dy \over 2 \pi i}\eta(y) \xi(y^\prime)\right).
\label{Resf}\ee   
Let us remind that  the restricted Fock spaces 
is in the kernel of $a_{X_1,0}+
a_{X_2,0}$. The above expression (\ref{Resf}) 
is slightly different from the proposed
one in Ref.  \cite{Kon94} in two aspects: 
First, we have displaced ${\cal O}^{[s]}$ to the 
left of  $\eta(y)$ 
and $\xi(y^\prime)$. Second, We have displaced 
$\xi(y^\prime)$ to
the right of $\eta(y)$. 
The reason for
the first choice is that by doing so we arrange for the smallest
number of poles from the OPE's to contribute to the $y$
contour integral.
This is particularly useful in the limit $q\rightarrow 0$,
where the leading term of the trace reduces to that of the matrix
element. The latter will be easily computed since only the pole
at $y=y^\prime$  contributes   to the $y$ integral. 
The reason for
the second modification is merely due to 
simplification purposes of
the trace formula over the $X_2$ and $X_3$ parts of the vertex
operators. 

Later we will illustrate explicitly  these statements
through the simplest example of the 1-point  correlation 
function of the spin-1/2 
Heisenberg model. Obviously, this is already known through
Baxter's methods, 
 but we will
rederive its leading term through this general bosonization 
method, which is applicable to any spin. 

\newpage

\section{\bf Correlation functions and form factors}

\subsection{\bf N-point  correlation functions }

Here, we build on the previous chapters to evaluate 
 the  N-point
 correlation functions of local operators. 
By this, we mean their vacuum to vacuum expectation values.
We attempt to make this evaluation as explicit as possible 
but there are some expressions which can be written in
complete explicit forms only once the value of the local 
spin ($=k/2$) is fixed to some number. 

Since a general  local operator $L$ acting  on 
$V^{(k)\otimes n}$, which
is embedded in $V^{(k)\otimes \infty}$ can be expressed as a 
linear combination of  local operators of the form 
$E_{i_nj_n}\otimes \cdots \otimes E_{i_1j_1}$, it is sufficient 
to consider just the latter basis operators. Here
$E_{ij}$ are the  $k\times k$ unit matrices.  

As (\ref{vvo2}) shows, a single basis local operators has
many   correlation functions depending on the  sectors
it is acting in.
Each sector is characterized by the boundary conditions imposed 
on both ends of the chain of spins. These boundary conditions, 
in turn, are in one to one correspondence with 
the highest weights
$\lambda$, that is, with the vacuum states $|vac>_\lambda$. 
Therefore, the correlation function of a basis operator 
in the sector $\lambda_m=(k-m)\Lambda_0+m\Lambda_1$, is given by
\bac
&&<E_{i_nj_n}\otimes \cdots \otimes 
E_{i_1j_1}>^{(\lambda_m)}_{z_n,\cdots,z_1}\nn\\
&&=A(q){{\rm Tr}_{V(\lambda_m)}(\mu^{-d }
\nu^{-\alpha_0} \vt_{\l^{(1)}
,V,i_1}^{\l_m}(z_1) 
\cdots \cdots \vt_{\l^{(n)} ,V,i_n}^{\l^{(n-1)}}(z_n) 
\vt_{\l^{(n-1)},j_n}^{\l^{(n)},V}(z_n)\cdots\cdots
\vt_{\l_m,j_1}^{\l^{(1)} ,V}(z_1) )\over 
{\rm Tr}_{V(\l_m)}(\mu^{-d }\nu^{-\alpha_0})}\nn\\
\ea
where we have used (\ref{vvo1}), (\ref{vvo2}),  
$\rho=2d+\alpha_0/2$,
and
\bac
\mu&=&q^4,\nn\\
\nu&=&q,\nn\\
A(q)&=&{1\over g_{\lambda_m} g_{\lambda^{(1)}}\cdots\cdots
g_{\lambda^{(n-1)}}},\nn\\
\lambda^{(1)}&=&\sigma(\lambda_m),\nn\\
\sigma(\l^{(l)})&=&\l^{(l-1)}.
\label{co3}
\ea
Using (\ref{vvo3}), we rewrite 
the above correlation function as 
\be
<E_{i_nj_n}\otimes \cdots \otimes 
E_{i_1j_1}>^{(\lambda_m)}_{z_n,\cdots,z_1}
=B(q){T^{(\lambda_m)}_{j_{2n},\cdots, 
j_1}(z_{2n},\cdots,z_1|\mu,\nu)\over 
{\rm Tr}_{V(\l_m)}(\mu^{-d }\nu^{-\alpha_0})}, 
\ee
where 
\bac
&&T^{(\lambda_m)}_{j_{2n},\cdots, 
j_1}(z_{2n},\cdots,z_1|\mu,\nu)=
{\rm Tr}_{V(\lambda_m)}(\mu^{-d }\nu^{-\alpha_0}  
\vt_{\l^{(2n-1)},j_{2n}}^{\l_m,V}(z_{2n})\cdots
\vt_{\l_m,j_1}^{\l^{(1)} ,V}(z_1) ),\nn\\
&&\lambda^{(1)}=\sigma(\lambda_m),\nn\\
&&\l^{(n+l)}=\l^{(n-l)},\quad 1\leq l\leq n,\nn\\
&&z_{n+l}=q^{-2}z_{n-l+1},\quad 1\leq l\leq n,\nn\\
&&j_{n+l}=k-i_{n-l+1},\quad 1\leq l\leq n,\nn\\
&&B(q)=A(q){\prod_{l=1}^n {\cal C}^{(k)}_{+,k-i_l}\over 
\prod_{s=1}^n {\cal C}^{(k)}_{+,s}},
\quad {\cal C}^{(k)}_{+,s}={\cal C}^{(k)}_{+,k-m} 
\>\>\>{\rm if}\> \>\>s\in
2\N;\quad {\cal C}^{(k)}_{+,s}={\cal C}^{(k)}_{+,m} 
\>\>\>{\rm if}\>\>\> s\in
2\N+1.\nn\\
\label{co1}\ea

The evaluation of this correlation functions reduces then 
to the evaluations of the traces 
$T^{(\lambda_m)}_{j_{2n},\cdots 
j_1}(z_{2n},\cdots,z_1|\mu,\nu)$ and  
${\rm Tr}_{V(\l_m)}(\mu^{-d }\nu^{-\alpha_0})$.
Let us first focus on the calculation
of the latter one, which is just the character of 
$V(\lambda_m)$. 
From  (\ref{vvo5})-(\ref{vvo7}), this can be expressed in 
general as
\be
{\rm Tr}_{V(\Lambda_m)}(\mu^{-d}\nu^{-\alpha_0})=
\sum_{s\in \Z}(-1)^s {\rm Tr}_{\~F^{[s]}_{\~n,\bar n}}
{\cal O}^{[s]},
\ee
where  
\bac
&&{\cal O}^{[s]}=\mu^{-d}\nu^{-\alpha_0},\nn\\
&&{\rm Tr}_{\~F^{[s]}_{\~n,\bar n}}{\cal O}^{[s]}
={\rm Tr}_{F^{[s]}_{\~n,\bar n}}(\mu^{-d}\nu^{-\alpha_0}
\oint {dy\over 2\pi i}\eta(y)\xi(y^\prime))\nn\\
&&={\rm Tr}_{F^{X_1,[s]}_{\~n,\bar n}}
\left(\mu^{-d^{X_1}}\nu^{-a_{X_1,0}}\right)
{\rm Tr}_{F^{X_2,X_3,[s]}_{\~n,\bar n}}
\left(\mu^{-d^{X_2}-d^{X_3}}\nu^{-a_{X_2,0}}
\oint {dy\over 2\pi i}\eta(y)\xi(y^\prime)\right).
\ea
Here 
\bac
d^{X_1}&=&-\sum_{n>0}{n^2a_{X_1,-n}a_{X_1,n}\over 
[2n][(k+2)n]}-{a_{X_1,0}(a_{X_1,0}+2)\over 4(k+2)},\nn\\
d^{X_2}&=&\sum_{n>0}{n^2a_{X_2,-n}a_{X_2,n}\over 
[2n]^2}+{a_{X_2,0}(a_{X_2,0}-2)\over 8},\nn\\
d^{X_3}&=&
-\sum_{n>0}{n^2a_{X_3,-n}a_{X_3,n}\over 
[2n]^2}-{a_{X_3,0}(a_{X_3,0}+2)\over 8}.
\ea
Let us note that our physical situation corresponds to
$m=\~n-1$ and $\bar n=0$.

Using  (\ref{vvo6}), the general 
expression of a trace of a product of vertex operators 
 (\ref{ap5}) and the  Fourier double
transformation formula (\ref{ap6}) given in the appendix 
we find 
\bac
&&{\rm Tr}_{F^{X_1,[s]}_{\~n,\bar n}}(\mu^{-d^{X_1}}
\nu^{-a_{X_1,0}})={\nu\mu^{1\over 24}
\over \eta(\mu)}\mu^{{PP^\prime\over 4}(s-
{\~nP^\prime-\bar n P\over PP^\prime})^2}
\nu^{P(s-{\~nP^\prime-\bar n P\over PP^\prime})},\quad s\in 2\Z,
\nn\\
&&{\rm Tr}_{F^{X_1,[s]}_{\~n,\bar n}}(\mu^{-d^{X_1}}
\nu^{-a_{X_1,0}})={\nu\mu^{1\over 24}
\over \eta(\mu)}\mu^{{PP^\prime\over 4}(s-1
+{\~nP^\prime+\bar n P\over PP^\prime})^2}\nu^{P
(s-1+{\~nP^\prime+\bar n P\over PP^\prime})},\quad s\in 2\Z+1,
\nn\\
&&{\rm Tr}_{F^{X_2,X_3,[s]}_{\~n,\bar n}}
(\mu^{-d^{X_2}-d^{X_3}}\nu^{-a_{X_2,0}}
\oint {dy\over 2\pi i}\eta(y)\xi(y^\prime))
={q^3\nu^{-2} \mu^{-{1\over 24}}\eta(\mu)
\over\theta_\mu(\nu^{-2})}.
\ea
Here
the eta function $\eta$ and the Jacobi theta function 
$\theta_p(z)$ are given in the appendix (\ref{ap7}) 
and (\ref{ap8}). 
Putting all  contributions together one arrives at \cite{Kon94} 
\be
{\rm Tr}_{V(\lambda_m)}(\mu^{-d}\nu^{-\alpha_0})=
{q^3\nu^{-1}\over \theta_\mu(\nu^{-2})}\sum_{s\in \Z}
(\mu^{PP^\prime(s-
{\~nP^\prime-\bar n P\over 2PP^\prime})^2}\nu^{2P
(s-{\~nP^\prime-\bar n P\over 2PP^\prime})}-
\mu^{PP^\prime(s+
{\~nP^\prime+\bar n P\over 2PP^\prime})^2}\nu^{2P
(s+{\~nP^\prime+\bar n P\over 2PP^\prime})}).
\ee
Up to an overall  factor this expression  coincides  
with the
Weyl-Kac character formula for  the $\uq$ irreducible
highest weight module.

Let us now turn to the evaluation of the more complicated 
$T^{(\lambda_m)}_{j_{2n},\cdots 
j_1}(z_{2n},\cdots,z_1|\mu,\nu)$ trace (\ref{co1}). 
For this purpose, 
we introduce a useful c-function $f$ which will
be a building block for most of the expressions that we will
encounter henceforth. 
Let $A=:\prod_{i=1}^n A_i:$ 
and $B=:\prod_{j=1}^m B_j:$ be two normal ordered
composite operators expressed as products of normal
ordered elementary operators $A_i, B_j$, respectively. 
By elementary operator
we mean any operator of type $E^\pm_\epsilon(z)$, $S_\delta(z)$, 
$\~\Phi_{k}^{(k)}(z)$,   $\~\Psi_{1}^{(1)}(z)$, $\eta(z)$, or
$\xi(z)$.
The function $f(A,B)$ is then simply defined by
\be
:A::B:=:AB:f(A,B)=:AB:\prod_{i=1}^n\prod_{j=1}^m f(A_i,B_j).
\ee

Due to  (\ref{vvo4}),
the above trace can be re-expressed as 
\be
T^{(\lambda_m)}_{j_{2n},\cdots 
j_1}(z_{2n},\cdots,z_1|\mu,\nu)=C_{j_{2n},\cdots, 
j_1}(z_{2n},\cdots,z_1)F^{(\lambda_m)}_{j_{2n},\cdots, 
j_1}(z_{2n},\cdots,z_1|\mu,\nu),
\label{co5}
\ee
with 
\bac
&&F^{(\lambda_m)}_{j_{2n},\cdots, 
j_1}(z_{2n},\cdots,z_1|\mu,\nu)=
{\rm Tr}_{V(\lambda_m)}(\mu^{-d }\nu^{-\alpha_0}  
\~\Phi^{(r_{2n})}_{k,j_{2n}}(z_{2n})\cdots
\~\Phi^{(r_1)}_{k,j_1}(z_1) ),\nn\\
&&C_{j_{2n},\cdots, 
j_1}(z_{2n},\cdots,z_1)=
g_{\lambda^{(2n-1)}}^{\lambda_m,V}(z_{2n})\cdots
g_{\lambda_m}^{\lambda^{(1)},V}(z_{1}),\nn\\
&&r_l=k-m\>\>\>{\rm if}\>\>\> l\in 2\N;\quad 
r_l=m\>\>\>{\rm if}\>\>\> l\in 2\N+1.
\ea
Now the problem reduces further to the calculation  
of the normalization functions 
$g_{\lambda^{(i-1)}}^{\lambda^{(i)},V}(z_{i})$ and the trace
$F^{(\lambda_m)}_{j_{2n},\cdots, 
j_1}(z_{2n},\cdots,z_1|\mu,\nu)$. To simplify the solution 
to this problem 
 we use the table of OPE's (\ref{ap9})-(\ref{ap10}) 
given in the appendix to
 re-write each component     
$\~\Phi^{(r_{l})}_{k,j_{l}}(z_{l})$  as
\bac
&&\~\Phi^{(r_{l})}_{k,j_{l}}(z_{l})=
\sum_{\delta_1^{(l)},\cdots,\delta_{r_l}^{(l)}}
\sum_{\epsilon_1^{(l)},\cdots,\epsilon_{k-j_l}^{(l)}}
\int_{q^k pz}^{q^kz\infty} d_{p}t_1^{(l)}\cdots  
\int_{q^k pz}^{q^kz\infty} d_{p}t_{r_l}^{(l)} 
\oint {d\xi_1^{(l)}\over 2\pi i}
\cdots \oint {d\xi_{k-j_l}^{(l)}\over 2\pi i}\nn\\
&&\times
H_{\delta_1^{(l)},\cdots,\delta_{r_l}^{(l)},\epsilon_1^{(l)},
\cdots, \epsilon_{k-j_l}^{(l)}}
(z_l^{(l)},t_1^{(l)},\cdots,t_{r_l}^{(l)},\xi_1^{(l)},
\cdots,\xi_{k-j_l}^{(l)},
{\cal C}_l)D_l,
\label{co6}
\ea
where 
\be
D_l=:\~\Phi_{k}^{(k)}(z_{l})
S_{\delta_1^{(l)}}(t_1^{(l)})\cdots
S_{\delta_{r_l}^{(l)}}(t_{r_l}^{(l)})
E^{-}_{\epsilon_{1}^{(l)}}(\xi_1^{(l)})\cdots 
E^{-}_{\epsilon_{k-j_l}^{(l)}}(\xi_{k-j_l}^{(l)}):,
\ee
and $H_{\delta_1^{(l)},\cdots,\delta_{r_l}^{(l)},
\epsilon_1^{(l)},
\cdots, \epsilon_{k-j_l}^{(l)}}
(z_l^{(l)},t_1^{(l)},\cdots,t_{r_l}^{(l)},\xi_1^{(l)},
\cdots,\xi_{k-j_l}^{(l)},
{\cal C}_l)$ is a c-function which can be obtained from 
(\ref{vo10}) and the OPE's given in the Appendix. Moreover, 
${\cal C}_l$ defines the contours of the above integrals
depending on the domains of convergence of the various OPE's
involved in the relation (\ref{vo10}). 
Unfortunately, it does not seem
possible  to write down a simple closed formula for
this latter function. However, once the level $k$ 
and the vacuum parameter $m$ are fixed numerically 
then it becomes
relatively simpler to write down such functions. Therefore,
we assume that they are known functions with the help
of the OPE table. 
Taking this into account and relation (\ref{vo12}) 
the functions 
$g_{\lambda^{(i-1)}}^{\lambda^{(i)},V}(z_{i})$ are then 
determined
from the normalizations
\bac
&&g_{\lambda^{(i-1)}}^{\lambda^{(i)},V}(z_{i}) 
\~\Phi^{(m)}_{k,k-m}(z_i)|m;0,0>=|k-m;0,0>\otimes v_{k-m}^{(k)}+
\cdots,\quad i=1,3,\cdots 2n-1,\nn\\
&&g_{\lambda^{(i-1)}}^{\lambda^{(i)},V}(z_{i}) 
\~\Phi^{(k-m)}_{k,m}(z_i)|k-m;0,0>=|m;0,0>\otimes v_{m}^{(k)}+
\cdots,\quad i=2,4,\cdots,2n.
\ea
Using the OPE's table given in the Appendix we find
\bac
g_{\lambda^{(i-1)}}^{\lambda^{(i)},V}(z_{i})&=& h(m)^{-1},
\quad i=1,3,\cdots 2n-1,\nn\\
g_{\lambda^{(i-1)}}^{\lambda^{(i)},V}(z_{i})&=& h(k-m)^{-1},
\quad i=2,4,\cdots 2n,
\ea 
with
\bac
h(m)&=& \left(
\sum_{\delta_1,\cdots,\delta_{m}}
\sum_{\epsilon_1,\cdots,\epsilon_{m}}
\int_{q^k pz}^{q^kz\infty} d_{p}t_1\cdots  
\int_{q^k pz}^{q^kz\infty} d_{p}t_{m} 
\oint {d\xi_1\over 2\pi i}
\cdots \oint {d\xi_{m}\over 2\pi i}\right.\nn\\
&&\times
\left. H_{\delta_1,\cdots,\delta_{m},\epsilon_1,
\cdots, \epsilon_{m}}
(z_i,t_1,\cdots,t_{m},\xi_1,
\cdots,\xi_{m},
{\cal C}_i) 
(z_iq^k)^{mk\over 2(k+2)} q^{m\sum_{j=1}^m \epsilon_j}
\prod_{j=1}^m (t_j q^{-2})^{-{m\over k+2}} \right).\nn\\
\ea

Now the trace 
$F^{(\lambda_m)}_{j_{2n},\cdots, 
j_1}(z_{2n},\cdots,z_1|\mu,\nu)$ reduces to 
\bac
&&F^{(\lambda_m)}_{j_{2n},\cdots,j_1}(z_{2n},\cdots,z_1|\mu,\nu)
\nn\\
&&=\prod_{l=1}^{2n}
\left(\sum_{\delta_1^{(l)},\cdots,\delta_{r_l}^{(l)}}
\sum_{\epsilon_1^{(l)},\cdots,\epsilon_{k-j_l}^{(l)}}
\int_{q^k pz}^{q^kz\infty} d_{p}t_1^{(l)}\cdots  
\int_{q^k pz}^{q^kz\infty} d_{p}t_{r_l}^{(l)} 
\oint {d\xi_1^{(l)}\over 2\pi i}
\cdots \oint {d\xi_{k-j_l}^{(l)}\over 2\pi i}\right. \nn\\ 
&&\left. \times H_{\delta_1^{(l)},\cdots,\delta_{r_l}^{(l)},
\epsilon_1^{(l)},
\cdots, \epsilon_{k-j_l}^{(l)}}
(z_l^{(l)},t_1^{(l)},\cdots,t_{r_l}^{(l)},\xi_1^{(l)},
\cdots,\xi_{k-j_l}^{(l)},
{\cal C}_l)\right)\prod_{2n\geq l_1>l_2\geq
1}f(D_{l_1},D_{l_2})\nn\\
&&\times {\rm Tr}_{V(\lambda_m)}(
\mu^{-d }\nu^{-\alpha_0} {\cal
O}_{r_1,\cdots,r_{2n}}),
\ea
with ${\cal O}_{r_1,\cdots,r_{2n}}$ being defined by  
\be
{\cal O}_{r_1,\cdots,r_{2n}}=
:\prod_{l=1}^{2n}\~\Phi_{k}^{(k)}(z_{l})
S_{\delta_1^{(l)}}(t_1^{(l)})\cdots
S_{\delta_{r_l}^{(l)}}(t_{r_l}^{(l)})
E^{-}_{\epsilon_{1}^{(l)}}(\xi_1^{(l)})\cdots 
E^{-}_{\epsilon_{k-j_l}^{(l)}}(\xi_{k-j_l}^{(l)}):.
\ee
The function $f(D_{l_1},D_{l_2})$ can  be read off 
from the OPE's (\ref{ap9})-(\ref{ap10}) 
of elementary operators given in the appendix.

The problem then reduces further to 
the calculation of the trace ${\rm Tr}_{V(\lambda_m)}
(\mu^{-d }\nu^{-\alpha_0} {\cal O}_{r_{2n},\cdots,r_{1}})$. 
To this end, we 
apply the method  of the previous section, that is,
\be
{\rm Tr}_{V(\lambda_m)}
(\mu^{-d }\nu^{-\alpha_0} {\cal O}_{r_1,\cdots,r_{2n}})= 
\sum_{s\in\Z}(-1)^s {\rm Tr}_{F^{[s]}_{\~n,\bar n}}
\left(\mu^{-d }\nu^{-\alpha_0}
{\cal O}_{r_1,\cdots,r_{2n}}^{[s]} \oint
{dy \over 2 \pi i}\eta(y) \xi(y^\prime)\right).
\ee  
Using (\ref{vvo8}), we define the  ${\cal O}^{[s]}$ as
\bac
{\cal O}_{r_1,\cdots,r_{2n}}^{[0]}&=&{\cal O}_{r_{2n},
\cdots,r_{1}},\nn\\
{\cal O}_{r_1,\cdots,r_{2n}}^{[s]}&=&{\cal A}^{nsp}_
{\~\Phi,k}{\cal O}_{r_{2n},\cdots,r_{1}},
\quad {\rm if}\>\>\> s\in 2\Z,\nn\\
{\cal O}_{r_1,\cdots,r_{2n}}^{[s]}&=&
{\cal A}^{n(k+sp-p)}_{\~\Phi,k}{\cal
O}_{k-r_{2n},\cdots,k-r_{1}},\quad {\rm if}\>\>\> s\in 2\Z+1.
\label{co2}\ea
This means that it is sufficient 
to calculate the following trace for a fixed even number $s$: 
\be
R(s,r_{2n},\cdots,r_{1})=
Tr_{F^{[s]}_{n,n^\prime}}
\left(\mu^{-d }\nu^{-\alpha_0}{\cal O}_{r_{2n},\cdots,r_{1}} 
\oint
{dy \over 2 \pi i}\eta(y) \xi(y^\prime)\right).
\ee
The modes of $X_2$ and $X_3$ satisfy
the same Heisenberg algebra but with opposite signatures. 
Therefore one expects many simplifications if they are treated
simultaneously. For this reason, it is convenient to split the
trace into the product of two parts. One involves only the
modes of $X_1$ and the other involves only those of 
$X_2$ and $X_3$, i.e.,
\be
R(s,r_{2n},\cdots,r_{1})=R^{X_1}(s,r_{2n},\cdots,r_{1})
R^{X_2,X_3}(r_{2n},\cdots,r_{1})
\ee
where
\bac
R^{X_1}(s;r_{2n},\cdots,r_{1})&=&
{\rm Tr}_{F^{X_1,[s]}_{\~n,\bar n}}
\left(\mu^{-d^{X_1}}\nu^{-a_{X_1,0}}
{\cal O}^{X_1}_{r_{2n},
\cdots,r_{1}}\right),\nn\\
R^{X_2,X_3}(r_{2n},\cdots,r_{1})&=&
{\rm Tr}_{F^{X_2,X_3,[s]}_{\~n,\bar n}}
\left(\mu^{-d^{X_2}-d^{X_3}
}\nu^{-a_{X_2,0}}
{\cal O}^{X_2,X_3}_{r_{2n},\cdots,r_{1}} \oint
{dy \over 2 \pi i}\eta(y) \xi(y^\prime)\right),\nn\\  
\ea
with
\bac
{\cal O}^{X_1}_{r_1,\cdots,r_{2n}}&=&
:\prod_{l=1}^{2n}
\~\Phi_{k}^{(k)}(z_{l})
S^{X_1}_{\delta_1^{(l)}}(t_1^{(l)})\cdots
S^{X_1}_{\delta_{r_l}^{(l)}}(t_{r_l}^{(l)})
E^{-,X_1}_{\epsilon_{1}^{(l)}}(\xi_1^{(l)})\cdots 
E^{-,X_1}_{\epsilon_{k-j_l}^{(l)}}(\xi_{k-j_l}^{(l)}):,\nn\\
{\cal O}^{X_2,X_3}_{r_1,\cdots,r_{2n}}&=&
:\prod_{l=1}^{2n}
S^{X_2,X_3}_{\delta_1^{(l)}}(t_1^{(l)})\cdots
S^{X_2,X_3}_{\delta_{r_l}^{(l)}}(t_{r_l}^{(l)})
E^{-,X_2,X_3}_{\epsilon_{1}^{(l)}}(\xi_1^{(l)})\cdots 
E^{-,X_2,X_3}_{\epsilon_{k-j_l}^{(l)}}(\xi_{k-j_l}^{(l)}):.\nn\\
\ea
Here we have introduced the notations
\bac
&&E^{-,X_1}_\epsilon(\xi)=\exp(\partial X_1^{(\epsilon)}
(q^{-2}\xi;-{k+2\over 2})),\nn\\
&&E^{-,X_2,X_3}_\epsilon(\xi)=
\exp(X_2(2|q^{(\epsilon-1)(k+2)}\xi;-1)+
X_3(2|q^{(\epsilon-1)(k+1)-1}\xi;0)),\nn\\
&&S_\delta^{X_1}(t)=\exp(-X_1(k+2|q^{-2}t;-{k+2\over 2})),\nn\\
&&S_\delta^{X_2,X_3}(t)=
\exp(-X_2(2|q^{-k-2}t;-1)-X_3(2|q^{-k-2+\delta}t;0)).
\ea
Using the general traces given in 
relations (\ref{ap4}) and (\ref{ap5}) of the 
appendix, we find then 
\bac
R^{X_1}(s;r_{2n},\cdots,r_{1})&=&
T_{\ell_s}^{X_1}(2n,0,N,M,n_4|X,Y,Z,W),\nn\\
R^{X_2,X_3}(r_{2n},\cdots,r_{1})&=&
T^{X_2,X_3}(A,B,C,D,G,H|\~X,\~Y,\~Z,\~W),
\ea
where
\bac
\ell_s&=& \ell_{n-Ps,n^\prime},\nn\\
N&=&\sum_{l=1}^{2n}r_l=nk,\nn\\
M&=&2kn-\sum_{l=1}^{2n}j_{l}=nk+\sum_{l=1}^n i_l-
\sum_{l=1}^n j_l,\nn\\
t_{a+\sum_{i=1}^{l-1}r_i}&=&t^{(l)}_a,\quad 1\leq a\leq r_l,
\nn\\
\xi_{a+\sum_{i=1}^{l-1}(k-j_i)}&=&\xi^{(l)}_a,\quad 1\leq a\leq 
k-j_l,\nn\\
X&=&\{z_1,\cdots,z_{2n}\},\nn\\
Y&=&\{\},\nn\\
Z&=&\{t_1,\cdots,t_N\},\nn\\
W&=&\{\xi_1,\cdots,\xi_M\},\nn\\
A&=&\{\delta_1-k-2,\cdots,\delta_N-k-2\},\nn\\
B&=&\{(\epsilon_1-1)(k+1)-1,\cdots,(\epsilon_M-1)(k+1)-1\},\nn\\
C&=&\{-k-2,\cdots,-k-2\},\nn\\
D&=&\{-1,\cdots,-1\},\nn\\
G&=&\{(\epsilon_1-1)(k+2),\cdots,(\epsilon_M-1)(k+2)\},\nn\\
H&=&\{-1,\cdots,-1\},\nn\\
\~X&=&\{t_1,\cdots,t_N\},\nn\\
\~Y&=&\{\xi_1,\cdots,\xi_M\},\nn\\
\~Z&=&\{t_1,\cdots,t_N\},\nn\\
\~W&=&\{\xi_1,\cdots,\xi_M\}.
\ea
Though these expressions are cumbersome they still allow us to
deduce a simple selection rule. 
Indeed, from the trace formulas in the appendix, we have 
\bac
R^{X_1}(s;r_{2n},\cdots,r_{1})&=&
T_{\ell_s}^{X_1}(2n,0,N,M,n_4|X,Y,Z,W)=
\delta_{2nk,2N}(\cdots),\nn\\
R^{X_2,X_3}(r_{2n},\cdots,r_{1})&=&
T^{X_2,X_3}(A,B,C,D,G,H|\~X,\~Y,\~Z,\~W)= \delta_{N,M}(\cdots).
\ea
The first constraint $N=nk$ is trivially satisfied, and the
second one, which less trivial reads
\be
\sum_{l=1}^n i_l=
\sum_{l=1}^n j_l.
\ee
Thus, only local operators which satisfy  the latter sum rule,
i.e., conserve the spin,  
lead to nonvanishing  correlation functions. 
This selection rule is a generalization of the
spin $1/2$ case \cite{JiMi94}.

Let us now recapitulate the explicit 
result of the latter calculation.
The  correlation function 
$<E_{i_nj_n}\otimes \cdots \otimes 
E_{i_1j_1}>^{(\lambda_m)}_{z_n,\cdots,z_1}$ 
of the basis local operator $E_{i_nj_n}\otimes \cdots \otimes 
E_{i_1j_1}$ is given by
\bac
<E_{i_nj_n}\otimes &\cdots& \otimes 
E_{i_1j_1}>^{(\lambda_m)}_{z_n,\cdots,z_1}=
{B(q)\over Tr_{V(\lambda_m)}(\mu^{-d}\nu^{-\alpha_0})}
C_{j_{2n},\cdots, 
j_1}(z_{2n},\cdots,z_1)\nn\\
&&\times
\prod_{l=1}^{2n}
\{\sum_{\delta_1^{(l)},\cdots,\delta_{r_l}^{(l)}}
\sum_{\epsilon_1^{(l)},\cdots,\epsilon_{k-j_l}^{(l)}}
\int_{q^k pz}^{q^kz\infty} 
d_{p}t_1^{(l)}\cdots  
\int_{q^k pz}^{q^kz\infty} 
d_{p}t_{r_l}^{(l)} 
\oint {d\xi_1^{(l)}\over 2\pi i}
\cdots \oint {d\xi_{k-j_l}^{(l)}\over 2\pi i}\nn\\ 
&&\times H_{\delta_1^{(l)},\cdots,\delta_{r_l}^{(l)},
\epsilon_1^{(l)},
\cdots, \epsilon_{k-j_l}^{(l)}}
(z_l^{(l)},t_1^{(l)},\cdots,t_{r_l}^{(l)},\xi_1^{(l)},
\cdots,\xi_{k-j_l}^{(l)},
{\cal C}_l)\}\prod_{2n\geq l_1>l_2\geq
1}f(D_{l_1},D_{l_2})\nn\\
&&
\times\left(\sum_{s\in 2\Z} {\cal A}^{nsp}_{\~\Phi,k} 
R(s,r_{2n},\cdots, r_1)
-\sum_{s\in 2\Z+1} {\cal A}^{n(k+sp-p)}_{\~\Phi,k}
R(s,k-r_{2n},\cdots, k-r_1)\right),\nn\\
\ea
where all the various symbols are defined throughout the 
preceeding analysis and all the $z_i$ are set  to 1.
Moreover, for this correlation function to be non-vanishing 
 the selection rule $\sum_{l=1}^n i_l=
\sum_{l=1}^n j_l$ must be fulfilled.

Let us now briefly check, 
in the simplest case  $s=1/2$
that the leading order term of the 1-point correlation 
 function coincides
with the known result. 
More specifically, we consider the 
1-point correlation  function of the local operator $E_{00}$ in 
the sector $\Lambda_0$. This means that this case corresponds to
the specializations 
\bac
k&=&1,\nn\\
P&=&k+2=3,\nn\\
P^\prime&=&1,\nn\\
\~n&=&1,\nn\\
\bar n&=&0.
\ea

Following the general prescription 
this correlation function is defined by 
\be
<E_{00}>_{\Lambda_0}=P_0(q^{-2}z,z),
\ee
where 
\bac
P_0(z,w)&=&{1\over g_{\Lambda_0}}
{{\rm
Tr}_{V(\Lambda_0)}(\mu^{-d}\nu^{-\alpha_0}\~\Phi^{\Lambda_0}_{
\Lambda_1,V,0}(zq^2) \~\Phi^{\Lambda_1,V}_{\Lambda_0,0}(w))\over
{\rm Tr}_{V(\Lambda_0)}(\mu^{-d}\nu^{-\alpha_0})}\nn\\
&=&-{q\over g_{\Lambda_0}}
{{\rm Tr}_{V(\Lambda_0)}(\mu^{-d}\nu^{-\alpha_0}
\~\Phi^{\Lambda_0,V}_{\Lambda_1,1}(z) 
\~\Phi^{\Lambda_1,V}_{\Lambda_0,0}(w))\over
{\rm Tr}_{V(\Lambda_0)}(\mu^{-d}\nu^{-\alpha_0})}\nn\\
&=&-{q\over g_{\Lambda_0}}g_{\Lambda_1}^{\Lambda_0,V}(z)
g_{\Lambda_0}^{\Lambda_1,V}(w)
{{\rm Tr}_{V(\Lambda_0)}(\mu^{-d}\nu^{-\alpha_0}
\~\Phi^{(1)}_{1,1}(z)\~\Phi^{(0)}_{1,0}(w)
\over
{\rm Tr}_{V(\Lambda_0)}(\mu^{-d}\nu^{-\alpha_0})},
\ea
with 
\bac
\mu&=&q^4,\nn\\
\nu&=&q,\nn\\
g_{\Lambda_0}&=&{(q^4;q^4)_\infty\over (q^2;q^4)_\infty}.
\ea
Here the component vertex operators are explicitly given by
\bac
\~\Phi^{(1)}_{1,1}(z)&=&\int_{qpz}^{qz\infty}d_ptS(t)\Phi_{k}
(z),\nn\\
\~\Phi^{(0)}_{1,0}(z)&=&\oint{d\xi\over 2\pi i}[\Phi_{k}(z),
E^-(\xi)]_q,\nn\\
g_{\Lambda_1}^{\Lambda_0,V}(z)&=& {q^{-{5\over 6}}z^{1\over
2}\over B_p({2\over 3};{2\over 3})},\nn\\
g_{\Lambda_0}^{\Lambda_1,V}(z)&=&1,
\ea
where $p=q^6$ and the quantum beta function $B_p(x,y)$ is
defined in the appendix (\ref{ap11}). The
normalization constants are determined from
\bac
&&\~\Phi^{\Lambda_1,V}_{\Lambda_0}(z))|\Lambda_0>=|\Lambda_1>+
\cdots,\nn\\
&&\~\Phi^{\Lambda_0,V}_{\Lambda_1}(z))|\Lambda_1>=|\Lambda_0>+
\cdots.
\ea
The latter relations translate into 
\bac
&&g_{\Lambda_0}^{\Lambda_1,V}(z)\~\Phi^{(0)}_{1,1}(z)
|\Lambda_0>=
|\Lambda_1>+\cdots,\nn\\
&&g_{\Lambda_1}^{\Lambda_0,V}(z)\~\Phi^{(1)}_{1,0}(z)
|\Lambda_1>=|\Lambda_0>+\cdots.
\ea
From the Lefschetz formula and (\ref{co2}) we have
\be
Tr_{V(\Lambda_0)}(\mu^{-d}\nu^{-\alpha_0}
\~\Phi^{(1)}_{1,1}(z)\~\Phi^{(0)}_{1,0}(w)=
\sum_{s\in \Z}
(-1)^s 
Tr_{F^{[s]}_{\~n,\bar n}}\left({\cal O}^{[s]} \oint
{dy \over 2 \pi i}\eta(y) \xi(y^\prime)\right),
\ee
with
\bac
&&{\cal O}={\cal O}^{[0]}=\mu^{-d}\nu^{-\alpha_0}
\~\Phi^{(1)}_{1,1}(z)\~\Phi^{(0)}_{1,0}(w),\nn\\
&&{\cal O}^{[s]}={\cal A}_\ell^{1+2n_1+(s-1)P}
\mu^{-d}\nu^{-\alpha_0}\~\Phi^{(0)}_{1,1}(z)
\~\Phi^{(1)}_{1,0}(w),\quad s\in 2\Z+1,\nn\\
&&{\cal O}^{[s]}={\cal A}_\ell^{sP}
\mu^{-d}\nu^{-\alpha_0}
\~\Phi^{(1)}_{1,1}(z)\~\Phi^{(0)}_{1,0}(w),\quad s\in 2\Z.
\ea
Since the vacuum state $|\Lambda_0>$ has the lowest degree with
respect to -$d$, which 
is zero, 
clearly the leading term in $q$ (in the limit $q\rightarrow 0$) 
of the 1-point function must coincide 
with the leading term of the $s=0$ term. The latter, in turn,  
is obtained from the leading term of the following 
matrix element of vertex operators:
\bac
&&<\Lambda_0|\~\Phi^{\Lambda_0,V}_{\Lambda_1,1}(z) 
\~\Phi^{\Lambda_1,V}_{\Lambda_0,0}(w))|\Lambda_0>_{\~F}\nn\\
&=& 
g_{\Lambda_1}^{\Lambda_0,V}(z)g_{\Lambda_0}^{\Lambda_1,V}(w)
<0|\int_{qpz}^{qz\infty}d_ptS(t)\Phi_{\ell}z)
\oint{d\xi\over 2\pi i}[\Phi_{\ell},
E^-(\xi)]_q\oint {dy\over 2\pi i}\eta(y) \xi(y^\prime)|0>\nn\\
&=&-q^4{(q^6w/z;q^4)_\infty\over (q^4w/z;q^4)_\infty}.
\ea   
This means that
\be
P_0(zq^{-2},z)\sim {q^2\over
g_{\Lambda_0}}.
\ee
and therefore the leading term of the physical 1-point 
correlation function is given by
\be
P_0(q^{-2}z,z)\sim  O(q^2),
\ee
and  coincides with the known result 
\cite{Bax73,Jimal92,Idz93}. 
Note that up to a factor $q^{-3}$ originating from the OPE of
$\eta(y)$ and $\xi(y^\prime)$ this matrix element is equal 
to the matrix element of the vertex operators acting on the 
unrestricted Fock space, i.e., 
\bac
<\Lambda_0|\~\Phi^{\Lambda_0,V}_{\Lambda_1,1}(z) 
\~\Phi^{\Lambda_1,V}_{\Lambda_0,0}(w))|\Lambda_0>_{F}
&=& 
g_{\Lambda_1}^{\Lambda_0,V}(z)g_{\Lambda_0}^{\Lambda_1,V}(w)\nn\\
&&\times<0|\int_{qpz}^{qz\infty}d_ptS(t)\Phi_{1}(z)
\oint{d\xi\over 2\pi i}[\Phi_{1}(z),
E^-(\xi)]_q|0>\nn\\
&=&-q{(q^6w/z;q^4)_\infty\over (q^4w/z;q^4)_\infty}\nn\\
&=&q^{-3}<\Lambda_0|\~\Phi^{\Lambda_0,V}_{\Lambda_1,1}(z) 
\~\Phi^{\Lambda_1,V}_{\Lambda_0,0}(w))|\Lambda_0>_{\~F}.
\ea   
The reason is that the vertex operators are BRST invariant and
the highest weight states $|\Lambda_i>$ are BRST states and
therefore the above relation is  expected. 

\subsection{\bf N-point form factors of  local operators}

In this section we derive  the N-point form factors of local
operators. Although the physically interesting ones 
are those where local operators act on a single site, we
consider the most general case where local operators act on
a set of neighboring sites. 
Let us recall that the physical type II vertex operators 
are the
spin-1/2 ones. They are bosonized as 
follows: 
\be
\~\Psi_{\lambda_{\ell}}^{V^{(1)},\lambda_{\ell\pm 1}}(z)=
h_{\lambda_{\ell}}^{V^{(1)},\lambda_{\ell\pm 1}}(z)
\sum_{m=0}^1 v_m^{(1)}\otimes
\~\Psi^{(r_{\pm})}_{1,m}(z),
\label{co4}
\ee
with $r_+=0$ and $r_-=1$, and
\bac
\~\Psi^{(0)}_{1,m}(z)&=&\~\Psi^{(1)}_{m}(z),\nn\\
\~\Psi^{(1)}_{1,m}(z)&=&
\int_0^{q^{k-1}z} d_{p}t   
\~\Psi^{(1)}_{m}(z) S(t).
\ea
The bosonizations of the bare components are derived from 
(\ref{vvo9}), that is,
\bac
\~\Psi^{(0)}_{1,0}(z)&=&\~\Psi^{(1)}_{0}(z)=
\exp({X_1(1;2,k+2|q^{k-2}z;-{k+2\over 2})+
X_2(2|q^{-2}z;0)+X_3(2|q^{-2}z;0)}),\nn\\
\~\Psi^{(0)}_{1,1}(z)&=&\oint{du\over 2\pi i}
{[\~\Psi^{(1)}_{0}(z),E^+(u)]}_q,\nn\\
\~\Psi^{(1)}_{1,0}(z)&=&\int_0^{q^{k-1}z} d_p t 
\~\Psi^{(1)}_{0}(z)S(t),\nn\\
\~\Psi^{(1)}_{1,1}(z)&=&\int_0^{q^{k-1}z} 
d_p t \oint{du\over 2\pi i}
{[\~\Psi^{(1)}_{0}(z),E^+(u)]}_q S(t).
\ea 
Moreover, one can easily check that
\bac
{[\~\Psi^{(\ell)}_{\ell}(z),\eta_0]}&=&0,\nn\\
{[\~\Phi^{(\ell)}_{\ell}(z),a_{X_2,0}+a_{X_3,0}]}&=&0,
\ea
and hence the screened components intertwine the restriced 
Fock spaces as
\be
\~\Psi^{(r)}_{1,m}(z): \~F_{\ell_1}\rightarrow   
V^{(1)}(z)\otimes \~F_{\ell_2}.
\ee 
Just as in the case of Type I vertex operators, it is 
convenient to write the bosonized bare components of the
vertex operators in the following form:
\bac
\~\Psi^{(0)}_{1,0}(z)&=&
\sum_{\epsilon^\prime=\pm 1} \oint{du\over 2\pi i}
I_{\epsilon^\prime}^{(0),(1)}(z,u|
C):\~\Psi^{(1)}_{0}(z)E^+_{\epsilon^\prime}(u):\nn\\
\~\Psi^{(1)}_{1,0}(z)&=&
\sum_{\delta^\prime=\pm 1}\int_{0}^{q^{k-1}z} d_p t
I_{\delta^\prime}^{(1),(0)}(z,t|
C):\~\Psi^{(1)}_{0}(z)
S_{\delta^\prime}(t):\nn\\
\~\Psi^{(1)}_{1,1}(z)&=&
\sum_{\delta^\prime=\pm 1}\sum_{\epsilon^\prime=\pm 1}
\int_{0}^{q^{k-1}z} d_p t \oint{du\over 2\pi i}
I_{\delta^\prime,\epsilon^\prime}^{(1),(1)}(z,t,u|
C):\~\Psi^{(1)}_{0}(z)
S_{\delta^\prime}(t)E^+_{\epsilon^\prime}(u):
\label{co8}
\ea
The various c-functions $I$ and corresponding  
contours $C$ are given
in the appendix (\ref{ap2}). 

As mentioned previously 
any local operator can be written as a linear
combination of the basis local operators 
$E_{i_nj_n}\otimes \cdots\otimes E_{i_1j_1}$, and therefore 
we consider
only the following  
general form factor:
\bac
&&<\lambda_m|E_{i_nj_n}\otimes \cdots
\otimes E_{i_1j_1}|\rho_{n^\prime},\cdots,
\rho_1>_{1-\alpha_{n^\prime},\cdots,
1-\alpha_{1}}^{\beta_{n^\prime},\cdots,\beta_{1}}\nn\\
&&={{\rm Tr}_{V(\lambda_m)}(\mu^{-d}\nu^{-\alpha_0}
E_{i_nj_n}\otimes \cdots\otimes E_{i_1j_1}
\Psi^{V^*{\lambda^\prime}_{n^\prime}}_{
{\lambda^\prime}_{n^\prime-1},1-\alpha_{n^\prime}}
(\rho_{n^\prime})\cdots 
\Psi^{V^*{\lambda^\prime}_{1}}_{{\lambda^\prime}_{0},
1-\alpha_{1}}
(\rho_{1}))\over {\rm
Tr}_{V(\lambda_m)}(\mu^{-d}\nu^{-\alpha_0})}\nn\\
&&=B(q)\prod_{i=1}^{n^\prime}d_{\alpha_{i},\beta_i}
{{\rm Tr}_{V({\lambda}_m)}(\mu^{-d }\nu^{-\alpha_0}  
\vt_{\l^{(2n-1)},j_{2n}}^{\l^{(2n)},V}(z_{2n})\cdots
\vt_{\l^{(0)},j_1}^{\l^{(1)} ,V}(z_1)\Psi^{
{\lambda^\prime}_{n^\prime}}_{
{\lambda^\prime}_{n^\prime-1},\alpha_{n^\prime}}
(\bar \rho_{n^\prime})\cdots 
\Psi^{{\lambda^\prime}_{1}}_{{\lambda^\prime}_{0},\alpha_{1}}
(\bar \rho_{1}))\over
{\rm Tr}_{V({\lambda}_m)}(\mu^{-d}\nu^{-\alpha_0})},\nn\\
\ea  
where 
\bac
\Psi^{V^*{\lambda^\prime}_{i}}_{{\lambda^\prime}_{1-i},
1-\alpha_{i}}
(\rho_{i})&=&d_{\alpha_{i},\beta_i}
\Psi^{{\lambda^\prime}_{i}}_{{\lambda^\prime}_{1-i},\alpha_{i}}
(\bar \rho_{i}),\quad 1\leq i\leq n^\prime,\nn\\
d_{1,1}&=&1,\quad d_{1,-1}=-q^{-1},\quad d_{0,1}=-q,\quad 
d_{0,-1}=1,\nn\\
\bar \rho_{i}&=&q^2\rho_{i},\nn\\
\alpha_i&=&0,1;\nn\\
{\lambda^\prime}_i&=&{\lambda^\prime}_{1-i}+\beta_i(
\Lambda_1-\Lambda_0),\quad
\beta_i=\pm 1,\quad 1\leq i\leq n^\prime,\nn\\
{\lambda^\prime}_{0}&=&{\lambda^\prime}_{n^\prime}=\lambda_m,
\label{cco7}
\ea
and all other notations are defined in (\ref{co3}). 
We have already
evaluated the trace in the denominator of the above form factor.
Let us then focus on the trace in the numerator, i.e.,
\bac
T_1&=& 
{\rm Tr}_{V(\lambda_m)}(\mu^{-d }\nu^{-\alpha_0}  
\vt_{\l^{(2n-1)},j_{2n}}^{\l_m,V}(z_{2n})\cdots
\vt_{\l^{(0)},j_1}^{\l^{(1)} ,V}(z_1)\Psi^{
\lambda_{n^\prime}}_{\lambda_{n^\prime-1},\alpha_{n^\prime}}
(\bar \rho_{n^\prime})\cdots 
\Psi^{\lambda^\prime_{1}}_{\lambda^\prime_{0},\alpha_{1}}
(\bar \rho_{1}))\nn\\
&=&C_{j_{2n},\cdots,j_1}(z_{2n},\cdots,z_1)
\left( \prod_{i=1}^{n^\prime}h^{
\lambda_{i}}_{\lambda_{i-1}}
(\bar \rho_{i})\right)\times T_2
\ea
where we have used (\ref{co4}) and (\ref{co5}) 
to reduce the trace down 
to 
\be
T_2={\rm Tr}_{V(\lambda_m)}(\mu^{-d }\nu^{-\alpha_0}  
\~\Phi^{(r_{2n})}_{k,j_{2n}}(z_{2n})\cdots
\~\Phi^{(r_1)}_{k,j_1}(z_1)
\Psi^{(r^\prime_{n^\prime})}_{1,\alpha_{n^\prime}}
(\bar \rho_{n^\prime})\cdots 
\Psi^{(r^\prime_{1})}_{1,\alpha_{1}}
(\bar \rho_{1})),
\ee
where  $r^\prime_i=(1-\beta_i)/2$.
Because any order of products of vertex operators 
can be obtained from an initial one through the commutation
relations (\ref{vo7}), let us assume for simplicity that 
in the above trace, type II vertex operators 
are ordered in the
following manner:
\bac
&&\Psi^{(r^\prime_{n^\prime})}_{1,\alpha_{n^\prime}}
(\bar \rho_{n^\prime})\cdots 
\Psi^{(r^\prime_{1}}_{1,\alpha_{1})}
(\bar \rho_{1}))=
\Psi^{(0)}_{1,0}
(\bar \rho_{n^\prime})\cdots \Psi^{(0)}_{1,0}
(\bar \rho_{a+b+c+1})\nn\\
&&\times \Psi^{(0)}_{1,1}
(\bar \rho_{a+b+c})\cdots \Psi^{(0)}_{1,1}
(\bar \rho_{a+b+1})
\Psi^{(1)}_{1,0}
(\bar \rho_{a+b})\cdots \Psi^{(1)}_{1,0}
(\bar \rho_{a+1})
 \Psi^{(1)}_{1,1}
(\bar \rho_{a})\cdots \Psi^{(1)}_{1,1}
(\bar \rho_{1}).
\ea
We can then simply normal order the latter operator as
\bac
&&\Psi^{(r^\prime_{n^\prime})}_{1,\alpha_{n^\prime}}
(\bar \rho_{n^\prime})\cdots 
\Psi^{(r^\prime_{1}}_{1,\alpha_{1}}
(\bar \rho_{1})= 
\prod_{i=1}^{a}\left(
\sum_{\delta^{(a)\prime}_i,\epsilon^{(a)\prime}_i}
\int_{0}^{q^{k-1}z} d_p t^{(a)}_i \oint{du^{(a)}_i\over 2\pi i}
I_{\delta^{(a)\prime}_i,\epsilon^{(a)\prime}_i}^{(1),(1)}
(\bar \rho_i,t^{(a)}_i,u^{(a)}_i|C)\right)\nn\\
&&\times \prod_{i=a+1}^{b}\left(
\sum_{\delta^{(b)\prime}_i} 
\int_{0}^{q^{k-1}z} d_p t^{(b)}_i
I_{\delta^{(b)\prime}_i}^{(1),(0)}
(\bar \rho_i,t^{(b)}_i|C)\right)
\prod_{i=b+1}^{c}\left(
\sum_{\epsilon^{(c)\prime}_i}
\oint{du^{(c)}_i\over 2\pi i}
I_{\epsilon^{(c)\prime}_i}^{(0),(1)}
(\bar \rho_i,u^{(c)}_i|C)\right)\nn\\ 
&&\times \prod_{1\leq i^\prime
<i\leq n^\prime} f(K_i,K_{i^\prime}){\cal O}_{
r^\prime_{n^\prime},\cdots,r^\prime_1},
\ea
where we have introduced the following notations for the
normal ordered composite operators:
\bac
{\cal O}_{
r^\prime_{n^\prime},\cdots,r^\prime_1}&=&:\prod_{i=1}^{n^\prime}
K_i:,\nn\\
K_i&=&:\Psi^{(1)}_{0}(\bar \rho_{i})S_{\delta^{(a)\prime}_i}
(t^{(a)}_i)
E^+_{\epsilon^{(a)\prime}_i}(u^{(a)}_i):,
\quad 1\leq i\leq a,\nn\\
K_i&=&:\Psi^{(1)}_{0}(\bar \rho_{i})S_{\delta^{(b)\prime}_i}
(t^{(b)}_i):,
\quad a+1\leq i\leq a+b,\nn\\
K_i&=&:\Psi^{(1)}_{0}(\bar \rho_{i})
E^+_{\epsilon^{(a)\prime}_i}(u^{(c)}_i):,
\quad a+b+1\leq i\leq a+b+c,\nn\\
K_i&=&\Psi^{(1)}_{0}(\bar \rho_{i}),
\quad a+b+c+1\leq i\leq n^\prime.
\ea
Using (\ref{co6}) and the latter normal ordering, 
we express  $T_2$ as
\bac
T_2&=&
\prod_{l=1}^{2n}
\left(\sum_{\delta_1^{(l)},\cdots,\delta_{r_l}^{(l)}}
\sum_{\epsilon_1^{(l)},\cdots,\epsilon_{k-j_l}^{(l)}}
\int_{q^k pz}^{q^kz\infty} 
d_{p}t_1^{(l)}\cdots  
\int_{q^k pz}^{q^kz\infty} 
d_{p}t_{r_l}^{(l)} 
\oint {d\xi_1^{(l)}\over 2\pi i}
\cdots \oint {d\xi_{k-j_l}^{(l)}\over 2\pi i}\right.\nn\\ 
&&\times \left. H_{\delta_1^{(l)},\cdots,\delta_{r_l}^{(l)},
\epsilon_1^{(l)},
\cdots, \epsilon_{k-j_l}^{(l)}}
(z_l^{(l)},t_1^{(l)},\cdots,t_{r_l}^{(l)},\xi_1^{(l)},
\cdots,\xi_{k-j_l}^{(l)},
{\cal C}_l)\right)\prod_{2n\geq l_1>l_2\geq
1}f(D_{l_1},D_{l_2})\nn\\
&&\times \left\{
\prod_{i=1}^{a}\left(
\sum_{\delta^{(a)\prime}_i,\epsilon^{(a)\prime}_i}
\int_{0}^{q^{k-1}z} d_p t^{(a)}_i \oint{du^{(a)}_i\over 2\pi i}
I_{\delta^{(a)\prime}_i,\epsilon^{(a)\prime}_i}^{(1),(1)}
(\bar \rho_i,t^{(a)}_i,u^{(a)}_i|C)\right)\right.\nn\\
&&\times  \prod_{i=a+1}^{b}\left(
\sum_{\delta^{(b)\prime}_i} 
\int_{0}^{q^{k-1}z} d_p t^{(b)}_i
I_{\delta^{(b)\prime}_i}^{(1),(0)}
(\bar \rho_i,t^{(b)}_i|C)\right)
\prod_{i=b+1}^{c}\left(
\sum_{\epsilon^{(c)\prime}_i} \oint{du^{(c)}_i\over 2\pi i}
I_{\epsilon^{(c)\prime}_i}^{(0),(1)}
(\bar \rho_i,u^{(c)}_i|C)\right)\nn\\ 
&&\times \left. \prod_{1\leq i^\prime
<i\leq n^\prime} f(K_i,K_{i^\prime})\right\} 
f({\cal O}_{r_{2n},\cdots,r_{1}}, {\cal O}_{
r^\prime_{n^\prime},\cdots,r^\prime_1})\times T_3.
\ea
In this expression which term is the integrand of which 
integral should be clear from the various indices.  
The trace $T_3$ is given by 
\bac
T_3&=& 
{\rm Tr}_{V(\lambda_m)}(
\mu^{-d }\nu^{-\alpha_0} 
{\cal O}_{r_{2n},\cdots,r_{1};r^\prime_{n^\prime},\cdots,
r^\prime_{1}})\nn\\
&=&\sum_{s\in \Z}(-1)^s 
{\rm Tr}_{F^{[s]}_{n,n^\prime}}(\mu^{-d }\nu^{-\alpha_0}  
{\cal O}^{[s]}_{r_{2n},\cdots,r_{1};r^\prime_{n^\prime},\cdots,
r^\prime_{1}}\oint {dy\over 2\pi i} \eta(y) \xi(y^\prime)).
\ea
Here the operators 
${\cal O}^{[s]}_{r_{2n},\cdots,r_{1};r^\prime_{n^\prime},\cdots,
r^\prime_{1}}$ are
defined through (\ref{vvo8}) and (\ref{vvo10} as follows: 
\bac
{\cal O}^{[0]}_{r_{2n},\cdots,r_{1};r^\prime_{n^\prime},\cdots,
r^\prime_{1}}&=&
{\cal O}_{r_{2n},\cdots,r_{1};r^\prime_{n^\prime},\cdots,
r^\prime_{1}}=:
{\cal O}_{r_{2n},\cdots,r_{1}}
{\cal O}_{
r^\prime_{n^\prime},\cdots,r^\prime_1}:,\nn\\
{\cal O}^{[s]}_{r_{2n},\cdots,r_{1};r^\prime_{n^\prime},\cdots,
r^\prime_{1}}&=&
{\cal A}^{n(k+sp-p)}_{\~\Phi,k}
 {\cal A}_{\~\Psi,1}^{spn^\prime/2}
{\cal O}_{r_{2n},\cdots,r_{1};r^\prime_{n^\prime},\cdots,
r^\prime_{1}},\quad s\in 2\Z\nn\\
{\cal O}^{[s]}_{r_{2n},\cdots,r_{1};r^\prime_{n^\prime},\cdots,
r^\prime_{1}}&=&{\cal A}^{n(k+sp-p)}_{\~\Phi,k}
{\cal A}_{\~\Psi,1}^{
(2n-n^\prime+(s-1)p+1)n^\prime/2+
2\sum_{l=1}^{n^\prime-1}lr^\prime_{l+1}} 
{\cal O}_{1-r_{2n},\cdots,1-r_{1};1-r^\prime_{n^\prime},\cdots,
1-r^\prime_{1}},\nn\\
&&\quad s\in 2\Z+1.
\ea
Let us note that the last relation  is satisfied if this  
condition on the various parameters
$r^\prime_l$ holds:
\be
2\sum_{l=1}^{n^\prime} r_l^{\prime}=n^\prime.
\ee

Here also, in order to evaluate the above trace it is sufficient 
to calculate the following one for a fixed even number $s$: 
\be
S(s,r_{2n},\cdots,r_{1};
r^\prime_{n^\prime},\cdots,r^\prime_{1})=
Tr_{F^{[s]}_{\~n,\bar n}}
\left({\cal O}_{r_{2n},\cdots,r_{1};
r^\prime_{n^\prime},\cdots,r^\prime_{1}} \oint
{dy \over 2 \pi i}\eta(y) \xi(y^\prime)\right).
\ee
Again,  it is convenient to split the
trace into the product of two parts. One involves only the
modes of $X_1$ and the other involves only those of 
$X_2$ and $X_3$, i.e.,
\be
S(s,r_{2n},\cdots,r_{1};
r^\prime_{n^\prime},\cdots,r^\prime_{1})=
S^{X_1}(s,r_{2n},\cdots,r_{1};
r^\prime_{n^\prime},\cdots,r^\prime_{1})
S^{X_2,X_3}(r_{2n},\cdots,r_{1};
r^\prime_{n^\prime},\cdots,r^\prime_{1}),
\ee
where
\bac
S^{X_1}(s,r_{2n},\cdots,r_{1};
r^\prime_{n^\prime},\cdots,r^\prime_{1})&=&
Tr_{F^{X_1,[s]}_{\~n,\bar n}}
\left(\mu^{-d^{X_1}}\nu^{-a_{X_1,0}}
{\cal O}^{X_1}_{r_{2n},\cdots,r_{1};
r^\prime_{n^\prime},\cdots,r^\prime_{1}}\right),\nn\\
S^{X_2,X_3}(r_{2n},\cdots,r_{1};
r^\prime_{n^\prime},\cdots,r^\prime_{1})&=&
Tr_{F^{X_2,X_3,[s]}_{\~n,\bar n}}
\left(\mu^{-d^{X_2}-d^{X_3}}\nu^{-a_{X_2,0}}
{\cal O}^{X_2,X_3}_{r_{2n},\cdots,r_{1};
r^\prime_{n^\prime},\cdots,r^\prime_{1}} \oint
{dy \over 2 \pi i}\eta(y) \xi(y^\prime)\right),\nn\\  
\ea
and
\bac
{\cal O}^{X_1}_{r_{2n},\cdots,r_{1};
r^\prime_{n^\prime},\cdots,r^\prime_{1}}&=&
:{\cal O}^{X_1}_{r_{2n},\cdots,r_{1}}
\prod_{i=1}^{n^\prime}
\~\Psi^{(1),X_1}_0(\bar \rho_{i})
\times \prod_{i=1}^{a} 
S^{X_1}_{\delta^{(a)\prime}_i}
(t^{(a)}_i)
\times \prod_{i=a+1}^{b} 
S^{X_1}_{\delta^{(b)\prime}_i}
(t^{(b)}_i):,\nn\\
{\cal O}^{X_2,X_3}_{r_{2n},\cdots,r_{1};
r^\prime_{n^\prime},\cdots,r^\prime_{1}}
&=&:{\cal O}^{X_2,X_3}_{r_1,\cdots,r_{2n}}
\prod_{i=1}^{n^\prime}\~\Psi^{(1),X_2,X_3}_0(\bar \rho_{i})
\times \prod_{i=1}^{a} 
S^{X_2,X_3}_{\delta^{(a)\prime}_i}
(t^{(a)}_i)
\times \prod_{i=a+1}^{b} 
S^{X_2,X_3}_{\delta^{(b)\prime}_i}
(t^{(b)}_i) \nn\\
&&\times \prod_{i=1}^{a} 
E^+_{\epsilon^{(a)\prime}_i}(u^{(a)}_i) 
\times \prod_{i=a+b+1}^{a+b+c} 
E^+_{\epsilon^{(c)\prime}_i}(u^{(c)}_i):, 
\ea
with
\bac
\~\Psi^{(1),X_1}_0(\bar \rho)&=&
\exp(X_1(1;2,k+2|q^{k-2}\bar \rho;-{k+2\over 2})),\nn\\
\~\Psi^{(1),X_2,X_3}_0(\bar \rho)&=&
\exp(X_2(2|q^{-2}\bar \rho;0)+X_3(2|q^{-2}\bar \rho;0)).
\ea
Using the latter expressions and  the general 
trace formulas 
in the appendix we  find then
\bac
S^{X_1}(s,r_{2n},\cdots,r_{1};
r^\prime_{n^\prime},\cdots,r^\prime_{1})
&=&
T_{\ell_s}^{X_1}
(2n,n^\prime,N+n_s,M,n_4|X,Y,Z,W),\nn\\
S^{X_2,X_3}(r_{2n},\cdots,r_{1};
r^\prime_{n^\prime},\cdots,r^\prime_{1})
&=&
T^{X_2,X_3}(A,B,C,D,G,H|\~X,\~Y,\~Z,\~W),
\ea
with   
\bac
\ell_s&=& \ell_{\~n-Ps,\bar n}=\~n-Ps-1-\bar n{P\over P^\prime},\nn\\
N&=&\sum_{l=1}^{2n}r_l=nk,\nn\\
M&=&2kn-\sum_{l=1}^{2n}j_{l}=nk+\sum_{l=1}^n i_l-
\sum_{l=1}^n j_l,\nn\\
t_{g+\sum_{i=1}^{l-1}r_i}&=&t^{(l)}_g,\quad 1\leq g\leq r_l,\nn\\
t_{N+i}&=&t_i^{(a)},\quad 1\leq i\leq a,\nn\\
t_{N+a+i}&=&t_i^{(b)},\quad 1\leq i\leq b,\nn\\
\xi_{g+\sum_{i=1}^{l-1}(k-j_i)}&=&\xi^{(l)}_g,\quad 1\leq g\leq 
k-j_l,\nn\\
u_i&=&u^{(a)}_i,\quad 1\leq i\leq a,\nn\\
u_{a+i}&=&u^{(c)}_i,\quad 1\leq i\leq c,\nn\\
\delta_{g+\sum_{i=1}^{l-1}r_i}&=&
\delta^{(l)}_g,\quad 1\leq g\leq r_l,\nn\\
\delta_{N+i}&=&\delta_i^{(a)\prime},\quad 1\leq i\leq a,\nn\\
\delta_{N+i}&=&\delta_i^{(b)\prime},\quad a+1\leq i\leq b,\nn\\
\epsilon^\prime_{i}&=&\epsilon^{(a)\prime}_i,\quad 
1\leq i\leq a,\nn\\
\epsilon^\prime_{a+i}&=&\epsilon^{(c)\prime}_i,
\quad 1\leq i\leq c,\nn\\
X&=&\{z_1,\cdots,z_{2n}\},\nn\\
Y&=&\{\bar \rho_1,\cdots,\bar \rho_{n^\prime}\},\nn\\
Z&=&\{t_1,\cdots,t_{N+a+b}\},\nn\\
W&=&\{\xi_1,\cdots,\xi_M\},\nn\\
A&=&\{\delta_1-k-2,\cdots,\delta_{N+a+b}-k-2,
\epsilon^\prime_{1}-k-2,\cdots \epsilon^\prime_{a+c}-k-2
\},\nn\\
B&=&\{(\epsilon_1-1)(k+1)-1,\cdots,(\epsilon_M-1)(k+1)-1,-2
,\cdots,-2\},\nn\\
C&=&\{-k-2,\cdots,-k-2\},\nn\\
D&=&\{-1,\cdots,-1,1,\cdots,1\},\nn\\
G&=&\{(\epsilon_1-1)(k+2),\cdots,(\epsilon_M-1)(k+2),-2,
\cdots,-2\},\nn\\
H&=&\{-1,\cdots,-1,0,\cdots,0\},\nn\\
\~X&=&\{t_1,\cdots,t_{N+a+b}, 
u_1,\cdots,
u_{a+c}\},\nn\\
\~Y&=&\{\xi_1,\cdots,\xi_M,\bar \rho_1,\cdots,\bar \rho_{n^\prime}\},
\nn\\
\~Z&=&\{t_1,\cdots,t_{N+a+b},u_1,\cdots,u_{a+c}\},\nn\\
\~W&=&\{\xi_1,\cdots,\xi_M,\bar \rho_1,\cdots,\bar \rho_{n^\prime}\}.
\ea
The dimensions of the latter sets and any other details should
be  clear from the contexts. 
Just as in the previous case of the correlation functions, 
though the above 
expression for the N-point form factor is very complicated,
one
can still deduce from it some simple selection rules. 
Indeed, the general trace formulas 
(\ref{ap4}) and (\ref{ap5}) in the appendix imply
\bac
S^{X_1}(s,r_{2n},\cdots,r_{1};
r^\prime_{n^\prime},\cdots,r^\prime_{1})
&=&\delta_{n^\prime-2(a+b),0}(\cdots),\nn\\
S^{X_2,X_3}(r_{2n},\cdots,r_{1};
r^\prime_{n^\prime},\cdots,r^\prime_{1})
&=&\delta_{\sum_{l=1}^{n} i_l+{n^\prime\over 2},
\sum_{l=1}^n j_l+a+c}
(\cdots).
\ea
Therefore the following selection rules must be satisfied 
in order for the form factors to be non-vanishing:
\bac
a+b&=&n^\prime/2,\nn\\
\sum_{l=1}^{n} i_l+{n^\prime\over 2}&=&\sum_{l=1}^n j_l+
a+c.
\ea
Since $a+b$ is an integer number, $n^\prime$ is an even 
integer number. These selection rules are also consistent 
with those of
the spin 1/2 case.

Let us now recapitulate the main result 
of this calculation. 
The form factor
$<\lambda_0|E_{i_nj_n}\otimes \cdots
\otimes E_{i_1j_1}|\rho_{n^\prime},\cdots,
\rho_1>_{1-\alpha_{n^\prime},\cdots,
1-\alpha_{1}}^{\beta_{n^\prime},\cdots,\beta_{1}}$ of the
basis local operator 
$E_{i_nj_n}\otimes \cdots
\otimes E_{i_1j_1}$ is given by
\bac
&&<\lambda_0|E_{i_nj_n}\otimes \cdots
\otimes E_{i_1j_1}|\rho_{n^\prime},\cdots,
\rho_1>_{1-\alpha_{n^\prime},\cdots,
1-\alpha_{1}}^{\beta_{n^\prime},\cdots,\beta_{1}}\nn\\
&&={B(q)\prod_{i=1}^{n^\prime}\left(d_{\alpha_{i},\beta_i}
h^{\lambda_{i}}_{\lambda_{i-1}}(\bar \rho_{i})\right)
C_{j_{2n},\cdots,j_1}(z_{2n},\cdots,z_1)
\over Tr_{V({\lambda^\prime}_0)}(\mu^{-d}\nu^{-\alpha_0})}\nn\\
&&\times
\left\{\prod_{l=1}^{2n}
\{\sum_{\delta_1^{(l)},\cdots,\delta_{r_l}^{(l)}}
\sum_{\epsilon_1^{(l)},\cdots,\epsilon_{k-j_l}^{(l)}}
\int_{q^k pz}^{q^kz\infty} 
d_{p}t_1^{(l)}\cdots  
\int_{q^k pz}^{q^kz\infty} 
d_{p}t_{r_l}^{(l)} 
\oint {d\xi_1^{(l)}\over 2\pi i}
\cdots \oint {d\xi_{k-j_l}^{(l)}\over 2\pi i}\right.\nn\\ 
&&\times \left. H_{\delta_1^{(l)},\cdots,\delta_{r_l}^{(l)},
\epsilon_1^{(l)},
\cdots, \epsilon_{k-j_l}^{(l)}}
(z_l^{(l)},t_1^{(l)},\cdots,t_{r_l}^{(l)},\xi_1^{(l)},
\cdots,\xi_{k-j_l}^{(l)},
{\cal C}_l)\}\prod_{2n\geq l_1>l_2\geq
1}f(D_{l_1},D_{l_2})\right\}\nn\\
&&\times \left\{
\prod_{i=1}^{a}\left(
\sum_{\delta^{(a)\prime}_i,\epsilon^{(a)\prime}_i}
\int_{0}^{q^{k-1}z} d_p t^{(a)}_i \oint{du^{(a)}_i\over 2\pi i}
I_{\delta^{(a)\prime}_i,\epsilon^{(a)\prime}_i}^{(1),(1)}
(\bar \rho_i,t^{(a)}_i,u^{(a)}_i|C)\right) \right.\nn\\
&&\times  \prod_{i=a+1}^{b}\left(
\sum_{\delta^{(b)\prime}_i} 
\int_{0}^{q^{k-1}z} d_p t^{(b)}_i
I_{\delta^{(b)\prime}_i}^{(1),(0)}
(\bar \rho_i,t^{(b)}_i|C)\right)
\prod_{i=b+1}^{c}\left(
\sum_{\epsilon^{(c)\prime}_i} \oint{du^{(c)}_i\over 2\pi i}
I_{\epsilon^{(c)\prime}_i}^{(0),(1)}
(\bar \rho_i,u^{(c)}_i|C)\right)\nn\\ 
&&\times \left. \prod_{1\leq i^\prime
<i\leq n^\prime} f(K_i,K_{i^\prime})\right\}\times 
f({\cal O}_{r_{2n},\cdots,r_{1}}, {\cal O}_{
r^\prime_{n^\prime},\cdots,r^\prime_1})\nn\\
&& 
\times\left(\sum_{s\in 2\Z} {\cal A}^{n(k+sp-p)}_{\~\Phi,k}
 {\cal A}_{\~\Psi,1}^{spn^\prime/2}
S(s,r_{2n},\cdots,r_{1};
r^\prime_{n^\prime},\cdots,r^\prime_{1})\right.\nn\\
&&\left. -\sum_{s\in 2\Z+1} 
{\cal A}^{n(k+sp-p)}_{\~\Phi,k}
{\cal A}_{\~\Psi,1}^{
(2n-n^\prime+(s-1)p+1)n^\prime/2+
2\sum_{l=1}^{n^\prime-1}lr^\prime_{l+1}} 
S(s,r_{2n},\cdots,r_{1};
r^\prime_{n^\prime},\cdots,r^\prime_{1})\right).\nn\\
\label{co7}
\ea  
Here also all the various symbols are c-functions and 
are defined
throughout the above analysis. Moreover, the selection rules 
leading to non-vanishing form factor can be rewritten in the
 form
\bac
\sum_{l=1}^n i_l+{n^\prime\over 2}&=&\sum_{l=1}^n j_l+
\sum _{i=1}^{n^\prime} \alpha_i,\nn\\
\sum_{i=1}^{n^\prime} \beta_i=0.
\ea
It is clear from the latter selection rules that since 
 $n^\prime$ must be even,  local
operators create spinons only in pairs, which is identical to
the spin $1/2$ case.

Now that we understand how the vertex operators are 
the building blocks of 
the form factors we can derive more relations satisfied by 
them using the
various basic relations among the vertex operators. 
For example, if
we define a general matrix element of a local operator 
${\cal O}$ as
\bac
&&_{\alpha_n,\cdots,\alpha_{m+1}}^{\beta_n,\cdots,\beta_{m+1}}
<\rho_n,\cdots,\rho_{m+1}|{\cal O}|\rho_{m},\cdots,
\rho_1>_{\alpha_{m},\cdots,
\alpha_{1}}^{\beta_{m},\cdots,\beta_{1}}\nn\\
&&={{\rm Tr}_{V(\lambda_n)}(\mu^{-d}\nu^{-\alpha_0}
\Psi^{V{\lambda}_{n}}_{
{\lambda}_{n-1},\alpha_{n}}
(\rho_{n})\cdots 
\Psi^{V{\lambda}_{m+1}}_{{\lambda}_{m},
\alpha_{m+1}}(\rho_{m+1})
{\cal O}
\Psi^{V^*{\lambda}_{m}}_{
{\lambda}_{m-1},\alpha_{m}}
(\rho_{m})\cdots 
\Psi^{V^*{\lambda}_{1}}_{{\lambda}_{0},
\alpha_{1}}(\rho_{1}))\over {\rm Tr}_{V(\lambda_n)}
(\mu^{-d}\nu^{-\alpha_0})}.\nn\\
\ea
then using the first relation of (\ref{cco7}) we arrive at 
the following identity: 
\bac
&&_{\alpha_n,\cdots,\alpha_{m+1}}^{\beta_n,\cdots,\beta_{m+1}}
<\rho_n,\cdots,\rho_{m+1}|{\cal O}|\rho_{m},\cdots,
\rho_1>_{\alpha_{m},\cdots,
\alpha_{1}}^{\beta_{m},\cdots,\beta_{1}}\nn\\
&&=\left(\prod_{i=m+1}^n d^{-1}_{\alpha_i,\beta_i}\right)
<{\cal O}|q^{-2} \rho_n,\cdots, q^{-2}\rho_{m+1},\rho_{m},
\cdots,\rho_1>_{1-\alpha_n,\cdots,1-\alpha_{m+1},
\alpha_{m},\cdots,
\alpha_{1}}^{\beta_n,\cdots,\beta_{m+1},\beta_{m},\cdots,
\beta_{1}}.\nn\\
\ea

This identity is particularly useful in situations where  matrix
elements of local operators mapping excited states to excited
states are needed.  It means that one needs
just to evaluate the form factors of local operators
 mapping excited states to the vacuum. An example of such a
situation is the evaluation of the dynamic correlation
functions at arbitrary nonzero temperature \cite{Fabal96}.

For completeness, let us now summarise some of the 
interesting results previously found in the
literature with respect to the form factors. 
It has been
shown in Ref. \cite{Idzal93} that the form factors satisfy
 the higher level q-KZ equation. This is given in a matrix
form. Let us keep most of the same notations as in this 
reference where
the form factor is defined in a matrix form (up to
normalization factor) as
\be
F_{\lambda_n,\cdots,\lambda_0}(\rho_n,\cdots,\rho_1)=
{\rm Tr}_{V(\lambda_0)}\left(q^{-2\rho} {\cal O}_{\lambda_0}
\Psi^{V^*_n\lambda_0}_{\lambda_{n-1}}(\rho_n)\cdots
\Psi^{V^*_1\lambda_1}_{\lambda_{0}}(\rho_1)\right).
\ee
Then using the commutation relations, the cyclic property
of a trace
it has been shown that the following q-KZ equation is
satisfied: 
\bac
&&
F_{\lambda_n,\cdots,\lambda_0}(\rho_n,\cdots,
q^4\rho_i,\cdots,\rho_1) \nn\\
&&=R^{*i-1i}(\rho_{i-1}/q^4\rho_i)^{-1}\cdots R^{*1i}(\rho_1/q^4
\rho_i)^{-1}\nn\\
&&\times \pi_{V^*}(q^{-2\bar \rho})
R^{*in}(\rho_{i}/\rho_n)\cdots R^{*ii+1}(\rho_i/
\rho_{i+1})\nn\\
&&\times \sum_{\lambda^\prime_n,\cdots,\lambda^\prime_{i+1},
\lambda^\prime_{i-1},\cdots,\lambda^\prime_1} 
q^{4(\Delta_{\lambda^\prime_1}-\Delta(\lambda_0))} 
F_{\lambda^\prime_1,\lambda^\prime_n,\cdots,\lambda_i,\cdots,
\lambda^\prime_1}(\rho_n,\cdots, \rho_1)\nn\\
&&\times
W\left.\left(\matrix{\lambda_0&\lambda_1\cr
\lambda_1^\prime&\lambda_2^\prime\cr}\right|
q^4 \rho_i/\rho_1\right)\cdots 
W\left.\left(\matrix{\lambda_{i-2}&\lambda_{i-1}\cr 
\lambda_{i-1}^\prime&
\lambda_i\cr}\right| q^4 \rho_i/\rho_{i-1}\right)\nn\\
&&\times
W\left.\left(\matrix{\lambda_i&\lambda_{i+1}\cr
\lambda_{i+1}^\prime&\lambda_{i+2}^\prime\cr}\right|
 \rho_i/\rho_{i+1}\right)\cdots 
W\left.\left(\matrix{\lambda_{n-1}&\lambda_{n}\cr 
\lambda_{n}^\prime&
\lambda^\prime_1\cr}\right| q^4 \rho_i/\rho_{i-1}\right),
\ea
where $\bar \rho$ is the restriction of $\rho=2d+\alpha_0/2$ 
to a  finite dimensional representation, $R^{ij}(z)$ is
$R(z)$ acting on $V_i\otimes V_j$, 
and similarly for $R^*(z)$. Moreover,
\bac
V^*&=&V^{(1)*a^{-1}},\nn\\
\Delta_\lambda&=&{(\lambda,\lambda+2\rho)\over 2(k+2)}.
\ea  
This means that the   explicit expression for the
form factors we found in (\ref{co7}) 
provides the unique physical 
solution for the q-KZ equation. This is  nontrivial 
because typically difference 
equations such as q-KZ admit an
infinite number of solutions which differ from each other 
by pseudo constants. 
Further  analiticity requirements based on
physical considerations are needed in order to single out
a unique solution. 
It would be interesting to try to solve direclty
the q-KZ equation consistently with the analyticity requirements
which can be read off from the explicit expressions above. 
This might lead
 to other  more compact forms for the form factors. 

The other interesting result that has been obtained in the
literature is the fact that the form factors satisfy 
a lattice version of Smirnov's axioms. This is an
important step towards making Smirnov's axioms universal by
including not only  integrable models of
  quantum field theory, 
as originally thought, but also lattice models with different
dispersion relations. 
Let us briefly remind Smirnov's 
axioms as applied to integrable quantum field theories
\cite{Smi92},
and as summarised in Ref. \cite{Kon94}:\\
Axiom 1: Form factors
$F(\beta_1,\cdots,\beta_n)_{\epsilon_1,\cdots,\epsilon_n}$ 
 have  the $S$-matrix symmetry
\be
F(\beta_1,\cdots,\beta_i,\beta_{i+1},\cdots,\beta_n)_
{\epsilon_1,\cdots,\epsilon_i,\epsilon_{i+1},\cdots,\epsilon_n}
S^{\epsilon^\prime_i\epsilon^\prime_{i+1}}_
{\epsilon_i\epsilon_{i+1}}=
F(\beta_1,\cdots,\beta_i,\beta_{i+1},\cdots,\beta_n)_
{\epsilon_1,\cdots,\epsilon^\prime_{i+1},
\epsilon^\prime_{i},\cdots,\epsilon_n}.
\ee 
\\
Axiom 2: Form factors satisfy the difference equation
\be
F(\beta_1,\cdots,\beta_n+2\pi i)_
{\epsilon_1,\cdots,\epsilon_n}  
=F(\beta_n,\beta_1,\cdots,\beta_{n-1})_
{\epsilon_n,\epsilon_1,\cdots,\epsilon_{n-1}}.  
\ee
Axiom 3: Form factors have the annihilition
property. This means that the residue of the 
$n$-particle form factor
$F(\beta_1,\cdots,\beta_n)_{\epsilon_1,\cdots,\epsilon_n}$ 
 can be expressed in terms of
$(n-2)$-particle ones. Moreover, they 
have simple poles at $\beta_i=\beta_j+i\pi$, $i>j$, with 
residues at $\beta_n=\beta_{n-1}+i\pi $ 
given by
\bac
&2\pi i& {\rm Res}
F(\beta_1,\cdots,\beta_n)_{\epsilon_1,\cdots,\epsilon_n}\nn\\  
&&= 
F(\beta_1,\cdots,\beta_{n-2})_{\epsilon^\prime_1,\cdots,
\epsilon^\prime_{n-1}} \delta_{\epsilon_n,-\epsilon_{n-1}}\nn\\
&&\times \left(
\delta_{\epsilon_1}^{\epsilon_1^\prime}\cdots 
\delta_{\epsilon_{n-2}}^{\epsilon_{n-2}^\prime}-
S^{\epsilon^\prime_{n-1}\epsilon^\prime_{1}}_{
\tau_1 \epsilon_1}(\beta_{n-1}-\beta_1)\cdots 
S^{\tau_{n-3}\epsilon^\prime_{n-2}}_{
\epsilon_{n-1} \epsilon_{n-2}}(\beta_{n-1}-\beta_{n-2})\right).
\ea
These  axioms provide a powerful alternative 
method for the exact derivation of  form factors of local 
operators in
integrable quantum field theories. 
They also lead to the
q-KZ equation for form factors.

A crucial result was obtained in Ref. \cite{Pak93} for the spin
1/2 case where it was shown that Smirnov's axioms are
still satisfied
by the form factors of the Heisenberg model. Then, this result
has been generalized for the third axiom to the higher spin 
Heisenberg model. To describe this result let us use
the following "light" definition for a component   form factor,
which is sufficient for this purpose:
\be
F_{\lambda,\lambda_{n-1},\cdots,\lambda_1,
\lambda}^{\epsilon_n,\cdots \epsilon_1}
 (\rho_n,\cdots,\rho_1)=
{\rm Tr}_{V(\lambda)}\left(q^{-2\rho} {\cal O}_{\lambda}
\Psi^{\lambda}_{\lambda_{n-1},\epsilon_n}(\rho_n)\cdots
\Psi^{\lambda_1}_{\lambda_{0},\epsilon_1}(\rho_1)
\right).
\ee
Recall that $\epsilon =+$ and $\epsilon=-$ correspond to 
the components $m=0$ and $m=1$ respectively. 
That only simple poles appear in the factorization
property of the latter form factors is traced to the following
OPE's of type II vertex operators in 
the limit $q^{-2}\rho_2/\rho_1\rightarrow 1$:
\be
\Psi^{\mu}_{\lambda_{\pm},\epsilon_2}(\rho_2)
\Psi^{\lambda_{\pm}}_{\lambda,\epsilon_1}(\rho_1)
\sim {\cal N}_\lambda^{\lambda_\pm}(\epsilon_2,\epsilon_1)
\delta_{\epsilon_2,-\epsilon_1} \delta_{\lambda}^{\mu} 
{1\over 1-q^{-2} \rho_1/\rho_2} {\rm id} +O(1),
\ee
where the constants ${\cal N}_\lambda^{\lambda_\pm}
(\epsilon_2,\epsilon_1)$ are given by
\bac
{\cal N}_\lambda^{\lambda_-}(+,-)&=&
 -q {\cal N}_\lambda^{\lambda_-}(-,+)=
-q^4(1-p) {(q^2p;p)_\infty \over (q^{-2}p;p)_\infty}
{\xi(q^2;1,q^4)\over B_p(1-2s,-2s)},\nn\\
{\cal N}_\lambda^{\lambda_+}(+,-)&=&
 -q {\cal N}_\lambda^{\lambda_+}(-,+)=
q(1-p) {(q^2p;p)_\infty \over (q^{-2}p;p)_\infty}
{\xi(q^2;1,q^4)\over B_p(6s,-2s)}.
\ea
Using these OPE's, the cyclic property of a trace  
and the commutation relations of type II vertex operators 
(\ref{}), the
following residue formula at $q^{-2}\rho_2/\rho_1=1$ 
was derived in Ref. \cite{Kon94}:
\bac
&& {\rm Res} 
F_{\lambda,\cdots,\lambda_1,\lambda}^{\epsilon_n,\cdots 
\epsilon_1}
 (\rho_n,\cdots,\rho_1)\nn\\
&&=\delta_\lambda^{\delta_2}\delta_{\epsilon_2,-\epsilon_1}
{\cal N}_\lambda^{\lambda_1}(-\epsilon_1,\epsilon_1)
F_{\lambda,\cdots,\lambda_3,
\lambda}^{\epsilon_n,\cdots \epsilon_3} 
(\rho_n,\cdots,\rho_3)\nn\\
&&-{\cal N}_{\lambda_1}^{\lambda}(\epsilon_1,-\epsilon_1)
\delta_{\tau_{n-2},-\epsilon_1}
R^{\alpha_n \tau_{n-2}}_{\epsilon_n\tau_{n-3}}(\rho_n/\rho_2)
R^{\alpha_{n-1}\tau_{n-3}}_{\epsilon_{n-1}\tau_{n-4}}
(\rho_{n-1}/\rho_2)\cdots
R^{\alpha_4 \tau_{2}}_{\epsilon_3\tau_{1}}(\rho_4/\rho_2)
R^{\alpha_{3}\tau_{1}}_{\epsilon_{3}\epsilon_{2}}
(\rho_{3}/\rho_2)\nn\\
&&\times \sum_{\mu_1,\cdots,\mu_{n-3}}
W\left.\left(\matrix{\mu_{n-3}&\lambda_{n-1}\cr 
\lambda_1&\lambda\cr}\right| \rho_n/\rho_2\right)    
W\left.\left(\matrix{\mu_{n-4}&\lambda_{n-2}\cr 
\lambda_{n-3}&\lambda_{n-1}\cr}\right| \rho_{n-1}/\rho_2\right)
\cdots
W\left.\left(\matrix{\lambda_{1}&\lambda_{2}\cr 
\mu_{1}&\lambda_{3}\cr}\right| \rho_3/\rho_2\right)\nn\\
&&\times    
F_{\lambda_1,\mu_{n-3}\cdots,\mu_1,
\lambda_1}^{\alpha_n,\alpha_{n-1}\cdots \alpha_2,\alpha_3} 
(\rho_n,\cdots,\rho_3).
\ea
As noted in Ref. \cite{Kon94}, the main new ingredient of the
lattice version of the third axiom is that now the face
type Boltzmann weights $W\left.\left(\matrix{\lambda&\mu\cr 
\mu^\prime&\nu\cr}\right| \rho\right)$ enter the residue formula.

On can also  derive the higher spin lattice
version of the first and second axioms. 
Indeed, using the commutation relations of type II vertex
operators (\ref{vo7}) we find the lattice analogue of the 
first axiom
\bac
&&F_{\lambda,\cdots,
\lambda_{i+1},\lambda_i,\cdots,\lambda_1,\lambda}^
{\epsilon_n,\cdots, \epsilon_{i+1},\epsilon_i,\cdots,\epsilon_1}
(\rho_n,\cdots,\rho_{i+1},\rho_i,\cdots,\rho_1)=\nn\\
&&
\sum_{\mu=\lambda_{i_1\pm}}
F_{\lambda,\cdots,
\lambda_{i-1},\mu,\cdots,\lambda_1,\lambda}^
{\epsilon_n,\cdots,
\epsilon^\prime_{i},\epsilon_{i+1}^\prime,\cdots,\epsilon_1}
(\rho_n,\cdots,\rho_i,\rho_{i+1},\cdots,\rho_1)
R^{\epsilon^\prime_{i+1} \epsilon^\prime_i}_
{\epsilon_{i+1} \epsilon_i}(\rho_{i+1}/\rho_i)
W\left.\left(\matrix{\lambda_{i+1}&\lambda_i\cr
\mu&\lambda_{i+1}\cr}\right| \rho_{i+1}/\rho_i\right).\nn\\
\ea
Moreover, (\ref{vvo9}) and (\ref{vvo11}) imply
\be
q^{-2\rho}    
\Psi^{\mu}_{\lambda,\epsilon}(z)q^{2\rho}=
q^{\epsilon}\Psi^{\mu}_{\lambda,\epsilon}(zq^4).
\ee
This relation  and the cyclic property of a trace 
lead to the lattice version of the second axiom, i.e.,
\be
F_{\lambda,\lambda_{n-1},\cdots,\lambda_1,\lambda}^{\epsilon_n,
\epsilon_{n-1},\cdots 
\epsilon_1}
 (\rho_n,\cdots,\rho_1)=
q^{\epsilon_n} 
F_{\lambda_{n-1},\cdots,\lambda_1,\lambda,
\lambda_{n-1}}^{
\epsilon_{n-1},\cdots 
\epsilon_1,\epsilon_n}
 (\rho_{n-1},\cdots,\rho_1,\rho_n q^4).
\ee
Therefore, the crucial point is that Smirnov's axioms extend
to integrable  lattice models, which is 
 beyond their original scope of integrable 
quantum field theories.
However, when applied to lattice models they also involve in
general 
the face type Boltzmann weights $W$ in addition to the 
usual $R$ matrix.

We see that despite the fact that the final
expression for the form factors is quite complicated one can 
still obtain simple information from it, such as the selection 
rules and  the way spinons are created 
by local operators. 
One of the main lessons one learns from this type of analysis
is that the identification of the creation and annihilation
operators of eigenstates of a Hamiltonian 
is an important step towards the
complete solution of a model. The next crucial step is handling
the commutation relations of the operators which are typically
generalized commutation relations. Since the only
commutation relations we know how to handle completely are
those of Heisenberg and Clifford types we try then to express
the creation and annihilation of the eigenstates 
in terms of the generators of the latter algebras. 
If we succeed in 
overcoming this step then the main goal of evaluating the
relevant physical quantities like correlation functions becomes
a matter of technical details, as was shown here. Needless, to
say that one would like to extend this method to as many 
models as possible, and in particular the Heisenberg model
with nonzero temperature and/or nonzero magnetic field, the
8-vertex model, the Potts model etc... The more we learn
about the generalized commutation relations the more
we learn  about the exact solution of integrable models at
the level of correlation functions and form factors and not
just the spectrum.

\section{\bf Acknowlegment}

This work is supported by the 
NSF Grant \# PHY9309888.

\newpage
\section{Appendix}

\subsection {Operator product expansions (OPE's) of the 
elementary operators}

The elementary operators that appear in the evaluation of the
correlation functions and  form factors  are
the vertex operators 
  $\~\Phi^{(\ell)}_\ell(z)$ (type I) and 
$\~\Psi_{0}^{(\ell)}(z)$ (type II), the currents
$E^\pm_\epsilon(z)$, 
the screening currents 
$S_\delta(z)$, $\eta(z)$, 
and $\xi(z)$. They are bosonized
as follows:
\bac
\~\Phi^{(\ell)}_\ell(z)&=&:\exp{(X_1(\ell;2,k+2|q^{k}z;{k+2\over
2}))}:,\nn\\
\~\Psi_{0}^{(\ell)}(z)&=&
:\exp({X_1(\ell;2,k+2|q^{k-2}z;-{k+2\over 2})+
X_2(\ell;2,1|q^{-2}z;0)+X_3(\ell;2,1|q^{-2}z;0)}):,\nn\\
E^-_\epsilon(z)&=& :\exp(\partial X_1^{(\epsilon)}(q^{-2}z;
-{k+2\over 2})+X_2(2|q^{(\epsilon-1)(k+2)}z;-1)+
X_3(2|q^{(\epsilon-1)(k+1)-1}z;0)):,\nn\\
E^+_\rho(z)&=& :\exp(-X_2(2|q^{-k-2}z;1)-
X_3(2|q^{-k-2+\rho}z;0)):,\nn\\
S_\delta(z)&=& :\exp(-X_1(k+2|q^{-2}z;
-{k+2\over 2})-X_2(2|q^{-k-2}z;-1)-
X_3(2|q^{-k-2+\delta}z;0)):,\nn\\
\eta(z)&=&:\exp({X_3(2|q^{-k-2}z;0)}):,\nn\\
\xi(z)&=&:\exp({-X_3(2|q^{-k-2}z;0)}):.
\ea
where
\be
\partial X_1^{(\epsilon)}(q^{-2}z;
-{k+2\over 2})=\epsilon \{(q-q^{-1})\sum_{n=1}^\infty 
a_{X_1,\epsilon n}z^{-\epsilon n}q^{(2\epsilon -{k+2\over 2})n}
+a_{X_1,0}\log(q)\}.
\ee
The correlation functions and form factors are given in terms
of the c-function $f$  defined by: 
\be
:A::B:=:AB:f(A,B)=:AB:\prod_{i=1}^n\prod_{j=1}^m f(A_i,B_j),
\ee  
for the composite normal ordered operators 
$A=:\prod_{i=1}^n A_i:$ 
and $B=:\prod_{j=1}^m B_j:$, with $A_i$ and $B_j$ being 
elementary operators. 
Therefore to evaluate  $f(A,B)$  for any composite
operators $A$ and $B$, it is sufficient to know just
the OPE's among all elementary operators.
The latter OPE's are  presented below.

OPE's:
\bac
\~\Phi^{(\ell)}_\ell(z)\~\Phi^{(\ell)}_\ell(w)
&=&(zq^k)^{\ell^2/2(k+2)}
{(pq^{2(1-\ell)}wz^{-1};p;q^4)_\infty 
(pq^{2(1+\ell)}wz^{-1};p;q^4)_\infty\over
(pq^{2}wz^{-1};p;q^4)^2_\infty}
:\~\Phi^{(\ell)}_\ell(z)\~\Phi^{(\ell)}_\ell(w):,\nn\\
\~\Psi^{(\ell)}_{0}(z)\~\Psi^{(\ell)}_{0}(w)
&=&(zq^{k-2})^{\ell^2/2(k+2)}
{(q^{2(1-\ell)}wz^{-1};p;q^4)_\infty 
(q^{2(1+\ell)}wz^{-1};p;q^4)_\infty\over
(q^{2}wz^{-1};p;q^4)^2_\infty}:\~\Psi^{(\ell)}_{0}(z)
\~\Psi^{(\ell)}_{0}(w):,\nn\\
\~\Phi^{(\ell)}_\ell(z)\~\Psi^{(\ell)}_{0}(w)
&=&(zq^{k})^{\ell^2/2(k+2)}
{(q^{k+2(1-\ell)}wz^{-1};p;q^4)_\infty 
(q^{k+2(1+\ell)}wz^{-1};p;q^4)_\infty\over
(q^{k+2}wz^{-1};p;q^4)^2_\infty}\\
&&\times :\~\Phi^{(\ell)}_\ell(z)
\~\Psi^{(\ell)}_{0}(w):,\nn\\
\~\Psi^{(\ell)}_{0}(z)\~\Phi^{(\ell)}_\ell(w)
&=&(zq^{k-2})^{\ell^2/2(k+2)}
{(q^{k+2(3-\ell)}wz^{-1};p;q^4)_\infty 
(q^{k+2(3+\ell)}wz^{-1};p;q^4)_\infty\over
(q^{k+6}wz^{-1};p;q^4)^2_\infty}\\
&&\times :\~\Psi^{(\ell)}_{0}(z)
\~\Phi^{(\ell)}_\ell(w):.\nn\\
\label{ap9}
\ea

\bac
\~\Phi^{(\ell)}_\ell(z)E^-_1(w)&=&
:\~\Phi^{(\ell)}_\ell(z)E^-_1(w):,\nn\\
\~\Phi^{(\ell)}_\ell(z)E^-_{-1}(w)&=& 
{z-q^{-\ell-k-2}w\over z-q^{\ell-k-2}w} 
:\~\Phi^{(\ell)}_\ell(z)E^-_{-1}(w):,
\quad |z|>|q^{-\ell-k-2}w|,\nn\\
E^-_1(z)\~\Phi^{(\ell)}_\ell(w)&=& 
{q^\ell z-q^{k+2}w\over z-q^{\ell+k+2}w}
: E^-_1(z)\~\Phi^{(\ell)}_\ell(w):,\quad 
|z|>|q^{-\ell+k+2}w|,\nn\\
E^-_{-1}(z)\~\Phi^{(\ell)}_\ell(w)&=& 
q^{-\ell}:E^-_-(z)\~\Phi^{(\ell)}_\ell(w):,\nn\\
\~\Psi^{(\ell)}_{0}(z)E^+_1(w)&=&
:\~\Psi^{(\ell)}_{\ell}(z)E^+_1(w):,\nn\\
\~\Psi^{(\ell)}_{0}(z)E^+_{-1}(w)&=& 
{z-q^{\ell-k}w\over z-q^{-\ell-k}w} 
:\~\Psi^{(\ell)}_{0}(z)E^+_{-1}(w):,\quad |z|>|q^{-k-\ell}w|,
\nn\\
E^+_{1}(z)\~\Psi^{(\ell)}_{0}(w)&=& 
q^{-\ell}{z-q^{\ell+k}w\over z-q^{k-\ell}w} 
:E^+_{1}(z)\~\Psi^{(\ell)}_{0}(w):,\quad |z|>|q^{k-\ell}w|,
\nn\\
E^+_{-1}(z)\~\Psi^{(\ell)}_{0}(w)&=& 
q^{\ell} :E^+_{-1}(z)\~\Psi^{(\ell)}_{0}(w):,
\ea

\bac
\~\Phi^{(\ell)}_\ell(z)S_1(w)&=&(q^kz)^{-\ell/(k+2)}
{(q^\ell zw^{-1};p)_\infty \over (q^{-\ell}
zw^{-1};p)_\infty}:\Phi^{(\ell)}_\ell(z)S_1(w):,\quad 
|z|>|q^{-\ell}w|,\nn\\ 
\~\Phi^{(\ell)}_\ell(z)S_{-1}(w)&=&(q^kz)^{-\ell/(k+2)}
{(q^\ell zw^{-1};p)_\infty \over (q^{-\ell}
zw^{-1};p)_\infty}:\~\Phi^{(\ell)}_\ell(z)S_{-1}(w):,\quad 
|z|>|q^{-\ell}w|,\nn\\
S_1(z)\~\Phi^{(\ell)}_\ell(w)&=& (q^{-2}z)^{-\ell/(k+2)}
{(q^\ell pwz^{-1};p)_\infty \over (q^{-\ell}
pwz^{-1};p)_\infty}:S_1(z)\~\Phi^{(\ell)}_\ell(w):,\quad 
|z|>|q^{-\ell}p w|,\nn\\
S_{-1}(z)\~\Phi^{(\ell)}_\ell(w)&=& (q^{-2}z)^{-\ell/(k+2)}
{(q^\ell pwz^{-1};p)_\infty \over (q^{-\ell}
pwz^{-1};p)_\infty}:S_{-1}(z)\~\Phi^{(\ell)}_\ell(w):,\quad 
|z|>|q^{-\ell}p w|,\nn\\
\~\Psi^{(\ell)}_{0}(z)S_1(w)&=&(q^{k-2}z)^{-1/k+2}
{(q^{1-k}p wz^{-1};p)_\infty \over (q^{-1-k}p
wz^{-1};p)_\infty}:\~\Psi^{(\ell)}_{0}(z)S_1(w):,
\quad |z|>|q^{-k-1}p w|,\nn\\ 
\~\Psi^{(\ell)}_{0}(z)S_{-1}(w)&=&
(q^{k-2}z)^{-1/k+2}
{(q^{1-k} wz^{-1};p)_\infty \over (q^{-1-k}
wz^{-1};p)_\infty}:\~\Psi^{(\ell)}_{0}(z)S_1(w):,
\quad |z|>|q^{-k-1} w|,\nn\\
S_1(z)\~\Psi^{(\ell)}_{0}(w)&=&
(q^{-2}z)^{-1/k+2}
{(q^{1+k} wz^{-1};p)_\infty \over (q^{k-1}
wz^{-1};p)_\infty}
:S_1(z)\~\Psi^{(\ell)}_{0}(w):,\quad |z|>|q^{k-1} w|,\nn\\
S_{-1}(z)\~\Psi^{(\ell)}_{0}(w)&=&
(q^{-2}z)^{-1/k+2}
{(q^{1+k}p wz^{-1};p)_\infty \over (q^{k-1}p
wz^{-1};p)_\infty}
:S_{-1}(z)\~\Psi^{(\ell)}_{0}(w):,\quad |z|>|q^{k-1}p w|,\nn\\
\ea

\bac
E^+_{1}(z)E^+_1(w)&=&q{z-w\over z-wq^2}:E^+_{1}(z)E^+_1(w):,
\quad |z|>|q^{2} w|, \nn\\
E^+_{1}(z)E^+_{-1}(w)&=&q{z-wq^{-2}\over z-wq^2}
:E^+_{1}(z)E^+_1(w):,\quad |z|>|q^{2} w|,\nn\\
E^+_{-1}(z)E^+_1(w)&=&q^{-1}:E^+_{-1}(z)E^+_1(w):,\nn\\
E^+_{-1}(z)E^+_{-1}(w)&=&q^{-1}{z-w\over z-wq^2}
:E^+_{-1}(z)E^+_{-1}(w):,\quad |z|>|q^{2} w|,\nn\\
E^-_{1}(z)E^-_1(w)&=&q^{-1}{z-w\over z-wq^{-2}}
:E^-_{1}(z)E^-_1(w):,\quad |z|>|w|,\nn\\
E^-_{1}(z)E^-_{-1}(w)&=&q^{-1}{z-wq^2\over z-wq^{-2}}
:E^-_{1}(z)E^-_{-1}(w):,\quad |z|>|q^2 w|,\nn\\
E^-_{-1}(z)E^-_1(w)&=&q:E^-_{-1}(z)E^-_1(w):,\nn\\
E^-_{-1}(z)E^-_{-1}(w)&=&q{z-w\over z-wq^{-2}}
:E^-_{-1}(z)E^-_{-1}(w):,\quad |z|>|w|.
\ea

\bac
S_{1}(z)S_1(w)&=&z^{-{k\over k+2}}q^{{k-2\over k+2}} 
(z-w){(z^{-1}wq^{-2}p;p)_\infty \over 
(z^{-1}wq^{2};p)_\infty},\quad |z|>|w|,\nn\\
S{1}(z)S_{-1}(w)&=&z^{{2\over k+2}}q^{{k-2\over k+2}}
{(z^{-1}wq^{-2};p)_\infty \over 
(z^{-1}wq^{2};p)_\infty},\quad |z|>|q^{-2} w|, \nn\\
S_{-1}(z)S_1(w)&=&z^{{2\over k+2}}q^{-{k+6\over k+2}} 
{(z^{-1}wq^{-2}p;p)_\infty \over 
(z^{-1}wq^{2}p;p)_\infty},\quad |z|>|q^{-2}p w|, \nn\\
S_{-1}(z)S_{-1}(w)&=&
z^{-{k\over k+2}}q^{-{k+6\over k+2}} 
(z-w){(z^{-1}wq^{-2}p;p)_\infty \over 
(z^{-1}wq^{2};p)_\infty},\quad |z|>|w|.
\ea

\bac
E^+_{1}(z)S_{1}(w)&=&q:E^+(z)S_1(w):,\nn\\
E^+_1(z)S_{-1}(w)&=&q{z-q^{-2}w\over z-w}
:E^+_1(z)S_{-1}(w):,\quad |z|>|q^{-2} w|,\nn\\
E^+_{-1}(z)S_1(w)&=&q^{-1}{z-q^{2}w\over z-w}
:E^+_{-1}(z)S_1(w):,\quad |z|>|w|,\nn\\
E^+_{-1}(z)S_{-1}(w)&=&q^{-1}:E^+_{-1}(z)S_{-1}(w):,\nn\\
S_1(z)E^+_1(w)&=&q:S_1(z)E^+_1(w):,\nn\\
S_1(z)E^+_{-1}(w)&=&q{z-q^{-2}w\over z-w}
:S_1(z)E^+_{-1}(w):,\quad |z|>|q^{-2} w|,\nn\\
S_{-1}(z)E^+_1(w)&=&q^{-1}{z-q^{2}w\over z-w}
:S_{-1}(z)E^+_1(w):,\quad |z|>|w|,\nn\\
S_{-1}(z)E^+_{-1}(w)&=&q^{-1}:S_{-1}(z)E^+_{-1}(w):,\nn\\
E^-_1(z)S_1(w)&=&q^{-1}:E^-(z)S_+(w):,\nn\\
E^-_1(z)S_{-1}(w)&=&{q^{-1}z-q^{-k-1}w\over z-q^{-k-2}w}
:E^-_+(z)S_-(w):,\quad |z|>|q^{-k-2} w|,\nn\\
E^-_{-1}(z)S_1(w)&=&{qz-q^{k+1}w\over z-q^{k+2}w}
:E^-_-(z)S_+(w):,\quad |z|>|q^{k} w|,\nn\\
E^-_{-1}(z)S_{-1}(w)&=&q:E^-_{-1}(z)S_{-1}(w):,\nn\\
S_1(z)E^-_1(w)&=&q^{-1}:S_1(z)E^-_1(w):,\nn\\
S_1(z)E^-_{-1}(w)&=&{q^{-1}z-q^{-k-1}w\over z-q^{-k-2}w}
:S_1(z)E^-_{-1}(w):,|z|>|q^{-k-2} w|,\nn\\
S_{-1}(z)E^-_1(w)&=&{qz-q^{k+1}w\over z-q^{k+2}w}
:S_{-1}(z)E^-_1(w):,\quad |z|>|q^{k} w|,\nn\\
S_{-1}(z)E^-_1(w)&=&q:S_{-1}(z)E^-_1(w):.
\ea

\bac
\eta(z)\xi(w)&=&{q^{k+2}\over z-w}:\eta(z)\xi(w):,\quad
|z|>|w|,\nn\\
\eta(z)\~\Phi^{(\ell)}_{\ell}(w)&=&:\eta(z)
\~\Phi^{(\ell)}_{\ell}(w):
,\nn\\
\~\Phi^{(\ell)}_{\ell}(z)\eta(w)&=&
:\~\Phi^{(\ell)}_{\ell}(z)\eta(w):,\nn\\
\eta(z)\~\Psi^{(1)}_{0}(w)&=&q^{-k-2}(z-wq^k):\eta(z)
\~\Psi^{(1)}_{0}(w):
,\quad |z|>|q^{k} w|,\nn\\
\~\Psi^{(1)}_{0}(z)\eta(w)&=&q^{-2}(z-wq^{-k})
:\~\Psi^{(1)}_{0}(z)\eta(w):,\quad |z|>|q^{-k}w|,\nn\\
\eta(z)E^+_{\rho}(w)&=&{q^{k+2}\over
z-wq^{\rho}}:\eta(z)E^+_{\rho}(w):,\quad |z|>|q^{\rho} w|,\nn\\
E^+_{\rho}(z)\eta(w)&=&{q^{k+2}\over
zq^{\rho}-w}:E^+_{\rho}(z)\eta(w):,\quad |z|>|q^{-\rho} w|,\nn\\
\eta(z)E^-_{\epsilon}(w)&=&(
zq^{-k-2}-wq^{(\epsilon-1)(k+1)-1}):\eta(z)E^-_{\epsilon}(w):,
\quad |z|>|q^{\epsilon(k+1)} w|,\nn\\ 
E^-_{\epsilon}(z)\eta(w)&=& 
(zq^{(\epsilon-1)(k+1)-1}-wq^{-k-2})
:E^-_{\epsilon}(z)\eta(w):,\quad 
|z|>|q^{-\epsilon(k+1)} w|,\nn\\ 
\eta(z)S_{\delta}(w)&=&{q^{k+2}\over
z-wq^{\delta}}:\eta(z)S_{\delta}(w):,\quad
|z|>|q^{\delta}w|,\nn\\
S_{\delta}(z)\eta(w)&=&{q^{k+2}\over
zq^{\delta}-w}:\eta(z)S_{\delta}(w):,
\quad |z|>|q^{-\delta}w|,\nn\\
\xi(z)E^+_{\rho}(w)&=&q^{-k-2}(z-wq^{\rho})
:\xi(z)E^+_{\rho}(w):,\quad |z|>|q^{\rho} w|,\nn\\
E^+_{\rho}(z)\xi(w)&=&q^{-k-2}(
zq^{\rho}-w):E^+_{\rho}(z)\xi(w):,\quad |z|>|q^{-\rho}w|,\nn\\
\xi(z)E^-_{\epsilon}(w)&=&{1\over 
zq^{-k-2}-wq^{(\epsilon-1)(k+1)-1}}:\xi(z)E^-_{\epsilon}(w):,
\quad |z|>|q^{\epsilon(k+1)}w|,\nn\\ 
E^-_{\epsilon}(z)\xi(w)&=& 
{1\over zq^{(\epsilon-1)(k+1)-1}-wq^{-k-2}}
:E^-_{\epsilon}(z)\xi(w):,\quad 
|z|>|q^{-\epsilon(k+1)}w|,\nn\\ 
\xi(z)S_{\delta}(w)&=&q^{-k-2}(
z-wq^{\delta}):\xi(z)S_{\delta}(w):,\quad |z|>|q^{\delta}w|,
\nn\\
S_{\delta}(z)\eta(w)&=&q^{-k-2}(
zq^{\delta}-w):S_{\delta}(z)\eta(w):,\quad |z|>|q^{-\delta}w|.
\label{ap10}
\ea

\subsection{\bf Type II vertex operators}

The various c-functions $H$ and  contours $C$ defining the
bosonization of the physical type II vertex operators in
relation (\ref{co8}) are given by
\bac
I_{1}^{(1),(0)}(z,t|C)&=&
-{(q^{k-2}z)^{-{1\over k+2}}\over t(q-q^{-1})}
{(q^{-k+1}ptz^{-1};p)_\infty\over 
(q^{-k-1}ptz^{-1};p)_\infty},\>\>\>{\rm with}\>\>\>
C:\>\>\>|z|>|q^{-k-1}pt|,\nn\\
I_{-1}^{(1),(0)}(z,t|C)&=&
{(q^{k-2}z)^{-{1\over k+2}}\over t(q-q^{-1})}
{(q^{-k+1}tz^{-1};p)_\infty\over 
(q^{-k-1}tz^{-1};p)_\infty},\>\>\>{\rm with}\>\>\>
C:\>\>\>|z|>|q^{-k-1}t|,\nn\\
I_{1}^{(0),(1)}(z,u|C)&=&
{z q^k\over u(u-q^{k-1}z)},\>\>\>{\rm with}\>\>\>
C:\>\>\>|u|>|q^{k-1}z|,\nn\\
I_{-1}^{(0),(1)}(z,u|C)&=&
-{z q^{k+2}\over u(u-q^{k+1}z)},\>\>\>{\rm with}\>\>\>
C:\>\>\>|z|>|q^{-k-1}u|,\nn\\
I_{1,1}^{(1),(1)}(z,t,u|C)&=&
-{zq^{k+1}(q^{k-2}z)^{-{1\over k+2}}\over 
ut(u-q^{k-1}z)(q-q^{-1})}
{(q^{-k+1}tpz^{-1};p)_\infty\over 
(q^{-k-1}tpz^{-1};p)_\infty},\>\>\>{\rm with}\nn\\
C&:&\>\>\> |u|>|q^{k-1}z|,\>\>|z|>|q^{-k-1}pt|,\nn\\
I_{1,-1}^{(1),(1)}(z,t,u|C)&=&
{zq^{k+1}(q^{k-2}z)^{-{1\over k+2}}\over 
ut(u-q^{k+1}z)(q-q^{-1})}{u-q^2t\over u-t}
{(q^{-k+1}tpz^{-1};p)_\infty\over 
(q^{-k-1}tpz^{-1};p)_\infty},\>\>\>{\rm with}\nn\\
C&:&\>\>\> |z|>|q^{-k-1}u|,\>\>|z|>|q^{-k-1}pt|,\>\>
|u|>|t|,\nn\\
I_{-1,1}^{(1),(1)}(z,t,u|C)&=&
{zq^{k+1}(q^{k-2}z)^{-{1\over k+2}}\over 
ut(u-q^{k-1}z)(q-q^{-1})}{u-q^{-2}t\over u-t}
{(q^{-k+1}tz^{-1};p)_\infty\over 
(q^{-k-1}tz^{-1};p)_\infty},\>\>\>{\rm with}\nn\\
C&:&\>\>\> |u|>|q^{k-1}z|,\>\>|z|>|q^{-k-1}t|,\>\>
|u|>|q^{-2}t|,\nn\\
I_{-1,-1}^{(1),(1)}(z,t,u|C)&=&
-{zq^{k+1}(q^{k-2}z)^{-{1\over k+2}}\over 
ut(u-q^{k+1}z)(q-q^{-1})}
{(q^{-k+1}tz^{-1};p)_\infty\over 
(q^{-k-1}tz^{-1};p)_\infty},\>\>\>{\rm with}\nn\\
C&:&\>\>\> |z|>|q^{-k-1}u|,\>\>|z|>|q^{-k-1}t|,
\label{ap2}
\ea

\subsection{\bf Trace of vertex operators over Fock spaces}
Let $\{a_n, n\neq 0\}$ generate a Heisenberg algebra with
\be
[a_n,a_m]=h(n)\delta_{n+m,0}.
\ee
Define as usual  a vacuum vector $|0>$ by
\be
a_n|0>=0,\quad n>0.
\ee
Let $F$ be a Fock space generated from the successive
actions of the creation operators on $|0>$ then we have the
following   formula for the evaluation of a trace of vertex
operators
\bac
&&{\rm Tr}_{F}(\mu^{\sum_{n>0}na_{-n}a_n/h(n)}\exp(-\sum_{n>0}
a_{-n}A_n) \exp(\sum_{n>0}a_{n}A_{-n}))\nn\\
&&={1\over \prod_{n=1}^\infty (1-\mu^n)} \exp(
-\sum_{m=1}^\infty h(m)A_mA_{-m}{\mu^m\over 1-\mu^m}).
\label{ap3}
\ea
Here $A_n$ and $A_{-n}$, $n>0$, are given c-functions.

\subsection{\bf Applications of the trace formula}
Let us consider the general trace in the sector $X_1$ 
\bac
&&T^{X_1}_\ell(N_1,N_2,N_3,N_4,n_4|X,Y,Z,W)
\nn\\
&&= Tr_{\~F^{X_1}_\ell}
(\mu^{-d^{X_1}}\nu^{-a_{X_1,0}}
:\prod_{i=1}^{N_1} e^{X_1(k;2,k+2|q^kx_i;{k+2\over 2})}
\prod_{j=1}^{N_2} e^{X_1(1;2,k+2|
q^{k-2}y_j;-{k+2\over 2})}\nn\\
&&\prod_{r=1}^{N_3} e^{-X_1(k+2|q^{-2}z_r;-{k+2\over 2})}
\prod_{s=1}^{n_4} e^{\partial 
X_1^{+}(q^{-2}w_s;-{k+2\over 2})}
\prod_{s=n_4+1}^{N_4} e^{\partial 
X_1^{-}(q^{-2}w_s;-{k+2\over 2})}:)
\ea
where the sets $X,Y,Z,W$ are defined by
\bac
X&=&\{x_1,\cdots,x_{N_1}\},\nn\\
Y&=&\{y_1,\cdots,y_{N_2}\},\nn\\
Z&=&\{z_1,\cdots,z_{N_3}\},\nn\\
W&=&\{w_1,\cdots,w_{N_4}\},
\ea
Applying the general formula (\ref{ap3}) for the trace we find
\bac
T^{X_1}_\ell(N_1,N_2,N_3,N_4,n_4|X,Y,Z,W)&=&
{\delta_{kN_1+N_2-2N_3}\over
(\mu;\mu)_\infty}
G_1(X)G_2(Y)G_3(Z)G_4(X,Y)G_5(X,Z)\nn\\
&&\times G^{+}_6(X,W)G^{-}_6(X,W)
G_7(Y,Z)G^{+}_8(Y,W)G^{-}_8(Y,W)\nn\\
&&\times G^{+}_9(Z,W)G^{-}_9(Z,W)G_{10}
\label{ap4}
\ea
where
\bac
G_1(X)&=&\prod_{i\leq j}
{(q^{2k+2}p\mu x_ix_j^{-1};q^4;p;\mu)_\infty
(q^{2k+2}p\mu x_jx_i^{-1};q^4;p;\mu)_\infty\over
(q^{2}p\mu x_ix_j^{-1};q^4;p;\mu)_\infty^2
(q^{2}p\mu x_jx_i^{-1};q^4;p;\mu)_\infty^2}\nn\\
&&\times (q^{-2k+2}p\mu x_ix_j^{-1};q^4;p;\mu)_\infty
(q^{-2k+2}p\mu x_jx_i^{-1};q^4;p;\mu)_\infty,\nn\\ 
G_2(Y)&=&\prod_{i\leq j}
{(q^{4}\mu y_iy_j^{-1};q^4;p;\mu)_\infty
(q^{4}\mu y_jy_i^{-1};q^4;p;\mu)_\infty
(\mu y_iy_j^{-1};q^4;p;\mu)_\infty
(\mu y_jy_i^{-1};q^4;p;\mu)_\infty\over
(q^{2}\mu y_iy_j^{-1};q^4;p;\mu)_\infty^2
(q^{2}\mu y_jy_i^{-1};q^4;p;\mu)_\infty^2},\nn\\
G_3(Z)&=&\prod_{i\leq j}
{(q^{-2}\mu z_iz_j^{-1};p;\mu)_\infty
(q^{-2}p\mu z_jz_i^{-1};p;\mu)_\infty\over
(q^{2}\mu z_iz_j^{-1};p;\mu)_\infty
(q^{2}p\mu z_jz_i^{-1};p;\mu)_\infty},\nn\\
G_4(X,Y)&=&\prod_{i, j}
{(q^{3}p\mu x_iy_j^{-1};q^4;p;\mu)_\infty
(q^{-1}p\mu y_jx_i^{-1};q^4;p;\mu)_\infty\over 
(qp\mu x_iy_j^{-1};q^4;p;\mu)_\infty
(q^{-3}p\mu y_jx_i^{-1};q^4;p;\mu)_\infty}\nn\\
&&\times 
{(q^{5}\mu x_iy_j^{-1};q^4;p;\mu)_\infty
(q\mu y_jx_i^{-1};q^4;p;\mu)_\infty\over
(q^{7}\mu x_iy_j^{-1};q^4;p;\mu)_\infty
(q^{3}\mu y_jx_i^{-1};q^4;p;\mu)_\infty},\nn\\
G_5(X,Z)&=&\prod_{i,j}
{(q^{k}\mu p x_iz_j^{-1};p;\mu)_\infty
(q^{k}\mu z_jx_i^{-1};p;\mu)_\infty\over
(q^{-k}\mu p x_iz_j^{-1};p;\mu)_\infty
(q^{-k}\mu z_jx_i^{-1};p;\mu)_\infty},\nn\\
G^{+}_6(X,W)&=&\prod_{i=1}^{n_4}\prod_{j} 
{(q^{2}\mu  x_jw_i^{-1};\mu)_\infty\over 
(q^{2+2k}\mu  x_jw_i^{-1};\mu)_\infty},\nn\\
G^{-}_6(X,W)&=&\prod_{i=n_4+1}^{N_4}\prod_{j} 
{(q^{-2k-2}\mu  x_j^{-1}w_i;\mu)_\infty\over 
(q^{-2}\mu  x_j^{-1}w_i;\mu)_\infty},\nn\\
G(Y,Z)_7&=&\prod_{i,j}
{(q^{1-k}\mu z_iy_j^{-1};p;\mu)_\infty
(q^{1+k}\mu y_jz_i^{-1};p;\mu)_\infty\over
(q^{-k-1}\mu z_iy_j^{-1};p;\mu)_\infty
(q^{k-1}p\mu y_jz_i^{-1};p;\mu)_\infty},\nn\\
G^{+}_8(Y,W)&=&\prod_{i=1}^{n_4}\prod_{j} 
{(q^{-3}\mu  y_jw_i^{-1};\mu)_\infty\over 
(q^{-1}\mu  y_jw_i^{-1};\mu)_\infty},\nn\\
G^{-}_8(Y,W)&=&\prod_{i=n_4+1}^{N_4}\prod_{j} 
{(q^{-2k-3}\mu  y_j^{-1}w_i;\mu)_\infty\over 
(q^{-2k-1}\mu  y_j^{-1}w_i;\mu)_\infty},\nn\\
G^{+}_9(Z,W)&=&\prod_{i=1}^{n_4}\prod_{j} 
{(q^{-k}\mu  z_jw_i^{-1};\mu)_\infty\over 
(q^{-k-4}\mu  z_jw_i^{-1};\mu)_\infty},\nn\\
G^{-}_9(Z,W)&=&\prod_{i=n_4+1}^{N_4}\prod_{j} 
{(q^{-k}\mu  y_j^{-1}w_i;\mu)_\infty\over 
(q^{-k-4}\mu y_j^{-1}w_i;\mu)_\infty},\nn\\
G_{10}&=&\mu^{{\ell(\ell+2)\over 12}}\nu^{-\ell}
q^{(2n_4-N_4)\ell}\prod_{i=1}^{N_1}(q^kx_i)^{k\ell\over 2(k+2)}
\prod_{i=1}^{N_2}(q^{k-2}y_i)^{\ell\over 2(k+2)}
\prod_{i=1}^{N_3}(q^2z_i^{-1})^{\ell\over (k+2)}.
\ea

\bac
&&T^{X_2,X_3}(A,B,C,D,G,H|\~X,\~Y,\~Z,\~W)
\nn\\
&&= {\rm Tr}_{\~F^{X_2,X_3}}{\cal O}={\rm Tr}_{F^{X_2,X_3}}
\left({\cal O} \oint
{dy\over 2 \pi i}\eta(y) \xi(w_0)\right)|_{a_{X_2,0}+
a_{X_3,0}=0}\ea
with
\be
{\cal O}=\mu^{-d^{X_2}-d^{X_3}}\nu^{-a_{X_2,0}}
\prod_{i=1}^{\bar N} e^{-X_3(2|q^{a_i}x_i;0)}
\prod_{j=1}^{\bar N} e^{X_3(2|q^{b_j}y_j;0)}
\prod_{r=1}^{\bar M} e^{-X_2(2|q^{c_r}z_r;d_r)}
\prod_{s=1}^{\bar M} e^{X_2(2|q^{g_s}w_s;h_s)},
\ee
and the various sets of parameters are defined by
\bac
A&=&\{a_1,\cdots,a_{\bar N}\},\nn\\
B&=&\{b_1,\cdots,b_{\bar N}\},\nn\\
C&=&\{c_1,\cdots,c_{\bar M}\},\nn\\
D&=&\{d_1,\cdots,d_{\bar M}\},\nn\\
G&=&\{g_1,\cdots,a_{=M}\},\nn\\
H&=&\{h_1,\cdots,a_{\bar M}\},\nn\\
\~X&=&\{x_1,\cdots,x_{\bar N}\},\nn\\
\~Y&=&\{y_1,\cdots,y_{\bar N}\},\nn\\
\~Z&=&\{z_1,\cdots,z_{\bar M}\},\nn\\
\~W&=&\{w_1,\cdots,w_{\bar M}\},
\ea
Applying (\ref{ap3}) we find
\bac
&&T^{X_2,X_3}(A,B,C,D,G,H|\~X,\~Y,\~Z,\~W)\nn\\
&&={\delta_{\bar N,\bar M}\gamma\over
\prod_{n=1}(1-\mu^n)^2}G_1(\~X,w_0)
G_2(\~Y,\gamma)G_3(\~Z)G_4(\~W)G_5(\~X,\~Y,w_0,\nu)
G_6(\~Z,\~W)G_7(\~X,\~Y,w_0,\gamma)\nn\\
\label{ap5}
\ea
where the delta function $\delta_{\bar N,\bar M}$ ensures that the
above trace
is restricted to the space $\~F^{X_2,X_3}$, and
\bac
G_1(\~X,w_0)&=&\prod_{i,i^\prime=0}^{\bar N}(\mu
q^{a_i-a_{i^\prime}}x_ix_{i^\prime}^{-1};\mu)_\infty,\quad 
{\rm with}\>\>\> x_0=w_0,\quad a_0=-k-2,\nn\\
G_2(\~Y,\gamma)&=&\prod_{j,j^\prime=0}^{\bar N}(\mu
q^{b_j-b_{j^\prime}}y_jy_{j^\prime}^{-1};\mu)_\infty,\quad 
{\rm with}\>\>\> y_0=\gamma,\quad b_0=-k-2,\nn\\
G_3(\~Z)&=&\prod_{r,r^\prime=1}^{\bar M}(\mu
q^{c_r+d_r+d_{r^\prime}-c_{r^\prime}}
z_rz_{r^\prime}^{-1};\mu)^{-1}_\infty,\nn\\
G_4(\~W)&=&\prod_{s,s^\prime=1}^{\bar M}(\mu
q^{g_s+h_s+h_{s^\prime}-g_{s^\prime}}
w_sw_{s^\prime}^{-1};\mu)^{-1}_\infty,\nn\\
G_5(\~X,\~Y,w_0,\gamma)&=&\prod_{i,j=0}^{\bar N}(\mu
q^{a_i-b_{j}}x_iy_{j}^{-1};\mu)_\infty^{-1}
(\mu q^{b_j-a_{i}}y_jx_{i}^{-1};\mu)^{-1}_\infty,\nn\\
G_6(\~Z,\~W)&=&\prod_{r,s=1}^{\bar M}(\mu
q^{c_r+d_r+h_{s}-g_{s}}
z_rw_{s}^{-1};\mu)_\infty
(\mu q^{g_s+h_s+d_{r}-c_{r}}
w_sz_{r}^{-1};\mu)_\infty,\nn\\
G_7(\~X,\~Y,w_0,\gamma)&=&{q^{k+2}\over \gamma-w_0}
\prod_{i=1}^{\bar N}{q^{a_i}x_i-q^{-k-2}w_0\over 
q^{a_i}x_i-q^{-k-2}\gamma} 
\prod_{j=1}^{\bar M}{q^{b_j}y_j-q^{-k-2}\gamma\over 
q^{b_j}y_j-q^{-k-2}w_0},\nn\\
\gamma&=&w_0 \nu^{-2}\prod_{i=1}^{\bar N} q^{a_i-b_i}x_iy_i^{-1}
 \prod_{r=1}^{\bar M} q^{g_r-c_r}w_rz_r^{-1}.
\ea

\subsection{\bf Useful functions and integrals}

Eta function:
\be
\eta(z)=z^{1\over 24}\prod_{n=1}^\infty (1-z^n).
\label{ap7}
\ee
Jacobi theta function: 
\be
\theta_p(z)=(z;p)_\infty (pz^{-1};p)_\infty
(p;p)_\infty,\quad {\rm with}\quad 
(z;p)_\infty=\prod_{n=0}^\infty (1-zp^n).
\label{ap8}
\ee
Jackson integrals:
\bac
\int_0^cd_ptf(t)&=&c(1-p)\sum_{n=0}^\infty f(cp^n)p^n,\nn\\
\int_0^{c\infty}d_ptf(t)&=&c(1-p)
\sum_{n=-\infty}^\infty f(cp^n)p^n.
\label{ap12}\ea
q-gamma function:
\be
\Gamma_p(z)={(p;p)_\infty \over (p^z;p)_\infty}
(1-p)^{1-z}.
\ee
q-beta function:
\be
B_p(x,y)=\int_0^1 d_pt t^{x-1} 
{(pt;p)_\infty \over (p^y;p)_\infty}.
\label{ap11}
\ee
Basic hypergeometric series:
\be
F_p(a,b,c;z)={\Gamma_p(c)\over \Gamma_p(a) \Gamma_p(c-a)}
\int_0^1 d_pt t^{a-1}{(pt;p)_\infty (p^bzt;p)_\infty \over
(p^{c-a}t;p)_\infty (zt;p)_\infty}.
\ee

\subsection{\bf Useful relations}
\be
B_p(x,y)={\Gamma_p(x)\Gamma_p(y)\over \Gamma_p(x+y)}.
\ee
\be
\Gamma_p(z+1)={1-p^z\over 1-p}\Gamma_p(z).
\ee
\bac
\int_0^c d_pt f(t)&=&c\int_0^1 d_pt f(ct),\nn\\
\int_{cp}^{c\infty} d_pt f(t)&=&c\int_p^\infty d_pt f(ct),\nn\\
\int_p^\infty d_pt f(t)&=&\int_0^1 d_pt t^{-2} f(t^{-1}).
\ea
Fourier double transformation formula:
\be
\sum_{n\in\Z}\oint_{C_0}{dz\over 2 \pi i}z^k f(z)=f(1).
\label{ap6}
\ee
Triple Jacobi identity
\be
\sum_{n\in\Z}x^{n^2/2}y^{2n}=(x;x)_\infty 
(-x^{1\over 2} y^2;x)_\infty
(-x^{1\over 2} y^{-2};x)_\infty.
\ee

\subsection{\bf Face type Boltzmann weights}

The face type Boltzmann weights
$W\left.\left(\matrix{\lambda&\mu\cr
\mu^\prime&\nu\cr}\right|z\right)$ that appear in the commutation
relations of the vertex operators (\ref{vo7}) are given by
\bac
W\left. \left(\matrix{\lambda&\mu\cr
\mu^\prime&\nu\cr}\right| z\right)&=&
-z^{\Delta_\lambda+
\Delta_\nu-\Delta_\mu-\Delta_\mu^\prime-1/2}
{\xi(z^{-1};1,pq^4)\over \xi(z;1,pq^4)}
\times \hat W\left. \left(\matrix{\lambda&\mu\cr
\mu^\prime&\nu\cr}\right| z\right),\nn\\
\hat W\left. \left(\matrix{\lambda&\lambda_+\cr
\lambda_+&\lambda\cr}\right| z\right)&=&
{\theta_p(pq^2)\over \theta_p(pq^{-2n-2})}
{\theta_p(pq^{-2n-2}z)\over \theta_p(pq^{2}z)},\nn\\
\hat W\left.\left(\matrix{\lambda&\lambda_+\cr
\lambda_-&\lambda\cr}\right| z\right)&=&
q^{-1}
\times
{\Gamma_p((2n+2)s)\Gamma_p((2n+2)s)\over 
\Gamma_p((2n+4)s)\Gamma_p(2ns)}
{\theta_p(pz)\over \theta_p(pq^{2}z)},\nn\\
\hat W\left. \left(\matrix{\lambda&\lambda_-\cr
\lambda_+&\lambda\cr}\right| z\right)&=&
q^{-1}
\times
{\Gamma_p(1-(2n+2)s)\Gamma_p(1-(2n+2)s)\over 
\Gamma_p(1-(2n+4)s)\Gamma_p(1-2ns)}
{\theta_p(pz)\over \theta_p(pq^{2}z)},\nn\\
\hat W\left.\left(\matrix{\lambda&\lambda_-\cr
\lambda_-&\lambda\cr}\right| z\right)&=&
z^{-1}\times 
{\theta_p(pq^2)\over \theta_p(q^{2n+2})}
{\theta_p(q^{2n+2}z)\over \theta_p(pq^2z)} ,\nn\\
\hat W\left.\left(\matrix{\lambda&\lambda_\pm\cr
\lambda_\pm&\lambda_\pm\cr}\right| z\right)&=&1,\nn\\
\hat W\left.\left(\matrix{\lambda&\lambda_\pm\cr
\lambda_\pm&\lambda_\pm\cr}\right| z\right)&=&0 
\quad {\rm otherwise}.
\ea

\pagebreak


\begin{thebibliography}{10}


\bibitem{Bet31} H. Bethe, {\em Z. Phys.} {\bf 71}, 205 (1931).

\bibitem{Hei28} W. Heisenberg, {\em Z. Phys.} {\bf 49}, 619 (1928).

\bibitem{Hul38} L. Hulth\'en, {\em Arkiv Mat. Astron. Fysik A11} 
{\bf 26}, 1 (1938).

\bibitem{LiMa62} E.H. Lieb and D.C. Mattis, {\em J. Math. Phys.} 
{\bf 3}, 749 (1962).

\bibitem{ClPe62} J. des Cloizeaux and J.J. Pearson, {\em Phys. Rev.}
          {\bf 128}, 2131 (1962).

\bibitem{Gri64} R.B. Griffiths, {\em Phys. Rev. } {\bf 133}, A768 (1964). 

\bibitem{YaYa66} C.N. Yang and C.P. Yang, {\em Phys. Rev. } 
{\bf 150}, 321 (1966);
{\bf 150}, 327 (1966); {\bf 151}, 258 (1966).

\bibitem{Nie67} Th. Niemeijer, {\em Physica} {\bf 36}, 377 (1967).

\bibitem{LuPe69} A. Luther and I. Peschel, {\em Phys. Rev.} 
{\bf B12}, 2131 (1969).

\bibitem{Baral71} E. Barouch, B.M. McCoy and D.B. Abraham, 
{\em Phys. Rev.} {\bf A4}, 2331 (1971).

\bibitem{Gau71} M. Gaudin, {\em Phys. Rev. Lett.} 
{\bf 26}, 1301 (1971).

\bibitem{Tak71} M. Takahashi, {\em Prog. Theor. Phys.} {\bf 46}, 
401 (1971).

\bibitem{JoMc72}J.D. Johnson and B.M. McCoy,
{\em Phys. Rev.}, {\bf A6}, 1613 (1972).

\bibitem{Bax72}R. J. Baxter, {\em Ann. Phys.} 
{\bf 70}, 323 (1972).

\bibitem{Bax73}R. J. Baxter, {\em J. Stat. Phys.} 
{\bf 9}, 145 (1973).


\bibitem{FaTa79} L.D. Faddeev and 
L.A. Takhtajan, {\em Russ. Math. Surveys} {\bf 34}, 11 (1979).

\bibitem{Mulal81} G. M\"uller, H. Thomas, H. Beck, and J.C. 
Bonner, {\em  Phys. Rev.}
          {\bf B24}, 1429 (1981).

\bibitem{Bax82}
R.~J. Baxter.
\newblock {\em Exactly Solved Models in Statistical Mechanics}.
\newblock Academic, London, 1982.

\bibitem{McC83} B.M. McCoy, J.H.H. Perk and R.E. Shrock,
 {\em Nucl. Phys.}
         {\bf B220}, 35 (1983).

\bibitem{Babal83} O. Babelon, H.J. de Vega and C.M. Viallet, 
{\em Nuc. Phys.}{\bf B220}, 13 (1983).

\bibitem{Sog83}K. Sogo,
{\em Phys. Lett.} {\bf A104}, 51 (1983).

\bibitem{MuSh84} G. M\"uller and R. E. Shrock,
           {\em  Phys. Rev.} {\bf B29}, 288 (1984).

\bibitem{Rolal86} J.M.R. Roldan, B.M. McCoy and J.H.H. Perk, 
{\em Physica }
           {\bf 136A}, 255 (1986).


\bibitem{Koral93} V.E. Korepin, A.G. Izergin, and N.M. 
Bogoliubov,{\em The Quantum
           Inverse Scattering Method and Correlation Functions},
           Cambridge University Press, (1993).




\bibitem{Daval92}
B.~Davies, O.~Foda, M.~Jimbo, T.~Miwa, and A.~Nakayashiki,
\newblock {\em Comm. Math. Phys.}, 
{\bf 151}, 89 (1993).

\bibitem{Idzal93} M. Idzumi, T. Tokihiro, K. Iohara, M. Jimbo,
T. Miwa and T. Nakashima, {\it Int. J. Mod. Phys.} 
{\bf A8}, 1479 (1993).


\bibitem{JiMi94}
M. Jimbo and T. Miwa,
\newblock {\it Algebraic Analysis of Solvable Lattice Models},
\newblock{ American Mathematical Society}, (1994).







\bibitem{Koral94} V.E. Korepin, A.G. Izergin, F.H.L. 
Essler and D. 
Uglov, {\it Phys. Lett.} {\bf 190A}, 182 (1994).

\bibitem{Fleal95} A. Fledderjohann, M. Karbach, K.-H. M\"utter, 
and P. Wielath, {\it J. Phys.: Cond. Matter} {\bf 7} 8993 (1995).


\bibitem{Jin92}N.H. Jing, {\it On a trace of $q$ analog 
vertex operators}, in Quantum Groups, Spring Workshop on
Quantum Groups, World Scientific, Singapore (1992).

\bibitem{Idz93}M. Idzumi, {\it Calculation of 
Correlation Functions of the Spin-1 XXZ model by Vertex
Operators}, Ph.D. thesis, 1993. 

\bibitem{BeFe90}D. Bernard and G. Felder,
{\em Commu. Math. Phys.}, {\bf 127}, 145 (1990).  

\bibitem{Kon94}H. Konno, 
{\em Nucl. Phys.}, {\bf B432}, 457 (1994);
{\rm Mod. Phys. Lett.}, {\bf A}.

\bibitem{BoWe94}
A.H. Bougourzi and R.A. Weston,
\newblock {\em Nucl. Phys.}, {\bf B417}, 439 (1994).

\bibitem{Smi92} F.A. Smirnov,
\newblock {\it Form Factors in Completely Integrable Models
of Quantum Field Theory}, World Scientific, Singapore, (1992).

\bibitem{Dri85}
V.~G. Drinfeld,
\newblock {\em Soviet Math. Doklady}, {\bf 32}, 254 (1985).

\bibitem{Jim85}
M.~Jimbo,
\newblock {\em Lett. Math. Phys.}, {\bf 10}, 63 (1985).

\bibitem{Dri86}
V.~G. Drinfeld,
\newblock Proc. ICM, Am. Math. Soc., Berkeley, CA, (1986).

\bibitem{ChPr91}
V. Chari and A. Pressley.
\newblock {\em Comm. Math. Phys.}, {\bf 142}, 261 (1991).


\bibitem{Jimal92}
M.~Jimbo, K.~Miki, T.~Miwa, and A.~Nakayashiki,
\newblock {\em Phys. Lett.}, {\bf A168}, 256 (1992).

\bibitem{FrJi88}
I.B. Frenkel and N. Jing,
\newblock {\em Proc. Natl. Acad. Sci.}, {\bf 85}, 9373 (1988).

\bibitem{Jin96} N. Jing, 
\newblock {\em J. Alg.}, {\bf 182}, 448 (1996).




\bibitem{BoVi96}A.H. Bougourzi and L. Vinet,
{\em J. Math. Phys.}, {\bf  37}, 3548 (1996).

\bibitem{BoVi94}A.H. Bougourzi and L. Vinet,
{\em Lett. Math. Phys.}, {\bf 36}, 101 (1996).   

\bibitem{Bou93}A.H. Bougourzi,
{\em Nuc. Phys.}, {\bf B404}, 457 (1993). 

\bibitem{Boual93}
A. Abada, A.H. Bougourzi and M.A. El Gradechi,
{\em Mod. Phys. Lett.},  {\bf A8}, 715 (1993).

\bibitem{Mat92}
A. Matsuo, {\em Phys. Lett.}, {\bf B308}, 260 (1992);
{\em Commu. Math. Phys.} {\bf 160} 33 (1994).


\bibitem{Shi92}J. Shiraishi, {\em Phys. Lett.}, 
{\bf A171}, 243 (1992).\newline
A. Kato, Y-H. Quano and J. Shiraishi, 
{\em Comm. Math. Phys.}, {\bf 157}, 119 (1993).





\bibitem{BoWe95}A.H. Bougourzi and R. Weston,
{\em Int. J. Mod. Phys.}, {\bf A10}, 561 (1995).    


\bibitem{Fabal96}K. Fabricius, U. Low, and J. Stolze,
{\it Dynamic correlations of antiferromagnetic spin-1/2 XXZ chains at arbitrary
     temperature from complete diagonalization}, 
preprint cond-mat/9611077, 1996. 
   
\bibitem{Pak93}S. Pakuliak, 
{\it Annihilation Poles for Form Factors in $XXZ$ Model},
preprint RIMS-933, hep-th/9307090.  

\end{thebibliography}
\end{document}